\documentclass[opre,nonblindrev]{informs3Itai}

\OneAndAHalfSpacedXI

\usepackage{natbib}
 \bibpunct[, ]{(}{)}{,}{a}{}{,}%
 \def\bibfont{\small}%
 \def\newblock{\ }%
\usepackage{enumerate}
\usepackage{latexsym, bbm, amssymb, amsmath, mathtools, mathabx}
\usepackage{multirow}
\usepackage{graphicx}
\usepackage{float}
\usepackage{paralist}
\usepackage{pgfplots} 
\usepackage{ifpdf}
\ifpdf
    \usepackage[pdftex]{hyperref}
\else
    \usepackage{hyperref}
\fi
\usepackage{tikz}
\usetikzlibrary{patterns,snakes,shapes}
\usepackage[margin=1in]{geometry}

\newcommand{\frakS}{\mathfrak{S}}
\newcommand{\T}{\mathcal{T}}

\newcommand{\R}{\mathbb{R}}

\newcommand{\bbZ}{\mathbb{Z}}
\newcommand{\Z}{\mathbb{Z}}

\renewcommand{\P}{\mathbb{P}}
\newcommand{\A}{\mathcal{A}}
\newcommand{\E}{\mathbb{E}}
\newcommand{\F}{\mathcal{F}}
\newcommand{\Fbar}{\bar{F}}
\newcommand{\Mbar}{\bar M}

\newcommand{\1}{\mathbbm{1}}

\DeclareMathOperator{\Var}{Var}

\DeclarePairedDelimiter\abs{\lvert\,}{\,\rvert}%

\newcommand{\Von}{V_{\rm on}}
\newcommand{\Voff}{V_{\rm off}}
\newcommand{\Vna}{V_{\rm na}}
\newcommand{\na}{{\rm na}}
\newcommand{\id}{{\rm id}}

\newcommand{\p}{\mathfrak{p}}

\newcommand{\sigmaBR}{\sigma^{\rm br}}
\newcommand{\BR}{{\rm br}}

{\theoremstyle{EXkey}\newtheorem{notation}{Notation}}

\TheoremsNumberedThrough     
\ECRepeatTheorems

\EquationsNumberedThrough    

\let\oldmarginpar\marginpar
\renewcommand\marginpar[1]{\-\oldmarginpar[\raggedleft\scriptsize #1]%
{\raggedright\scriptsize #1}}

\begin{document}

\TITLE{Uniformly Bounded Regret \\ in the Multi-Secretary Problem}

\RUNTITLE{Uniformly bounded regret in the multi-secretary problem}

\ARTICLEAUTHORS{%
\AUTHOR{Alessandro Arlotto}
\AFF{The Fuqua School of Business, Duke University, 100 Fuqua Drive, Durham, NC, 27708, \\  \EMAIL{alessandro.arlotto@duke.edu}}

\AUTHOR{Itai Gurvich}
\AFF{Cornell School of ORIE and Cornell Tech, 2 West Loop Road, New York, NY, 10044. \\ \EMAIL{gurvich@cornell.edu}}
}

\RUNAUTHOR{A. Arlotto and I. Gurvich}

\ABSTRACT{In the secretary problem of \citet{Cay:ET1875} and \citet{Mos:SM1956},
$n$ non-negative, independent, random variables with common distribution
are sequentially presented to a decision maker who decides when to stop and collect the most recent realization.
The goal is to maximize the expected value of the collected element.
In the $k$-choice variant, the decision maker is allowed to make $k \leq n$ selections
to maximize the expected total value of the selected elements.
Assuming that the values are drawn from a known distribution with finite support,
we prove that the best regret---the expected gap between the optimal online policy
and its offline counterpart in which all $n$ values are made visible at time $0$---is uniformly bounded in
the number of candidates $n$ and the budget $k$.
Our proof is constructive: we develop an adaptive Budget-Ratio policy that achieves this performance.
The policy selects or skips values depending on where the ratio of the residual budget to the remaining time
stands relative to multiple thresholds that correspond to middle points of the distribution.
We also prove that being adaptive is crucial: in general, the minimal regret among non-adaptive policies grows like the square root of $n$.
The difference is the value of adaptiveness.
}
	
\KEYWORDS{multi-secretary problem, regret, adaptive online policy.}

\MSCCLASS{Primary:    90C39. 
    Secondary:  60C05, 
                60J05, 
                68W27, 
                68W40, 
                90C27, 
                90C40. 
}

\ORMSCLASS{Primary:
    Decision Analysis: Sequential, Theory;
    Dynamic programming: Markov.
Secondary:
    Probability: Markov processes;
    Analysis of Algorithms: Suboptimal algorithms.
}

\HISTORY{\emph{First version:} October 20, 2017. \emph{This version:} June 1, 2018. 
}

\maketitle
\vspace{-0.5cm}

\section{Introduction}

In the classic formulation of the secretary problem  a decision maker (referred to as ``she'')
is sequentially presented with $n$ non-negative, independent, values representing the ability of potential candidates
and must select one candidate (referred to as ``he'').
Every time a new candidate is inspected and his ability is revealed,
the decision maker must decide whether to reject or select the candidate, and her decision is irrevocable.
If the candidate is selected, then the problem ends;
if the candidate is rejected, then he cannot be recalled at a later time.
The decision maker knows the number of candidates $n$,
the distribution of the ability values in the population,
and her objective is to maximize the probability of selecting the most able candidate.
For any given $n$ the problem can be solved by dynamic programming, but
there is an asymptotically optimal heuristic that is remarkably elegant.
The decision maker observes the abilities of the first $n/e$ candidates
and selects the first candidate whose ability exceeds that of the current best candidate;
or the last candidate if no such candidate exists
\citep[see, e.g.,][]{Lindley:AS1961,ChowMorigutiRobbinsSamuels:IJM1964,GilbertMosteller:JASA1966,louchard2016sequential}

Several variations of this simple model have been introduced in the literature,
and we refer to \citet{Freeman:ISR1983} and \citet{Ferguson:SS1989} for a survey of extensions and references.
Relevant to us is the formulation in which the decision maker
seeks to maximize the expected ability of the selected candidate, rather
than maximizing the probability of selecting the best.
This problem was first considered by \citet{Cay:ET1875} and \citet{Mos:SM1956},
and it is a special case of the $k$-choice (multi-secretary) problem we study here.
In our formulation:
\begin{itemize}
    \item candidate abilities are independent, identically distributed, and supported on a finite set;
    \item the decision maker is allowed to select up to $k$ candidates ($k$ is the \emph{recruiting budget}); and
    \item the decision maker's goal is to maximize the expected total ability of the selected candidates.
\end{itemize}
This multi-secretary problem has applications in revenue management and auctions, among others. In the standard capacity-allocation revenue management problem
a firm sells $k$ items (e.g. airplane seats) to $n$ customers from a discrete set of fare classes over a finite horizon
and wishes to allocate the seats in the best possible way
\citep[see, e.g.,][]{kleywegt1998dynamic,TallvanR:KLU2004}.
In auctions, the decision maker observes arriving bids and must decide whether to allocate one of the available $k$ items to an arriving customer
\citep[see, e.g.,][]{kleinberg2005multiple,Babaioff2007knapsack}.

The performance of any online algorithm for this $k$-choice secretary problem
is bounded above by the offline (or posterior) sort:
the decision maker waits until all the ability values are presented, sorts them, and then picks the largest $k$ values.
As such, we define the regret of an online selection algorithm as the expected gap in performance
between the online decision and the offline solution,
and we prove that the optimal online algorithm has a regret that is bounded uniformly over $(n,k)$.
The constant bound depends only on the maximal element of the support
and on the minimal mass value of the ability distribution.

Our proof is constructive: we devise an adaptive policy---the Budget-Ratio (BR) policy---that achieves bounded regret.
The policy is adaptive in the sense that the actions are adjusted based on the remaining number of candidates to be inspected
and on the residual budget.
The proof that the BR policy achieves bounded regret is based on an interesting drift argument
for the BR process that tracks the ratio between the residual budget and the remaining number of candidates to be inspected.
Under our policy, BR sample paths are attracted to and then remain pegged to a certain value;
see in Figure \ref{fig:intro} in Section \ref{sec:heuristic}.
Drift arguments are typical in the study of stability of queues \citep[see, e.g.,][]{bramson2008stability},
but the proof here is somewhat subtle:
the drift is strong early in the horizon but it weakens as the remaining number of steps decreases.
Since ``wrong'' decisions early in the horizon are more detrimental,
this diminishing strength does not compromise the regret.
This result shows, in particular, that the budget ratio is a sufficient state descriptor for the development of near-optimal policies.
In other words, while there are pairs of (\emph{residual budget}, \emph{remaining number of candidates})
that share the same ratio, the same BR decisions, but different optimal actions,
these different decisions have little effect on the total expected reward.

We also show that {\em adaptivity is key}. While non-adaptive policies could have temporal variation in actions
(and hence could have different actions towards the end of the horizon) this variation is not enough:
the regret of non-adaptive policies is, in general, of the order of $\sqrt{n}$.
Specifically, non-adaptive policies tend to be too greedy and run out of budget too soon.
Non-adaptive policies introduce independence between decisions and, consequently,
``too much'' variance into the speed at which the recruiting budget is consumed.

Closely related to our work is the paper \cite{wu2015algorithms} that offers an elegant adaptive index policy that we re-visit in Section \ref{se:adaptive-index}.
Under this policy, the ratio of remaining budget to remaining number of steps is a martingale. When the initial ratio of budget to horizon, $k/n$,
is safely far from the jumps of the discrete distribution, this martingale property guarantees a bounded regret.
In general, however, it is precisely the martingale symmetry that increases the regret.
In their proof, \cite{wu2015algorithms} show that their policy achieves, up to a constant, a deterministic-relaxation upper bound.
This upper bound is not generally achievable and we must, instead, use the stochastic {\em offline} sort as a benchmark.
The detailed analysis of the {\em offline} sort gives rise to sufficient conditions that, when satisfied by an {\em online} policy,
guarantee uniformly bounded regret.

\begin{notation}[Notation.]
We use $\Z_+$ to denote the non-negative integers and $\R_+$ to denote the non-negative reals.
For $j \in \{1,2,\ldots\}$, we use $[j]$ to denote the set of integers $\{1,\ldots,j\}$
and we set $[j] = \emptyset$ otherwise.
Given the real numbers $x,y,z$, we set $(x)_+=\max\{0,x\}$, $x \wedge y = \min\{x,y\}$,
and we write $y = x  \pm z$ to mean that $\abs{ y - x} \leq z$.
Throughout, to simplify notation, we use $M \equiv M(x,y,z)$ to denote a Hardy-style constant dependent on $x,y$, and $z$
that may change from one line to the next.
\end{notation}

\section{The multi-secretary problem\label{sec:model}}

A decision maker is sequentially presented with $n$ candidates with abilities $X_1, X_2, \ldots, X_n$,
and, given a recruiting budget equal to $k$,
can select up to $k$ candidates to maximize the total expected ability
of those selected up to and including time $n$.
Of course, there is nothing to study if $k>n$:
the decision maker can take all candidates,
so it suffices to consider pairs $(n,k)$ that belong to the lattice triangle
$$
\T = \{(n,k) \in \Z_+^2: 0 \leq k \leq n\}.
$$

The abilities $X_1, X_2, \ldots, X_n$ are assumed to be independent across candidates and
drawn from a common cumulative distribution function $F$ supported on a finite set
$\A = \{a_m,a_{m-1}, \ldots , a_1: 0 < a_m < a_{m-1}<  \cdots < a_1\}$ of distinct real numbers.
We denote by $(f_m, f_{m-1}, \ldots, f_1)$ the probability mass function with $f_j = \P(X_1 = a_j)$ for all $j \in [m]$,
and we let $F$ be the cumulative distribution function
(so that $f_m=F(a_m)<F(a_{m-1})<\ldots < F(a_1)=1$)
and $\Fbar=1-F$ be the survival function.
Also, for future reference we choose a value $a_{m+1} < a_m$ with $f_{m+1} = 0$ so that $F(a_{m+1})=0$ and $\Fbar(a_{m+1}) = 1$.

The selection process unfolds as follows: suppose that at time $t \in [n]$
the residual recruiting budget is $\kappa$
and that the sum of the abilities of the selected candidates up to and including time $t-1$ is $w$.
If the candidate inspected (or ``interviewed'') at time $t$ has ability $X_t = a$,
then the decision maker may \emph{select} the candidate---increasing the cumulative ability to $w + a$ and reducing the residual budget to $\kappa-1$---or to \emph{reject} the candidate---leaving the accrued ability at $w$ and the remaining budget at $\kappa$.

A policy is \emph{feasible} if the number of selected candidates does not exceed the recruiting budget $k$.
It is \emph{online} if the decision with regards to the $t$th candidate is based only on its ability,
the abilities of prior candidates, and the history of the (possibly randomized) decisions up to time $t$.
All decisions are final: if the candidate interviewed at time $t$ is rejected, it is forever lost.
Vice versa, if the $t$th candidate is selected at time $t$, then that decision cannot be revoked at a later time.

Given a sequence $U_1, U_2, \ldots, U_n$ of independent random variables with the uniform distribution on $[0,1]$
that is also independent of $\{X_1, X_2, \ldots, X_n\}$,
we let $\F_0$ denote the trivial $\sigma$-field and,
for $t \in [n]$, we set $\F_t = \sigma\{X_1, U_1, X_2, U_2, \ldots, X_t, U_t\}$ be the $\sigma$-field
generated by the random variables $\{X_1, U_1, X_2, U_2, \ldots, X_t, U_t\}$. An \emph{online policy} $\pi$ is a sequence
of $\{\F_t: t \in [n]\}$-adapted binary random variables
$\sigma^\pi_1, \sigma^\pi_2, \ldots, \sigma^\pi_n$
where $\sigma^\pi_t=1$ means that the candidate with ability $X_t$
is selected. Since the uniform random variables $U_1, U_2, \ldots, U_t$ are included in $\F_t$,
then the time-$t$ decision, $\sigma^\pi_t$, can be randomized. A \emph{feasible online} policy requires that the number of selected candidates does not
exceed the recruiting budget, i.e., that $\sum_{t\in[n]}\sigma^\pi_t \leq k$, so
$$
\Pi(n,k)\equiv \left\{ (\sigma^\pi_1, \sigma^\pi_2, \ldots, \sigma^\pi_n) \in \{0,1\}^{n}: \sigma^\pi_t \in \F_t \text{ for all } t \in [n] \text{ and } \sum_{t\in [n]}\sigma^\pi_t \leq k\right\},
$$
is then the set of all feasible online policies.
For $\pi \in \Pi(n,k)$, we let $W^\pi_0, W^\pi_1, \ldots, W^\pi_n$ be the sequence of random
variables that track the accumulated ability: we set $W^\pi_0 = 0$, and for $r \in [n]$ we let
$$
W^\pi_r = \sum_{t\in [r]} X_t \sigma^\pi_t = W^\pi_{r-1} + X_r \sigma ^\pi_r.
$$
The expected ability accrued at time $n$ by policy $\pi$ is then given by
$$
\Von^{\pi}(n,k) = \E[W^\pi_n] = \E \left[\sum_{t\in [n]} X_t \sigma_t^{\pi}\right].
$$
For each $(n,k)\in\mathcal{T}$, the goal of the decision maker is to maximize the expected value:
\begin{equation*}
\Von^*(n,k)= \max_{\pi\in\Pi(n,k)}\Von^{\pi}(n,k).
\end{equation*}
For completeness, we include in Appendix \ref{sec:bellman} an analysis of the dynamic program.
The analysis of the Bellman equation confirms an intuitive property of ``good'' solutions: the optimal action should only depend on the remaining number of steps and the remaining budget
and not on the current level of accrued ability.
It also allows for a comparison of the Budget-Ratio policy we propose in Section \ref{sec:heuristic}
with the optimal policy.

{\em Non-adaptive} policies are an interesting subset of feasible online policies.
If the residual recruiting budget at time $t$ is positive,
then a non-adaptive policy selects the arriving candidate with ability $X_t = a_j$
with probability $p_{j,t} \in [0,1]$ independently of all the previous actions.
The probabilities $p_{j,t}$ are determined in advance: they can vary from one period to
the next, but this variation is not adapted in response
to the previous selection/rejection decisions.
A more complete description of non-adaptive policies appears in Section \ref{sec:nonadaptive}.
We let $\Pi_{\na} \subseteq \Pi(n,k)$ be the family of non-adaptive policies and define
$$
\Vna^{*}(n,k)=\sup_{\pi \in \Pi_{\na}} \Von^{\pi}(n,k)
$$
to be the optimal performance among these.

No feasible online policy---adapted or not---can do better than the offline, full-information,
counterpart in which all values are presented in advance. The expected ability accrued at time $n$ by the offline problem is given by
\begin{equation*}
  \Voff^*(n,k)
  = \E\left[ \max_{\sigma_1, \ldots, \sigma_n} \left\{ \sum_{t\in[n]} X_t \sigma_t: (\sigma_1, \ldots, \sigma_n) \in \{0,1\}^n \text{ and } \sum_{t\in[n]} \sigma_t \leq k\right\} \right].
\end{equation*}
Since the offline solution selects the best $k$ candidates {\em for each} realization of $X_1, \ldots, X_n$,
we have that
$$
\Von^{\pi}(n,k)\leq \Voff^*(n,k), \mbox{ for all } (n,k)\in\mathcal{T} \text{ and all } \pi \in \Pi(n,k).
$$
We use the value of the offline problem as a benchmark.
The optimal online policy trivially achieves the offline benchmark if $k=0$ or $k = n$.
The main results of this paper (Theorem \ref{thm:nonadaptive} and Corollary \ref{cor:regret}) are gathered in the following theorem.

\begin{theorem}[The regret of online feasible policies]\label{thm:main}
Let $\epsilon = \tfrac{1}{2} \min\{ f_m, f_{m-1}, \ldots, f_1\}$.
Then there exists a policy $\BR \in \Pi(n,k)$ and a constant $M \equiv M(\epsilon)$ such that
  \begin{equation*}
    \Voff^{*}(n,k) - \Von^{*}(n,k) \leq \Voff^{*}(n,k) - \Von^{\BR}(n,k) \leq a_1 M
    \quad  \text{ for all } (n,k) \in \T.
  \end{equation*}
Furthermore, if $(f_1 + \epsilon) n  \leq k \leq (1-f_m - \epsilon) n$,
then an optimal non-adaptive policy has regret that grows like $\sqrt{n}$.
Specifically, there is a constant $M \equiv M(\epsilon, m, a_1, \ldots, a_m)$
such that
$$
M \sqrt{n }\leq \Voff^{*}(n,k)-\Vna^{*}(n,k).
$$
\end{theorem}

We conclude this section with a discussion of the dependence of the regret
on the minimal mass $2\epsilon = \min\{f_m,\ldots,f_1\}$.
The following lemma shows that, in general,
the regret of the optimal policy grows with $\epsilon$.
The example was provided to us by \cite{kleinberg2017personalCommunications}.
The proof is in Appendix \ref{app:KleinbergExample}.

\begin{figure}[t!]

    \caption{\textbf{Dynamic-programming and Budget-Ratio regrets in \citeauthor{kleinberg2017personalCommunications}'s example}}\label{fig:kleinberg}

    \begin{center}
        \includegraphics[width=0.48\textwidth]{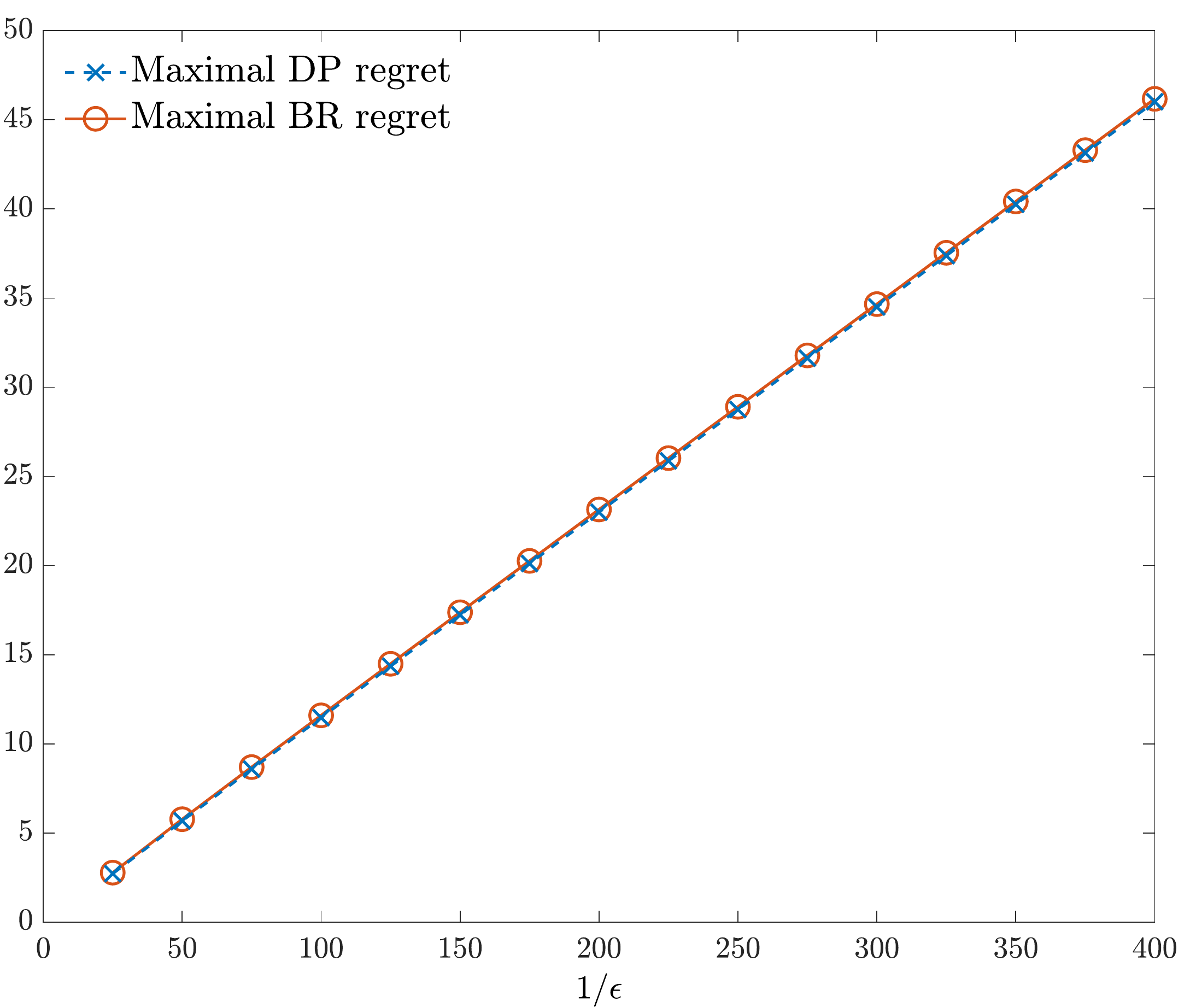}
        \hfill
        \includegraphics[width=0.48\textwidth]{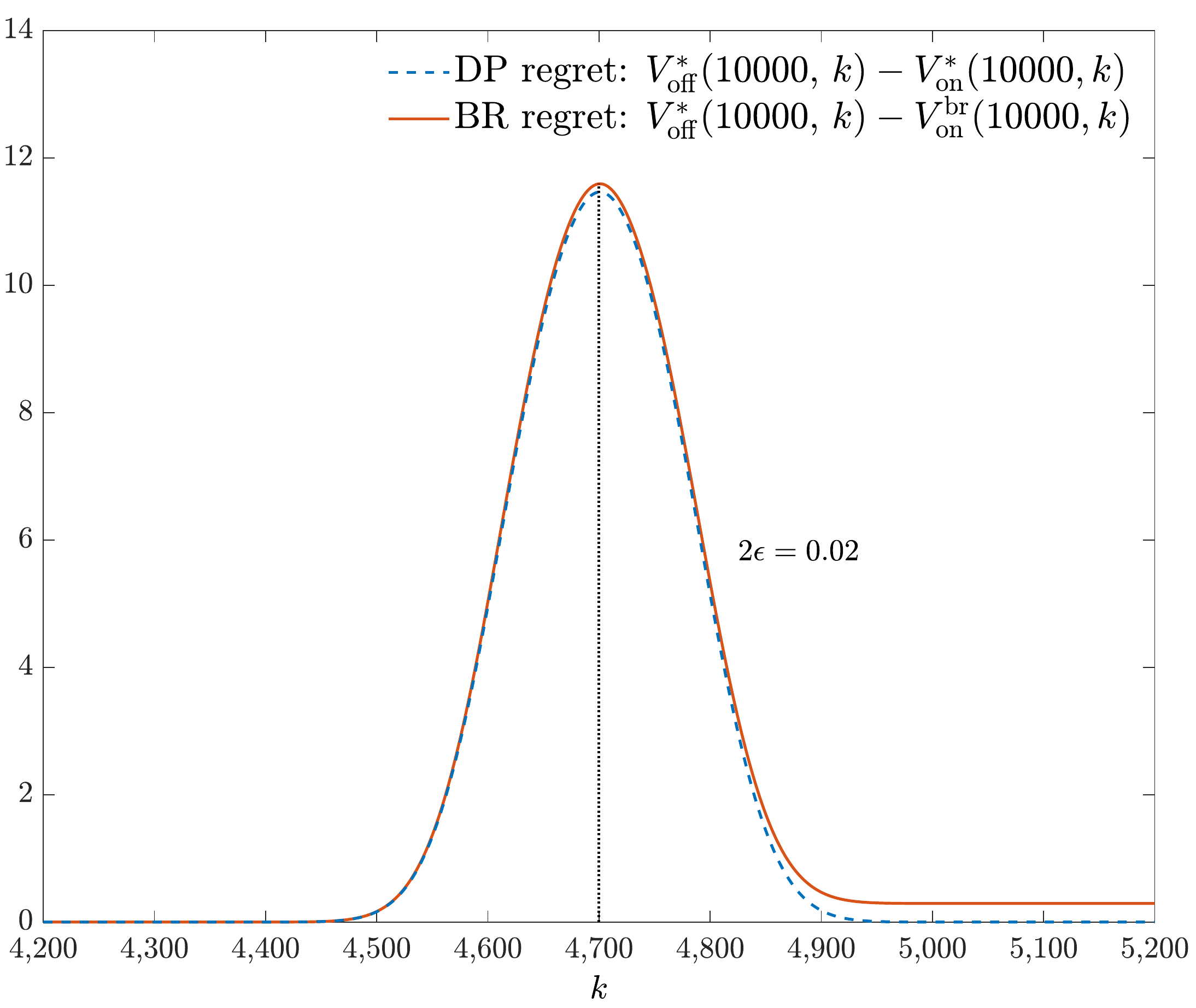}	
    \end{center}	

    \footnotesize{\emph{Notes.}
    \citeauthor{kleinberg2017personalCommunications}'s example on $\A = \{1,2,3\}$ and probability mass function
    $(\frac{1}{2} + 2 \epsilon, 2 \epsilon, \frac{1}{2} - 4\epsilon)$.
    On the left, we see that the regret grows linearly with $1/\epsilon$
    (as $n_{\epsilon}=\lceil 1/\epsilon^2 \rceil$ ).
    On the right, we take $\epsilon=0.01$ and $n=10,000$ and note that the regret of
    both the dynamic-programming and the Budget-Ratio policies is maximized at $k=4,700$,
    which corresponds to the threshold $T_2 = \frac{1}{2}(\Fbar(a_2)+\Fbar(a_3)) = 0.47$
    which will be key in our Budget-Ratio policy for this setting.}

\end{figure}

\begin{lemma} [Regret lower bound in \citeauthor{kleinberg2017personalCommunications}'s example]\label{lem:kleinberg}
Consider the three-point distribution on $\mathcal{A}=\{a_3,a_2,a_1\}$ with
probability mass function $(\frac{1}{2}+2\epsilon, 2\epsilon,\frac{1}{2} - 4\epsilon)$
and let $n_{\epsilon}=\lceil \frac{1}{\epsilon^2}\rceil$ and $k_{\epsilon}=\lceil n_{\epsilon}/2\rceil $.
Then, there exists a constant $\Gamma>0$ such that
$$
\frac{\Gamma}{\epsilon}
    \leq \Voff^{*}(n_{\epsilon},k_{\epsilon})- \Von^*(n_{\epsilon},k_{\epsilon})
    \leq \Voff^{*}(n_{\epsilon},k_{\epsilon})- \Von^{\BR}(n_{\epsilon},k_{\epsilon}).
$$
Thus, we also have that
$$
\frac{\Gamma}{\epsilon}
    \leq \sup_{(n,k)\in \mathcal{T}} \left\{ \Voff^*(n,k) - \Von^*(n,k) \right\}
    \leq \sup_{(n,k)\in \mathcal{T}} \left\{ \Voff^*(n,k) - \Von^{\BR}(n,k) \right\}.
$$
\end{lemma}

Thus, one cannot generally expect a bound that does not depend on the minimal mass.
Figure \ref{fig:kleinberg} shows that in this case too, the Budget-Ratio policy we develop in this paper
performs remarkably well relative to the optimal.
It is just that, as $\epsilon$ shrinks, the regret of both grows linearly with $1/\epsilon$.

\begin{remark}[On the regret with uniform random permutations]
We note here that there are regret upper bounds for a version of this multi-secretary problem
in which the values are given by a uniform random permutation of the integers $[n]$
instead of being from a random sample.
Within the uniform random permutation framework,
\cite{kleinberg2005multiple} proves that the minimal regret in this setting is of the order of $\sqrt{k}$
and provides an algorithm that achieves this lower bound.
\end{remark}

\section{The offline problem}

Denoting by $Z_j^r = \sum_{t \in [r]}\1(X_t=a_j)$ the number of candidates with ability $a_j$ inspected up to and including time $r$, the offline optimization problem has the compact representation
\begin{equation}\label{eq:Voff-star}
\Voff^{*}(n, k) = \E\left[  \varphi(Z_1^n, \ldots, Z_m^n, k)  \right ],
\end{equation}
where, for $(z_1, \ldots, z_m) \in \Z_+^m$ and $k \in \Z_+$,
\begin{eqnarray}\label{eq:LP}
\varphi(z_1, \ldots, z_m, k)  = & \underset{s_1, \ldots, s_m}{\max} &  \sum_{j\in[m]} a_j s_j\\
& \text{s.t.} &   0 \leq s_j\leq z_j \mbox{ for all } j\in[m] \notag\\
&  & \sum_{j\in [m]}s_j \leq k.  \notag
\end{eqnarray}

For a given realization of $Z_1^n,\ldots,Z_m^n$, the (trivial) optimal solution is to sort the values
and select the $k$ candidates with the largest abilities.
That is, one selects $\frakS_1^n = \min\{Z_1^n, k\}$ candidates with ability $a_1$.
Then, if there is any recruiting budget left, one selects candidates with ability $a_2$
until either selecting all of them or exhausting the remaining budget of $k-Z_1^n$,
i.e. one selects $\frakS_2^n = \min\{ Z_2^n, (k-Z_1^n)_+ \}$ candidates with ability $a_2$.
In general, the offline number of selected $a_j$-candidates is given by
\begin{equation}\label{eq:Sj-definition}
\frakS_j^n = \min\big\{ Z_j^n , ( k-\sum_{i\in [j-1]}Z_i^n )_+ \big\}
\quad \quad  \text{for each $j \in [m]$,}
\end{equation}
so that
\begin{equation}\label{eq:Voff}
\Voff^{*}(n,k)
=\sum_{j\in [m]} a_j \E\left[\frakS_j^n\right]= \sum_{j\in [m]} a_j \E\left[ \min\big\{ Z_j^n , ( k-\sum_{i\in [j-1]}Z_i^n )_+ \big\}\right].\\[10pt]
\end{equation}

The offline solution \eqref{eq:Sj-definition} has the appealing property
that---up to deviations that are constant in expectation---all
the action is within two ability levels.
That is, depending on the ratio $k/n$, there is an index $j_0 \in [m]$ such
that the offline sort algorithm rejects almost all of the candidates with ability strictly below $a_{j_0+1}$ and selects all the candidates with ability strictly above $a_{j_0}$.
The action index $j_0$ depends on the horizon length $n$ and on the recruiting budget $k$,
and it is given for any pair $(n,k)\in\mathcal{T}$ by
\begin{equation}\label{eq:j0}
j_0 (n,k) = \begin{cases}
              1 & \mbox{if } k/n < f_1 + f_2/2 \\
              j: \Fbar(a_{j}) +  f_j/2 \leq k/n < \Fbar(a_{j+2})- f_{j+1}/2
                        & \mbox{if } f_1 + f_2/2 \leq k/n < 1-f_m/2\\
              m & \mbox{if } 1 - f_m/2 \leq k/n.
            \end{cases}
\end{equation}
If $(n,k)$ are so that $j_0(n,k) = j$, then the two values at play are $a_j$ and $a_{j+1}$
and this suggests partitioning $\T$ into the sets
\begin{equation*}
  \T_j  = \{ (n, k) \in \T: j_0(n,k) = j \}, ~~  j \in \{1,2, \ldots, m\},
\end{equation*}
to obtain a helpful decomposition of the offline value \eqref{eq:Voff} .

\begin{proposition}[Offline sort decomposition]\label{prop:offline-decomposition}
    Let $\epsilon = \tfrac{1}{2} \min\{f_m, \ldots, f_1\}$ and fix $j \in \{1, \ldots, m\}$.
    For all $(n, k) \in \T_j$ one has the bounds
	\begin{equation}\label{eq:offline-lemma-count-bounds}
\sum_{i\in [j-1]}\E[Z_i^n] - \frac{1}{4 \epsilon} \leq  \sum_{i\in [j-1]}\E[\frakS_i^n]
	\quad \quad \text{and} \quad \quad
	\sum_{i=j+2}^{m}\E[\frakS_i^n]\leq \frac{1}{4 \epsilon}.
	\end{equation}
    Consequently, for all $(n, k) \in \T_j$ one has the decomposition
	\begin{equation}\label{eq:voff-decomposition}
	\Voff^{*}(n,k)
	= \sum_{i\in [j-1]} a_i\E[Z_i^n]
	+ a_{j}\E\left[\frakS_{j}^n\right]
	+ a_{j+1}\E\left[\frakS_{j+1}^n\right] \pm \frac{a_1}{4 \epsilon}.
	\end{equation}
\end{proposition}

The following lemma will be used repeatedly in the sequel.
All lemmas stated in the paper are proved in Appendix \ref{sec:proof-of-lemmas}.

\begin{lemma}\label{lm:multinomial-lemma}
	Let $B$ be a binomial random variable with $n$ trials
	and success probability $p$. Then, for any $\varepsilon > 0$,
	\begin{equation}\label{eq:binomial-lemma-E(X-k)bound}
    	\E[ (B- k)_+]  \leq \frac{1}{4\varepsilon}
    	\quad \text{ if } p+\varepsilon \leq \frac{k}{n},
        \quad \quad \mbox{ and } \quad \quad
    	\E[(k - B)_+] \leq \frac{1}{4\varepsilon}
        \quad \text{ if } \frac{k}{n} \leq p - \varepsilon.
    \end{equation}
\end{lemma}

The first use of Lemma \ref{lm:multinomial-lemma} is in the proof of Proposition \ref{prop:offline-decomposition}.

\proof{Proof of Proposition \ref{prop:offline-decomposition}.}
If $j=1$, the left inequality of \eqref{eq:offline-lemma-count-bounds} reduces to $-(4 \epsilon)^{-1}\leq 0$,
and there is nothing to prove.
For $1<j$ and $\iota < j$, we obtain from \eqref{eq:Sj-definition} that
$$
0 \leq Z_\iota^n - \frakS_\iota^n
  = \Big(\min\Big\{ Z_\iota^n,  \sum_{i\in [\iota]} Z_i^n - k\Big\}\Big)_+
  \leq \Big(\sum_{i\in [\iota]} Z_i^n - k\Big)_+
  \leq \Big(\sum_{i\in [j-1]} Z_i^n - k\Big)_+,
$$
By the definition of $\T_j$ and \eqref{eq:j0}, we have that
$
\Fbar(a_{j}) + \epsilon \leq k/n,
$
so, since the sum  $ B = \sum_{i\in [j-1]} Z_i^n$
is binomial with parameters $n$ and $\Fbar(a_{j})= f_{j-1}+ \ldots + f_1$,
the first inequality in \eqref{eq:offline-lemma-count-bounds}
follows from \eqref{eq:binomial-lemma-E(X-k)bound}.

For the second inequality in \eqref{eq:offline-lemma-count-bounds},
notice that if $j \in \{m-1, m\}$ then there is nothing to prove.
Otherwise, if $j \leq m-2$ we have for all $\iota \geq j+2$ (or, equivalently for $\iota-1\geq j+1$) that
$$
\frakS^n_\iota =\min\left\{Z_\iota^n,  \Big(k - \sum_{i\in [\iota-1]}Z_i^n\Big)_+\right\} \leq \Big(k - \sum_{i\in[\iota-1]}Z_i^n\Big)_+
      \leq \Big(k - \sum_{i\in[j+1]}Z_i^n \Big)_+.
$$
The sum $B = \sum_{i\in [j+1]}Z_i^n$ is a binomial random variable
with success probability $\Fbar(a_{j+2})= f_{j+1}+ f_{j}+ \ldots + f_1$,
so that since $\epsilon = \tfrac{1}{2} \min\{f_m, f_{m-1}, \ldots, f_1\}$ and $(n,k) \in \T_j$
we have that
$
k/n \leq \Fbar(a_{j+2}) - \epsilon,
$
and the second inequality in \eqref{eq:offline-lemma-count-bounds}
again follows from \eqref{eq:binomial-lemma-E(X-k)bound}.

To conclude the proof of the proposition, we recall \eqref{eq:Voff} and
rewrite $\Voff^{*}(n,k)$ as
\begin{equation}\label{eq:voff-decomposition-step0}
\Voff^{*}(n,k)
	  = - \!\! \sum_{i\in [j-1]}  \!\! a_i  \{ \E[Z^n_i] - \E\left[\frakS_i^n\right] \}
        + \!\! \sum_{i\in [j-1]} \!\!  a_i \E [Z_i^n]
        + a_{j}   \E\left[\frakS_j^n\right]
        + a_{j+1} \E\left[\frakS_{j+1}^n\right]
        + \!\! \sum_{i=j+2}^m \!\! a_i \left[\frakS_i^n\right].
\end{equation}
The definition of $\frakS_i^n$ in \eqref{eq:Sj-definition}, the monotonicity $a_m< a_{m-1}< \cdots <a_1$,
and the left inequality of \eqref{eq:offline-lemma-count-bounds} then give us that
\begin{equation}\label{eq:offline-decomposition-bound1}
	0 \leq \sum_{i\in [j-1]}  a_i  \{ \E[Z^n_i] - \E\left[\frakS_i^n\right] \}
      \leq a_1 \left( \sum_{i\in [j-1]} \E[Z^n_i] - \sum_{i\in [j-1]} \E\left[\frakS_i^n\right] \right)
      \leq \frac{a_1}{4 \epsilon}.
\end{equation}
Similarly, the monotonicity $a_m< a_{m-1}< \cdots <a_1$ and the right inequality of \eqref{eq:offline-lemma-count-bounds} then imply that
\begin{equation}\label{eq:offline-decomposition-bound2}
	0 \leq \sum_{i=j+2}^m  a_i \left[\frakS_i^n\right] \leq \frac{a_1}{4 \epsilon},
\end{equation}
so the decomposition \eqref{eq:voff-decomposition} follows after one estimates the first and the last
sum of \eqref{eq:voff-decomposition-step0} with the bounds in \eqref{eq:offline-decomposition-bound1} and \eqref{eq:offline-decomposition-bound2} respectively.
\halmos\endproof\vspace{0.25cm}

For any feasible online policy $\pi \in \Pi(n,k)$, we let
$S_j^{\pi,r}=\sum_{t=1}^{r}   \sigma_t^{\pi} \1(X_t=a_j)$
be the number of candidates with ability level $a_j$ that are
\emph{selected} by policy $\pi$ up to and including time $r$,
so the expected total ability accrued by policy $\pi$ can be written as
\begin{equation*}
\Von^{\pi}(n,k)
= \E\left[\sum_{t\in [n]} X_t \sigma_t^{\pi}\right]
=  \sum_{j\in [m]} a_j \E[S_j^{\pi,n}].
\end{equation*}

Proposition \ref{prop:offline-decomposition} suggests that, to have bounded regret, an online algorithm must be selecting almost all candidates with
abilities $a_1,a_2,\ldots, a_{j_0-1}$ and
rejecting almost all of the candidates with abilities $a_{j_0+2},\ldots,a_m$ if $j_0$
is such that $k/n\in [\Fbar(a_{j_0}) + \tfrac{1}{2}f_{j_0}, \Fbar(a_{j_0+2})- \tfrac{1}{2}f_{j_0+1})$.
This sufficient condition guides us in the development of the Budget-Ratio policy in Section \ref{sec:heuristic}.

\begin{proposition}[A sufficient condition]\label{prop:sufficient-condition-for-policy}
    Let $\epsilon = \tfrac{1}{2}\min\{f_m, \ldots, f_1\}$.
    For any given $(n,k) \in \T$ let $j_0 \equiv j_0(n,k)$ be the index defined by \eqref{eq:j0}
    and suppose there is a feasible online policy $\pi \equiv \pi(n,k) \in \Pi(n, k)$,
    a stopping time $\tau \equiv \tau(\pi) \leq n$, and a constant $M\equiv M(\epsilon)<\infty$ such that
    \begin{enumerate}[\rm (i)]	
		\item
		      $\sum_{j\in [j_0-1]}\E[S_j^{\pi, \tau}]=\sum_{j\in [j_0-1]}\E[Z_j^{\tau}]$,
		\item
        	  $\E[S_{j_0}^{\pi, \tau}]
               \geq  \E[\frakS_{j_0}^{\tau}] - M$,

		\item
		      $\E[S_{j_0+1}^{\pi, \tau}]
               \geq \E[\frakS_{j_0+1}^{\tau}] - M$,

        \item
              $\E[\tau]\geq n-M$.
	\end{enumerate}
    Then, one has that
	$$
     \sup_{(n,k) \in \T} \left\{ \Voff^{*}(n,k) -\Von^{\pi}(n,k) \right\}
     \leq 3 a_1 M + \frac{a_1}{4 \epsilon}.
    $$  	
\end{proposition}

The regret bound in Proposition \ref{prop:sufficient-condition-for-policy} has two components (summands). The first one accounts for three inefficiencies of the online policy $\pi$ with respect to the offline solution.
There is the cost for running out of budget too early---at most $M$ periods too early in expectation according to
condition (iv)---which cannot exceed $a_1 M$, and there are the costs that come from conditions (ii) and (iii)
that can each contribute to a maximal expected loss of $a_1 M$. Finally, the second summand accounts for the values that the offline solution might be choosing
and that are strictly smaller than $a_{j_0+1}$.
However, Proposition \ref{prop:offline-decomposition} tells us that their cumulative expected
value is at most $a_1 / (4 \epsilon)$.

\proof{Proof of Proposition \ref{prop:sufficient-condition-for-policy}.}
Since $\T = \bigcup_{j\in [m]} \T_j$, we have that
$$
\sup_{ (n,k) \in \T } \left\{ \Voff^{*}(n,k) -\Von^{\pi}(n,k) \right\}
= \max_{j\in[m]} \left\{ \sup_{ (n,k) \in \T_j } \left\{ \Voff^{*}(n,k) -\Von^{\pi}(n,k) \right\} \right\},
$$
and it suffices to verify that, for any $j\in [m]$,
$$
  \sup_{ (n,k) \in \T_j } \left\{ \Voff^{*}(n,k) -\Von^{\pi}(n,k) \right\}
    \leq 3 a_1 M + \frac{a_1}{4 \epsilon}.
$$
For any $(n, k) \in \T_j$ the definition \eqref{eq:j0} of the map $(n,k) \mapsto j_0(n,k)$ gives us that $j_0(n,k) = j$,
so for any policy $\pi \in \Pi(n,k)$ and any stopping time $\tau \leq n$
we have the lower bound
$$
 \sum_{i\in [j-1]} a_i \E[S_i^{\pi, \tau}]
 + a_{j} \E[S_{j}^{\pi, \tau}]
 + a_{j+1} \E[S_{j+1}^{\pi, \tau}]
 \leq  \sum_{j\in [m]} a_j \E[S_j^{\pi, n}] = \Von^{\pi}(n, k).
$$
In turn, for the policy $\pi \equiv \pi(n,k)$ and the stopping time $\tau \equiv \tau(\pi)$
given in the proposition, it follows from the requirements (i), (ii), and (iii) that there is a constant $M \equiv M(\epsilon)$ such that
$$
\sum_{i\in [j-1]} a_i \E[Z_i^{\tau}] + a_{j} \{ \E[\frakS_{j}^{\tau}] - M\} + a_{j+1} \{ \E[\frakS_{j+1}^{ \tau}]  - M\}
\leq \Von^{\pi}(n, k),
$$
which, since $a_m < a_{m-1} < \ldots < a_1$, gives
\begin{equation}\label{eq:Von-lower-bound}
\sum_{i\in [j-1]} a_i \E[Z_i^{\tau}]
    + a_{j} \E[\frakS_{j}^{\tau}] + a_{j+1} \E[\frakS_{j}^{\tau}]  - 2 a_1 M
\leq \Von^{\pi}(n, k).
\end{equation}

For $\tau \leq n$ we now rewrite the upper bound in \eqref{eq:voff-decomposition} as
\begin{align}\label{eq:voff-upperbound-tau}
	\Voff^{*}(n,k)
		& \leq  \sum_{i\in [j-1]} a_i\E[Z_i^\tau]
        	   + a_{j}\E\left[\frakS_{j}^\tau \right]
			   + a_{j+1}\E\left[\frakS_{j+1}^\tau \right]  + \frac{a_1}{4 \epsilon} \\
             & + \sum_{i\in [j-1]} a_i\E[Z_i^n - Z_i^\tau]
			   + a_{j}\E\left[\frakS_{j}^n - \frakS_{j}^\tau \right]
			   + a_{j+1}\E\left[\frakS_{j+1}^n - \frakS_{j+1}^\tau \right], 	 \notag
\end{align}
and we obtain from the definition \eqref{eq:Sj-definition} of the offline number of selected $a_i$-candidates
that the difference
$0 \leq \frakS_{i}^n - \frakS_{i}^\tau \leq Z_i^n - Z_i^\tau = \sum_{t = \tau+1}^n \1(X_t = a_i)$
for all $i \in [m]$.
Thus, if we recall that $a_m < a_{m-1} < \ldots < a_1$, we obtain the bound
$$
\sum_{i\in [j-1]} \!\!\! a_i\E[Z_i^n - Z_i^\tau]
	+ a_{j}\E\left[\frakS_{j}^n - \frakS_{j}^\tau \right]
	+ a_{j+1}\E\left[\frakS_{j+1}^n - \frakS_{j+1}^\tau \right]
\leq a_1 \E\left[ \sum_{i \in [m]} \sum_{t=\tau+1}^n \!\!\! \1(X_t = a_i) \right] \!\!
  =  a_1 \E[n - \tau],
$$
so when we use Property (iv) to estimate this last right-hand side
and replace this estimate in \eqref{eq:voff-upperbound-tau} we finally have that
\begin{equation}\label{eq:voff-upperbound-tau2}
\Voff^{*}(n,k)
 	 \leq  \sum_{i\in [j-1]} a_i\E[Z_i^\tau]
        	   + a_{j}\E\left[\frakS_{j}^\tau \right]
			   + a_{j+1}\E\left[\frakS_{j+1}^\tau \right]
               + a_1 M
               + \frac{a_1}{4 \epsilon}.
\end{equation}
The proof of the proposition then follows by combining the bounds \eqref{eq:Von-lower-bound}
and \eqref{eq:voff-upperbound-tau2}.
\halmos\endproof \vspace{0.25cm}

\begin{remark}[A deterministic relaxation.]
A further upper bound to the multi-secretary problem is provided by an intuitive deterministic relaxation. Given the linear program \eqref{eq:LP}, we consider its optimal value
\begin{equation}\label{eq:DR}
DR(n,k)=  \varphi ( \E[Z_1^n],\ldots, \E[Z_m^n], k)
\end{equation}
and we note that its optimal solution is given by
\begin{equation}\label{eq:DR-optimal-solution}
s^*_j = \min \big\{ \E[Z_j^n], (k - \sum_{i\in[j-1]}\E[Z_i^n])_+ \big\}
      = \min \big\{  n f_j , (k - n \Fbar(a_j))_+ \big\}
      \quad\quad \text{ for all } j\in [m].
\end{equation}
Since the expected number of $a_j$-candidates selected by the offline solution
$
\widehat s_j = \E[\frakS_j^n]= \E[ \min\{ Z_j^n, (k - \sum_{i\in[j-1]}Z_i^n)_+ \} ]
$
is a feasible solution for the deterministic relaxation \eqref{eq:DR},
we immediately have the bound
\begin{equation*}
\Voff^*(n,k) \leq DR(n,k)  \quad \text{for all pairs } (n,k) \in \T.
\end{equation*}
As central-limit-theorem intuition suggests,
the ``cost of randomness'' is at most of the order of $\sqrt{n}$. But it does not have to be that large: there is a subset
$\T' \subset \T$ for which the difference $DR(n,k) - \Voff^{*}(n,k)$
is bounded by a constant that does not depend on $(n, k) \in \T'$.
Members $(n,k)$ of $\T'$ are such that $k/n$ is ``safely'' away from the jump points of the distribution $F$
(Appendix \ref{sec:det_relax}).
For $(n,k)\in \T'$, benchmarking against the deterministic relaxation
is the same as benchmarking against the offline sort \citep[see also][]{wu2015algorithms}.
In general this is not the case. For instance, if one takes $k=\E[Z_1^n]=nf_1$ then
$DR(n,k)=a_1 n f_1$ but there exists $M<\infty$ such that
$a_1\E[\mathfrak{S}_1^n] = a_1\E[ \min\{ Z_1^n, k \} ]=a_1 nf_1+ \E[\min\{Z_1^n-nf_1,0\}]	
\leq a_1nf_1-M\sqrt{n}=DR(n,k)-M\sqrt{n}$.
The first inequality follows from the approximation of the centered binomial $Z_1^n-nf_1$
by a normal with standard deviation $\sqrt{nf_1(1-f_1)}$ as in Lemma \ref{lm:martingale} below.
\end{remark}

\section{The Budget-Ratio (BR) policy \label{sec:heuristic}}

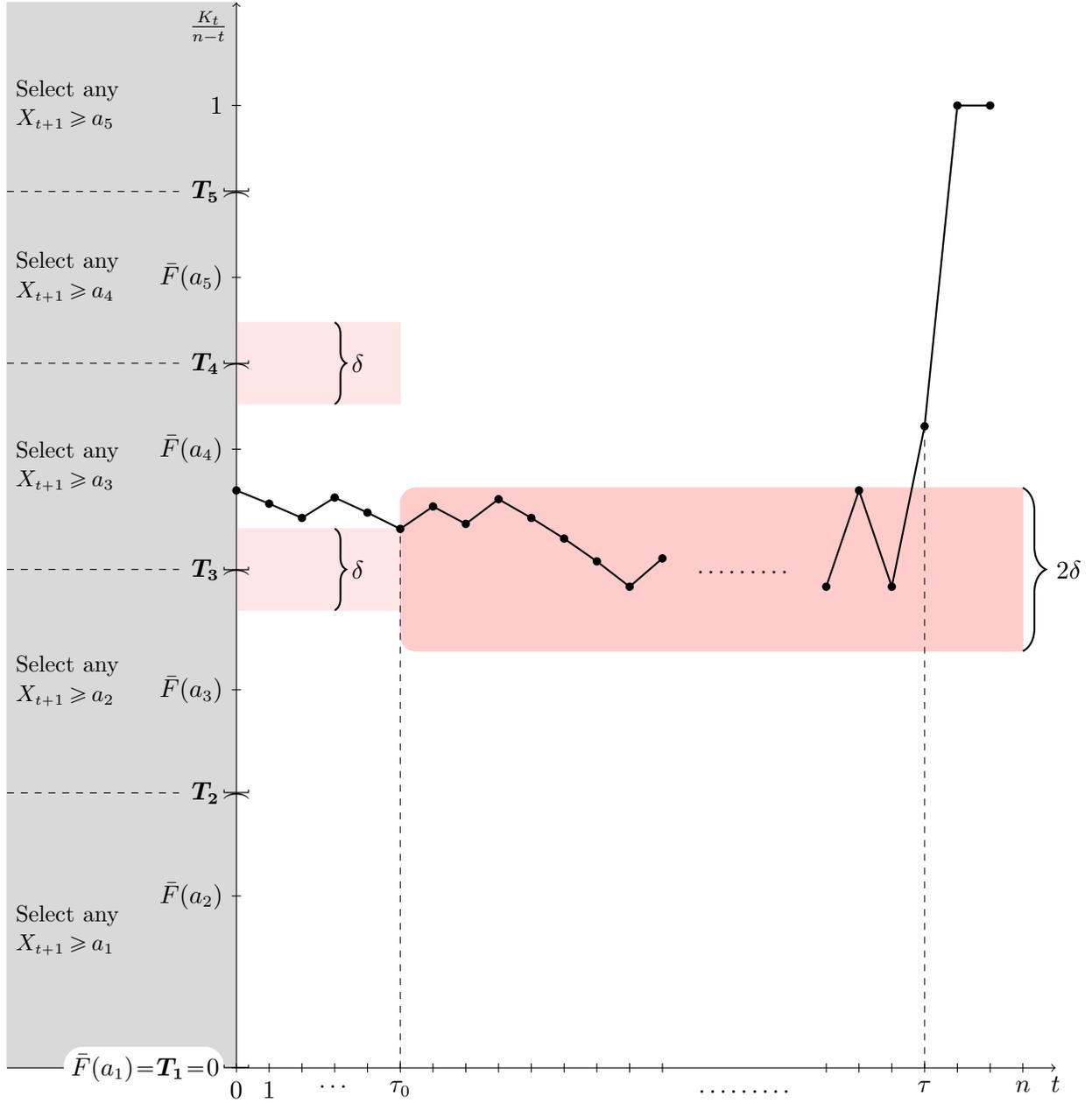
\begin{figure}[ht!]

\caption{\textbf{The BR policy: thresholds and dynamics.}}\label{fig:BR-policy}

\begin{center}
	\scalebox{0.995}{\begin{tikzpicture}[y=.525cm, x=.5cm]

    \draw [fill=gray!30,gray!30] (-7,0) rectangle (0,31);

    \draw [fill=red!10,red!10] (0,14.5-1.1875) rectangle (5,14.5+1.1875);

    \draw [thick,decorate,decoration={brace,amplitude=5pt,mirror},xshift=0pt,yshift=0pt]
    (3,14.5-1.1875) -- (3,14.5+1.1875) node [black,midway,xshift=10pt]
    {$\delta$};

    \draw [fill=red!10,red!10] (0,20.5-1.1875) rectangle (5,20.5+1.1875);

    \draw [thick,decorate,decoration={brace,amplitude=5pt,mirror},xshift=0pt,yshift=0pt]
    (3,20.5-1.1875) -- (3,20.5+1.1875) node [black,midway,xshift=10pt]
    {$\delta$};

    \draw (-7,0) -- (-5.25,0);
    \draw [fill=white,white, rounded corners = 8pt] (-5.25,-.6) rectangle (-.4,.6);

   \fill [red!20,draw]
        (5,14.5-2.375) --
        (24, 14.5-2.375) --
        (24,14.5+2.375) {[rounded corners=7pt]--
        (5,14.5+2.375) --
        cycle}
        {};

    \draw [thick,decorate,decoration={brace,amplitude=10pt,mirror},xshift=0pt,yshift=0pt]
    (24,14.5-2.375) -- (24,14.5+2.375) node [black,midway,xshift=20pt]
    {$2 \delta$};

	\draw [<->] (-0,31) node (yaxis) [anchor= north east] {$\tfrac{K_t}{n-t}$}
	|- (25,0) node (xaxis) [below] {$t$};
	
    \draw (0,2pt) -- (0,-2pt)
	node[anchor=north] {$0$};

    \draw (1,2pt) -- (1,-2pt)
	node[anchor=north] {$1$};

    \draw (2,2pt) -- (2,-2pt)
	node[anchor=north] {};

    \draw (3,2pt) -- (3,-2pt)
	node[anchor=north] {$\cdots$};

    \draw (4,2pt) -- (4,-2pt)
	node[anchor=north] {};

    \draw (5,2pt) -- (5,-2pt)
	node[anchor=north] {$\tau_0$};

    \draw (6,2pt) -- (6,-2pt)
	node[anchor=north] {};

    \draw (7,2pt) -- (7,-2pt)
	node[anchor=north] {};

    \draw (8,2pt) -- (8,-2pt)
	node[anchor=north] {};

    \draw (9,2pt) -- (9,-2pt)
	node[anchor=north] {};

    \draw (10,2pt) -- (10,-2pt)
	node[anchor=north] {};

    \draw (11,2pt) -- (11,-2pt)
	node[anchor=north] {};

    \draw (12,2pt) -- (12,-2pt)
	node[anchor=north] {};

    \draw (13,2pt) -- (13,-2pt)
	node[anchor=north] {};

    \draw (18,2pt) -- (18,-2pt)
	node[anchor=north] {};

    \draw (19,2pt) -- (19,-2pt)
	node[anchor=north] {};

    \draw (20,2pt) -- (20,-2pt)
	node[anchor=north] {};

    \draw (22,2pt) -- (22,-2pt)
	node[anchor=north] {};

    \draw (23,2pt) -- (23,-2pt)
	node[anchor=north] {};

    \draw (24,2pt) -- (24,-2pt)
	node[anchor=north] {$n$};


    \node[anchor=west,text width=2cm] at (-7,28){\small{Select any $X_{t+1} \geq a_5$}};

	\draw (.14,28) -- (-.14,28)
	node[anchor=east] {$1$};

    \node [rotate=90]	at (0,25.5) {[};
    \node [rotate=90]	at (0,25.425) {)};
    \node[anchor=east] at (-.25,25.5) {$\boldsymbol{T_5}$};

    \draw[dashed] (-7,25.5) -- (-1.75,25.5);

    \node[anchor=west,text width=2cm] at (-7,23){\small{Select any $X_{t+1} \geq a_4$}};	

	\draw (.14,23) -- (-.14,23)
	node[anchor=east] {$\Fbar(a_5)$};

    \node [rotate=90]	at (0,20.5) {[};
    \node [rotate=90]	at (0,20.425) {)};
    \node[anchor=east] at (-.25,20.5) {$\boldsymbol{T_4}$};

    \draw[dashed] (-7,20.5) -- (-1.75,20.5);

    \node[anchor=west,text width=2cm] at (-7,17.5){\small{Select any $X_{t+1} \geq a_3$}};

	\draw (.14,18) -- (-.14,18)
	node[anchor=east] {$\Fbar(a_4)$};

    \node [rotate=90]	at (0,14.5) {[};
    \node [rotate=90]	at (0,14.425) {)};
    \node[anchor=east] at (-.25,14.5) {$\boldsymbol{T_3}$};

    \draw[dashed] (-7,14.5) -- (-1.75,14.5);

    \node[anchor=west,text width=2cm] at (-7,11.25){\small{Select any $X_{t+1} \geq a_2$}};

	\draw (.14,11) -- (-.14,11)
	node[anchor=east] {$\Fbar(a_3)$};

    \node [rotate=90]	at (0,8) {[};
    \node [rotate=90]	at (0,7.925) {)};
    \node[anchor=east] at (-.25,8) {$\boldsymbol{T_2}$};

    \draw[dashed] (-7,8) -- (-1.75,8);

    \node[anchor=west,text width=2cm] at (-7,4){\small{Select any $X_{t+1} \geq a_1$}};

    \draw (.14,5) -- (-.14,5)
	node[anchor=east] {$\Fbar(a_2)$};

	\node [rotate=90]	at (0,0) {[};
	\draw (.14,0) -- (-.14,0)
	node[anchor=east] {$\Fbar(a_1) \!=\! \boldsymbol{T_1} \!=\! 0\,$};

    \draw[thick] (0, 28*18/30) -- (1,28*17/29) -- (2, 28*16/28) --
                 (3, 28*16/27) -- (4, 28*15/26) -- (5, 28*14/25) --
                 (6, 28*14/24) -- (7, 28*13/23) -- (8, 28*13/22) --
                 (9, 28*12/21) -- (10, 28*11/20) -- (11, 28*10/19) --
                 (12, 28*9/18) -- (13, 28*9/17);

    \node [fill, circle, scale=0.35] at (0, 28*18/30) {};
    \node [fill, circle, scale=0.35] at (1,28*17/29) {};
    \node [fill, circle, scale=0.35] at (2, 28*16/28) {};
    \node [fill, circle, scale=0.35] at (3, 28*16/27) {};
    \node [fill, circle, scale=0.35] at (4, 28*15/26) {};
    \node [fill, circle, scale=0.35] at (5, 28*14/25) {};
    \node [fill, circle, scale=0.35] at (6, 28*14/24) {};
    \node [fill, circle, scale=0.35] at (7, 28*13/23) {};
    \node [fill, circle, scale=0.35] at (8, 28*13/22) {};
    \node [fill, circle, scale=0.35] at (9, 28*12/21) {};
    \node [fill, circle, scale=0.35] at (10, 28*11/20) {};
    \node [fill, circle, scale=0.35] at (11, 28*10/19) {};
    \node [fill, circle, scale=0.35] at (12, 28*9/18) {};
    \node [fill, circle, scale=0.35] at (13, 28*9/17) {};

    \node at (15.5,28*157.5/306) {$\cdots \cdots \cdots$};

    \node[anchor=north] at (15.5,-6pt) {$\ldots \ldots \ldots$};

    \draw[thick] (18, 28*3/6) -- (19, 28*3/5) -- (20, 28*2/4) -- (21,28*2/3) -- (22, 28*2/2) --
                 (23, 28*1/1);

    \node [fill, circle, scale=0.35] at (18, 28*3/6) {};
    \node [fill, circle, scale=0.35] at (19, 28*3/5) {};
    \node [fill, circle, scale=0.35] at (20, 28*2/4) {};
    \node [fill, circle, scale=0.35] at (21,28*2/3) {};
    \node [fill, circle, scale=0.35] at (22, 28*2/2) {};
    \node [fill, circle, scale=0.35] at (23, 28*1/1) {};

    \draw[dashed] (21,0) -- (21,28*2/3);
	
    \draw (21,2pt) -- (21,-2pt)
	node[anchor=north] {$\tau$};

    \draw[dashed] (5,0) -- (5, 28*14/25);

	\end{tikzpicture}}
\end{center}

\footnotesize{\emph{Notes.}
            The $y$-axis has the thresholds of the BR policy for the 5-point distribution
            on $\A = \{a_5, a_4, a_3, a_2, a_1\}$ with the probability mass function
            $(f_5, f_4, f_3, f_2, f_1) = (\tfrac{5}{28}, \tfrac{5}{28}, \tfrac{7}{28}, \tfrac{6}{28}, \tfrac{5}{28})$.
            The plotted series is a sample path realization of the ratio $\{ \tfrac{K_t}{n-t}: 0 \leq t \leq n\}$
            which enters the ``orbit'' of the threshold $T_3$ at time $\tau_0$ (so $j(\tau_0)=3$) and exits at time $\tau$. Up until $\tau_0$ both thresholds $T_3$ and $T_4$ are in play.
            In this chart we take $\delta = \tfrac{19}{224}<\epsilon=\tfrac{5}{56}=\frac{1}{2}\min\{f_5,\ldots,f_1\}$.
            Notice that $\Fbar(a_3)=f_1+f_2$. When $\Fbar(a_3) <K_t/(n-t)<T_3$, the budget is, in expectation, sufficient to take some $a_3$ values but the policy will not do that until $T_3$ is crossed. This ``under-selection'' makes $K_t/(n-t)$ drift up toward $T_3$. When $T_3 <K_t/(n-t)<\Fbar(a_4)$, the budget is, in expectation, insufficient to take all $a_3$ values but the policy does select them. This ``over-selection'' makes $K_t/(n-t)$ drift towards $T_3$.}

\end{figure}

We now introduce an adaptive online policy that makes selection decisions depending on the ratio between the
remaining number of positions to be filled (the \emph{remaining budget}) and the remaining number of candidates to be inspected (the \emph{remaining time}),
and we refer to this policy as the \emph{Budget-Ratio (BR) policy}.

With $\pi=\text{br}$, the random variables $\sigmaBR_1, \sigmaBR_2, \ldots, \sigmaBR_n$
give us the sequence of selection decisions under the BR policy (see Section \ref{sec:model}),
and we let, for $t \in [n]$,
$$
K_t = k-\sum_{s\in[t]} \sigmaBR_s = K_{t-1} - \sigmaBR_t
$$
be the remaining budget after the $t$th decision ($K_0=k$).
We now introduce the thresholds $0 = T_1 < T_2 < \cdots < T_m < T_{m+1}= +\infty$ given by
$$
T_j = \frac{1}{2} ( \Fbar(a_{j}) + \Fbar(a_{j+1}))  \quad \quad  \text{ for each } j \in \{2, 3, \ldots, m\},
$$
so that the Budget-Ratio decision at time $t+1$
selects $X_{t+1}$ depending on its value and on the position of the ratio $K_t/(n-t)$ relative to these thresholds.
Specifically, at each decision time $t+1 \in \{1, \ldots n-1\}$, the BR policy
\begin{enumerate}[(i)]
  \item identifies the index $j\in[m]$ such that
        $$
            T_j \leq \frac{K_{t}}{n-t} < T_{j+1};
        $$
  \item selects $X_{t+1}$ if and only if $K_{t} >0$ and $X_{t+1} \geq a_{j}$;
        i.e., it sets
        $$
        \sigmaBR_{t+1} = \begin{cases}
                            1 & \mbox{if } K_{t} > 0 \text{ and } X_{t+1} > a_{j+1}\\
                            0 & \mbox{otherwise}.
                         \end{cases}
        $$
        Thus, the value $X_{t+1}$ is selected with probability $\P(\sigmaBR_{t+1}=1 | \F_t) = \Fbar(a_{j+1})\1(K_t > 0)$.
\end{enumerate}

Figure \ref{fig:BR-policy} gives a graphical representation of the selection regions of the BR policy. It encapsulates the source of its good performance. Consider an initial budget ratio $k/n \in [T_3,\Fbar(a_4))$ as in the figure.
In this range there is enough budget, in expectation, to select all of the
candidates with ability values $a_1$ and $a_2$.
The remaining budget after such selections should be used for candidates with ability
equal to $a_3$ or smaller. While the budget ratio is in the colored band, the policy will select all of the $a_1$ and $a_2$ values. It will also select $a_3$-candidates whenever the budget ration is above $T_3$. If we can guarantee that the budget ratio stays in the colored band almost until the end of the horizon, then the BR policy will select as many $a_3$ values as possible without giving up on any $a_1$ or $a_2$ values. This logic is formalized in the proof of Corollary \ref{cor:regret}.
The policy has two natural properties: (i) since $T_1 = 0$, the BR policy selects all $a_1$-valued candidates until exhausting the budget; and (ii) since $T_{m-1} < 1$ and $T_m = + \infty$, the BR policy selects all remaining values
as soon as the remaining budget is greater than or equal to the remaining number of time periods (i.e., if $n - t \leq K_t$).

\begin{figure}[ht!]
    \caption{\textbf{Simulated ratio and remaining budget averages}}\label{fig:heuristic}
    \begin{center}
        \includegraphics[width=.48\textwidth]{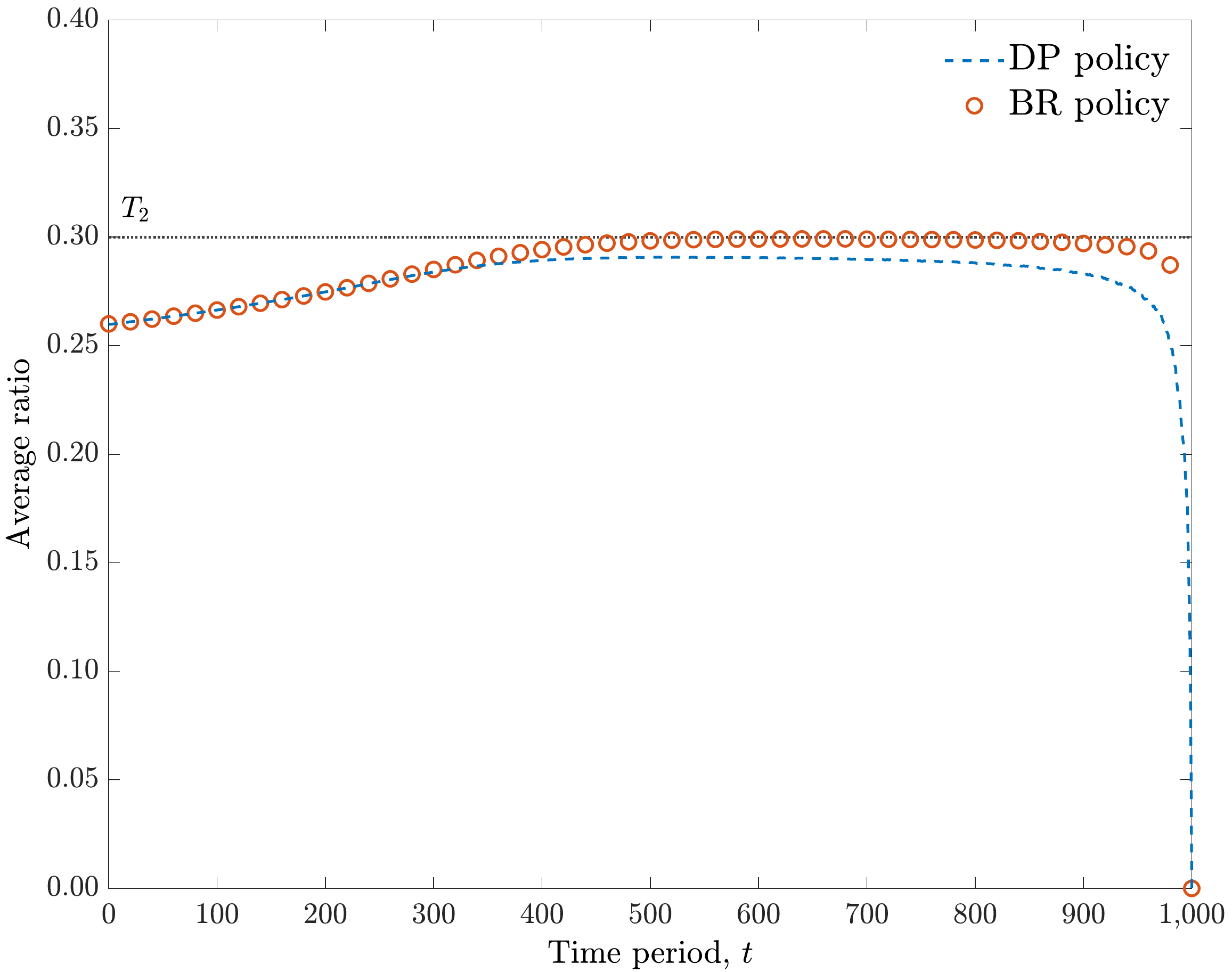}
        \hfill
        \includegraphics[width=.48\textwidth]{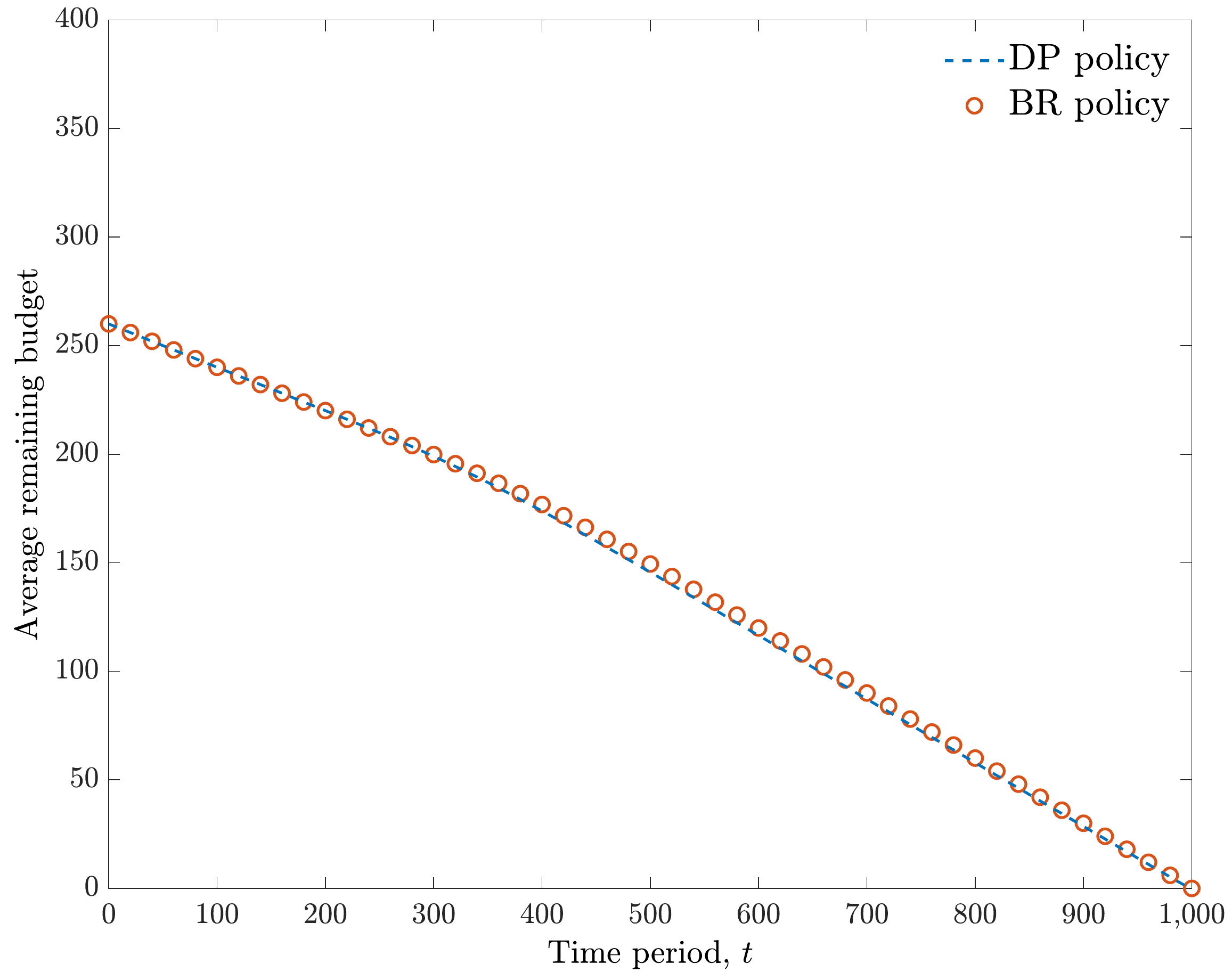}
        \\[7pt]
        \includegraphics[width=.48\textwidth]{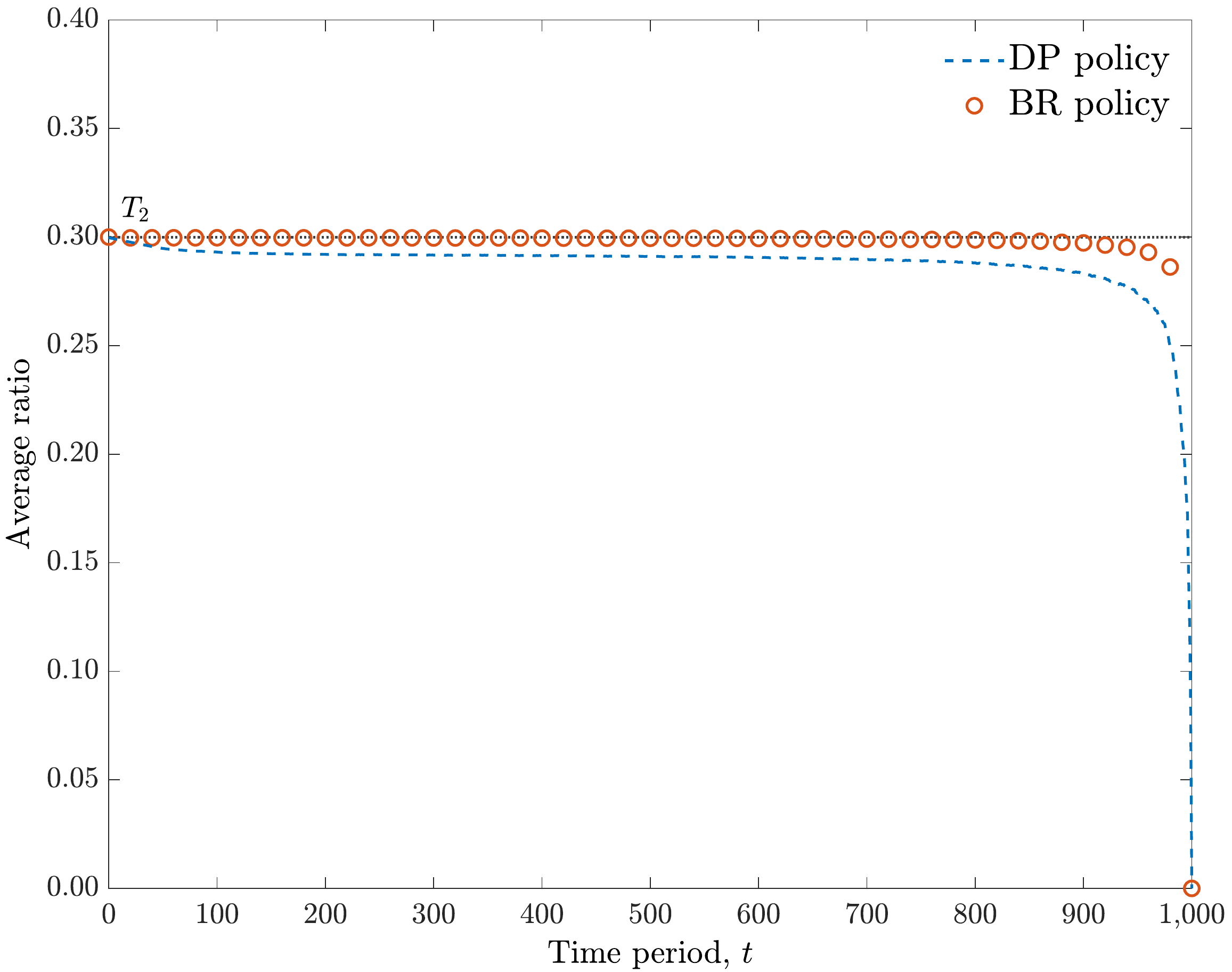}
        \hfill
        \includegraphics[width=.48\textwidth]{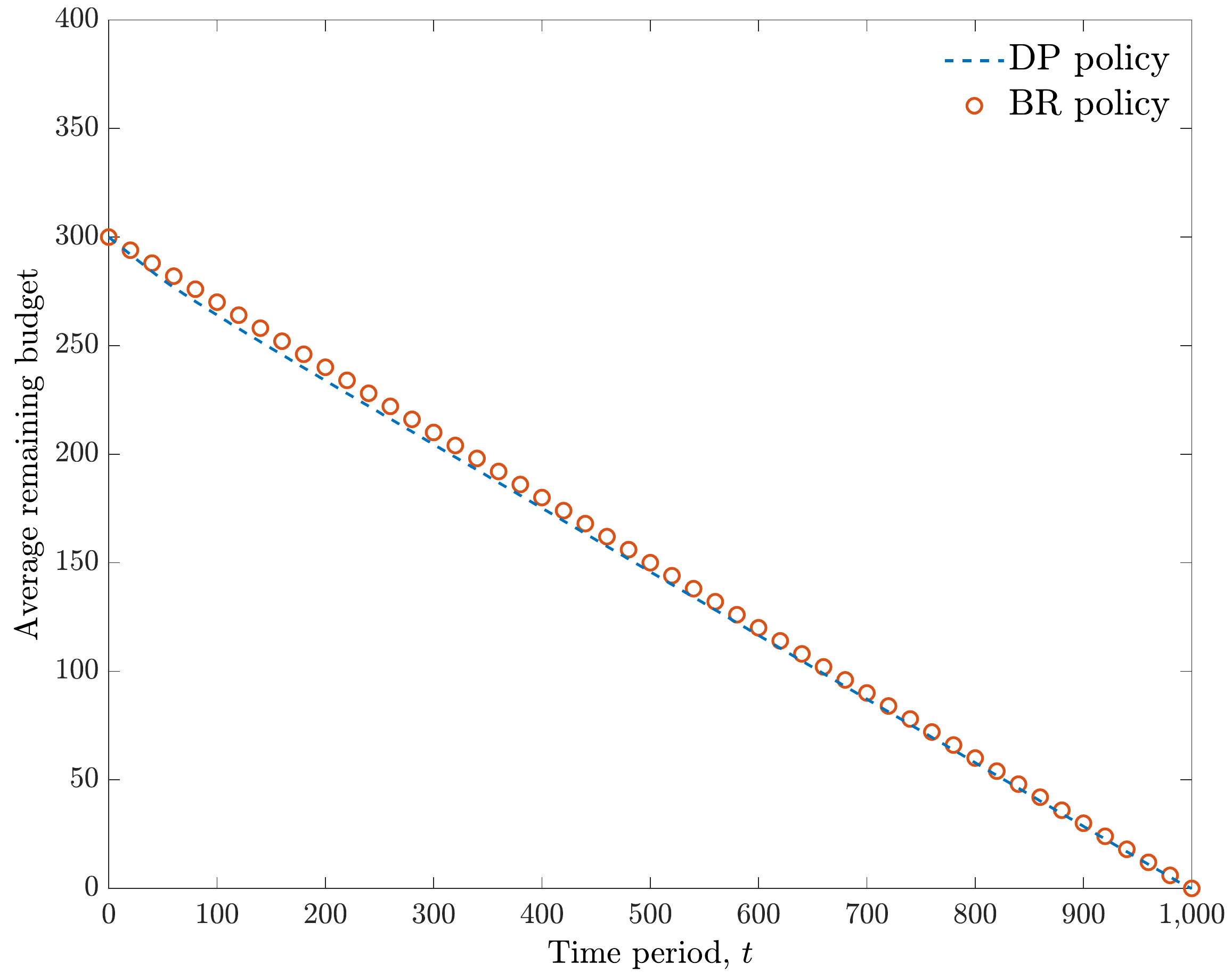}
        \\[7pt]
        \includegraphics[width=.48\textwidth]{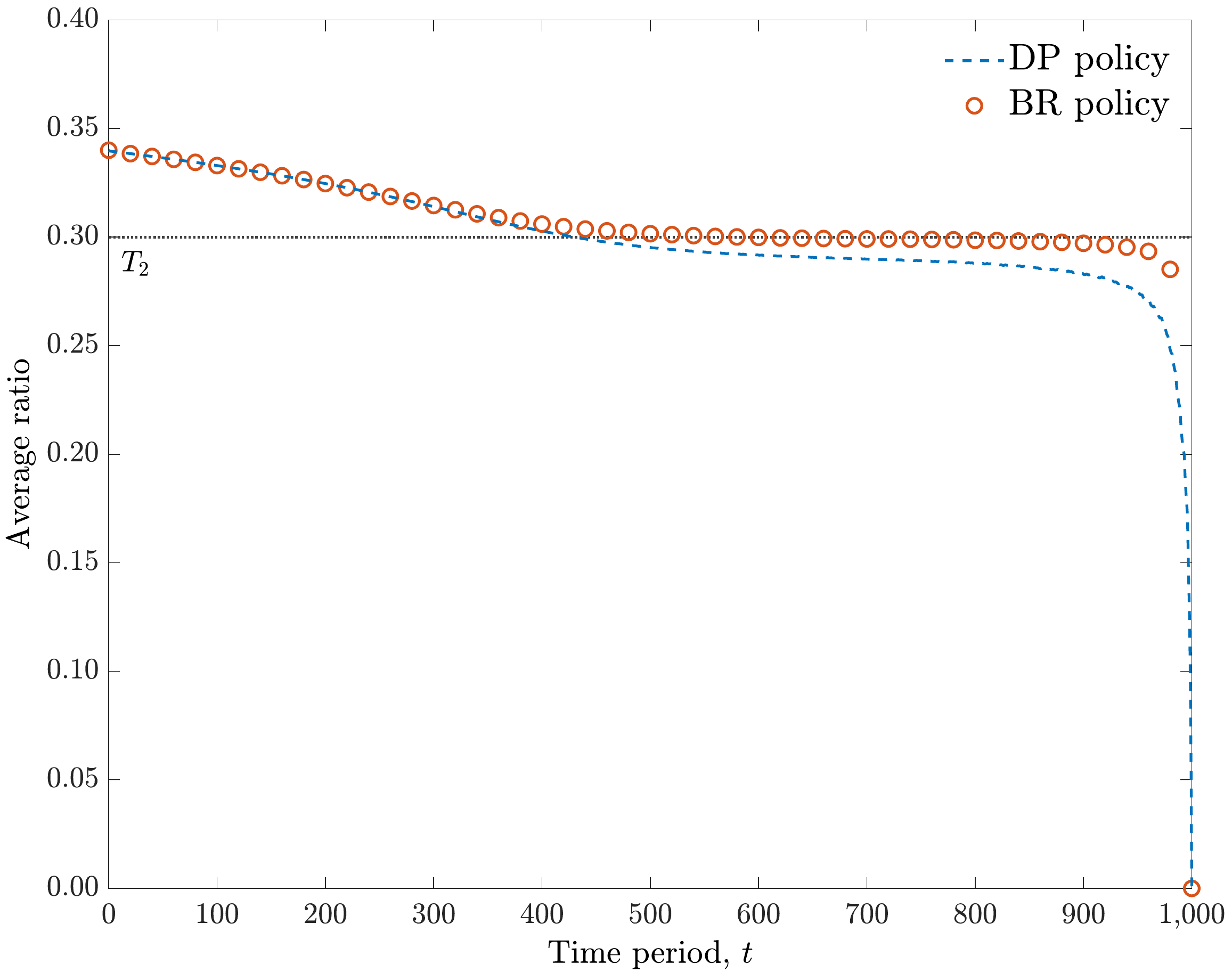}
        \hfill
        \includegraphics[width=.48\textwidth]{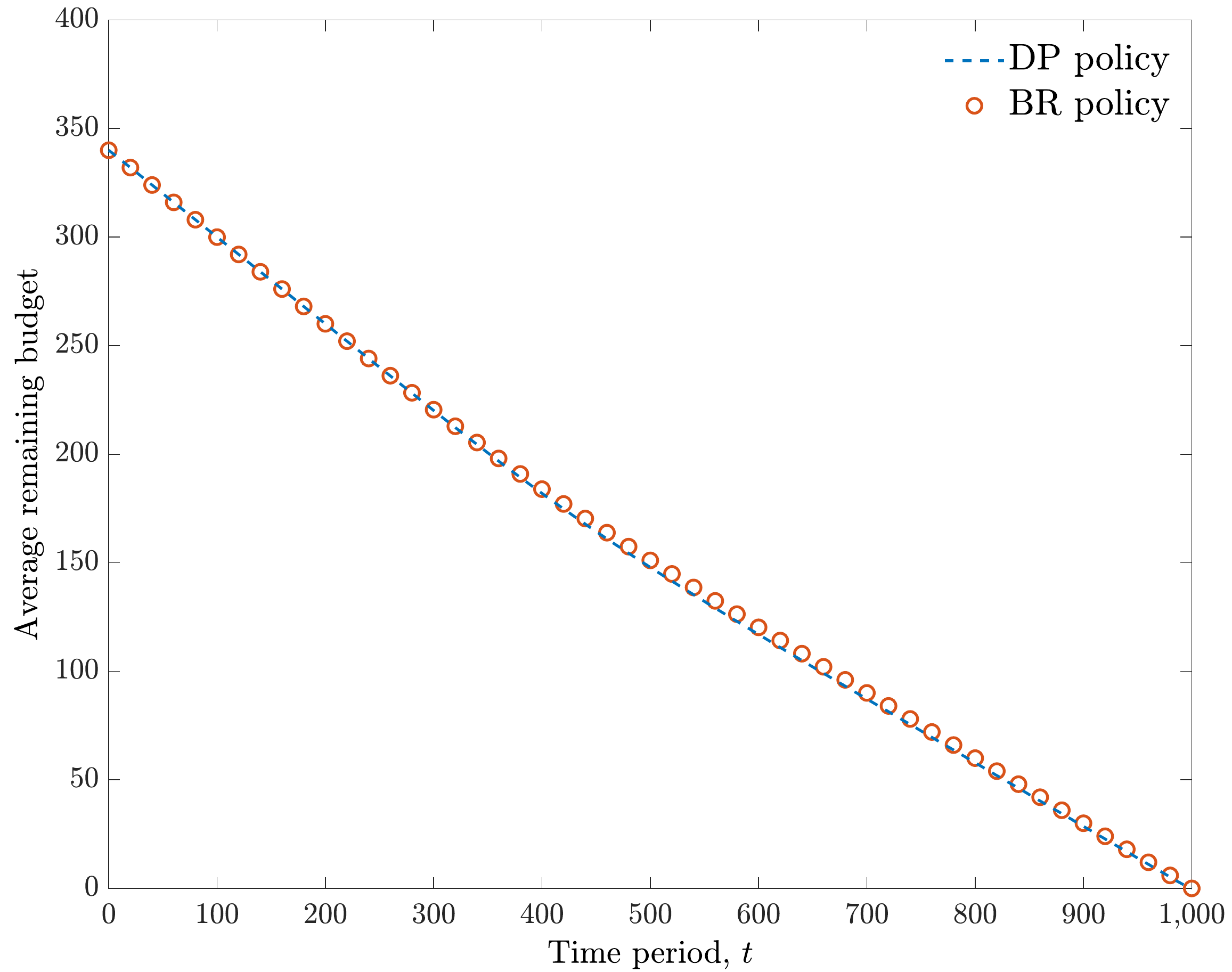}
    \end{center}
    \footnotesize{\emph{Notes.}
        Simulated averages---with $n=1,000$ and based on 10,000 trials---of
        the ratio $\tfrac{1}{n-t} \left( k - \sum_{s \in [t]} \sigma^\pi_s \right)$
        and of the remaining budget $ k - \sum_{s \in [t]} \sigma^\pi_s $
        when $\pi$ is the optimal dynamic-programming (DP) policy
        or the Budget-Ratio (BR) policy.
        The ability distribution is uniform on $\A = \{0.20, 0.65, 1.10, 1.55, 2.00\}$.
        In the plots, we vary the budget $k \in \{260,300,340\}$ starting with the smallest at the top
        and while keeping the ratio $k/n$ 
        within the orbit of $T_2 = 0.30$.
        }
\end{figure}

Next, we set
$$
\epsilon = \tfrac{1}{2}\min \{f_m, f_{m-1}, \ldots, f_1\},
$$
fix $ 0 < \delta < \epsilon$, and consider the stopping time
\begin{equation*}
  \tau_0
    =\inf \left\{t\geq 0: \left|\frac{K_t}{n-t} - T_j\right| \leq \frac{\delta}{2} \text{ for some } j \in [m]
            \text{ or } t \geq n - 2\delta^{-1}-1 \right\}.
\end{equation*}
If $\tau_0 < n - 2 \delta^{-1} - 1$, then $\tau_0$ is the first time that the ratio
$K_t/(n-t)$ enters the ``orbit'' of one of the thresholds; see Figure \ref{fig:BR-policy}.
We denote by $j(\tau_0)$ the index of the threshold that is within $\delta/2$ of the ratio $K_{\tau_0}/(n-\tau_0)$,
and we use $T_{j(\tau_0)}$ to denote the value of that threshold.
If $\tau_0 = n - 2 \delta^{-1} - 1$, then we set $j(\tau_0)=m+1$ and $T_{j(\tau_0)}=\infty$.
For all $t\leq n-2\delta^{-1}-1$ the jumps of $K_t/(n-t)$ satisfy
the absolute bound
$$
\left|\frac{K_t}{n-t}-\frac{K_{t+1}}{n-(t+1)}\right|\leq \frac{\delta}{2},
$$
so that, on the event $\tau_0<n-2\delta^{-1}-1$, we are guaranteed that
 $T_{j(\tau_0)}$ is either $T_{j}$ or $T_{j+1}$ when $j$ is such that $k/n\in [T_{j},T_{j+1})$.

After time $\tau_0$, we consider the process
\begin{equation*}
Y_u = K_{\tau_0 + u} - T_{j(\tau_0)} ( n - \tau_0 - u) \quad \quad  \text{for }u \in \{0,1, \ldots n-\tau_0\},
\end{equation*}
which serves a useful vehicle to study the behavior of the budget-ratio process
and, in particular, to track the deviations of the budget ratio from the threshold $T_{j(\tau_0)}$.
This is because
$$
\left|\frac{K_{\tau_0+u}}{n-\tau_0-u}- T_{j(\tau_0)}\right| > \delta\mbox{ if and only if }
|Y_u| > \delta ( n - \tau_0 - u).
$$
In words, the ratio is outside of the dark region in Figure \ref{fig:BR-policy} if and only if the deviation process $Y_u$ exceeds the ``moving target'' $\delta (n-\tau_0-u)$.

The process $\{ Y_u: 0 \leq u \leq n-\tau_0 \}$ is adapted to the increasing sequence of $\sigma$-fields $\{ \widehat \F_u \equiv \F_{\tau_0+u}: 0 \leq u \leq n-\tau_0\}$.
Importantly, the process $\{ Y_u: 0 \leq u \leq n-\tau_0 \}$ has the mean-reversal property that we alluded to
in the description of Figure \ref{fig:BR-policy}:
if $K_{\tau_0 + u } > 0$ and $Y_u \geq 0$, then we have that $T_{j(\tau_0)} \leq K_{\tau_0 + u } / (n - \tau_0 -u)$,
and the BR policy selects \emph{all} values strictly larger than $a_{j(\tau_0)+1}$, i.e.,
$\sigmaBR_{\tau_0 + u + 1}\geq \1(X_{\tau_0 + u + 1} > a_{j(\tau_0)+1} )$ so that
$$
\E[\sigmaBR_{\tau_0 + u + 1} | \widehat \F_u] \1(K_{\tau_0 + u } > 0, Y_u \geq 0) \geq \Fbar( a_{j(\tau_0)+1} ) \1(K_{\tau_0 + u } > 0, Y_u \geq 0),
$$
and we have the negative-drift property
\begin{eqnarray}\label{eq:drift-Yupositive}
  \E[ Y_{u+1} -  Y_u | \widehat \F_u] \1(K_{\tau_0 + u } > 0, Y_u \geq 0)
        & = & \{ -\E[\sigmaBR_{\tau_0 + u + 1} | \widehat \F_u] + T_{j(\tau_0)} \}  \1(K_{\tau_0 + u } > 0, Y_u \geq 0) \\
   & \leq & - \frac{1}{2} f_{j(\tau_0)}  \1(K_{\tau_0 + u } > 0, Y_u \geq 0). \notag
\end{eqnarray}
In the case that $Y_u < 0$ and $K_{\tau_0 + u } > 0$, we have that $ K_{\tau_0 + u } / (n - \tau_0 -u) < T_{j(\tau_0)}$, so that the BR policy skips all values smaller or equal to $a_{j(\tau_0)}$, i.e., $\sigmaBR_{\tau_0 + u + 1}\leq \1(X_{\tau_0 + u + 1} > a_{j(\tau_0)} )$ so that
$$
\E[\sigmaBR_{\tau_0 + u + 1} | \widehat \F_u] \1(K_{\tau_0 + u } > 0, Y_u < 0) \leq \Fbar( a_{j(\tau_0)} ) \1(K_{\tau_0 + u } > 0, Y_u < 0),
$$
and we have a strictly positive drift
\begin{eqnarray}\label{eq:drift-Yunegative}
  \E[ Y_{u+1} -  Y_u | \widehat \F_u] \1(K_{\tau_0 + u } > 0, Y_u < 0)
        & = & \{ -\E[\sigmaBR_{\tau_0 + u + 1} | \widehat \F_u] + T_{j(\tau_0)} \}  \1(K_{\tau_0 + u } > 0, Y_u < 0) \\
   & \geq &  \frac{1}{2} f_{j(\tau_0)}  \1(K_{\tau_0 + u } > 0, Y_u < 0). \notag
\end{eqnarray}
For completeness, we also note that if $ K_{\tau_0 + u } = 0$, then the drift is simply given by
\begin{equation}\label{eq:drift-K=0}
  \E[ Y_{u+1} -  Y_u | \widehat \F_u] \1(K_{\tau_0 + u } = 0) = T_{j(\tau_0)}\1(K_{\tau_0 + u } = 0).
\end{equation}

The bounds \eqref{eq:drift-Yupositive} and \eqref{eq:drift-Yunegative} show that, regardless of whether the ratio $K_{\tau_0 + u }/(n-\tau_0-u)$ is above or below the critical threshold, the BR policy pulls it towards the threshold. This preliminary drift analysis will be useful in showing that the stopping time
\begin{equation} \label{eq:taudefin}
\tau
= \inf\left\{t > \tau_0: \left|\frac{K_t}{n-t} - T_{j(\tau_0)} \right| > \delta   \text{ or } t \geq n - 2 \delta^{-1} - 1\right\},
\end{equation}
at which the ratio exits the critical orbit (see again the right side of Figure \ref{fig:BR-policy}) is suitably large.

\begin{theorem}[BR stopping time]
Let $\epsilon = \tfrac{1}{2}\min\{ f_m, f_{m-1}, \ldots, f_1 \}$.
Then there is a constant $M \equiv M(\epsilon)$ such that for all $(n,k)\in\mathcal{T}$, the stopping time $\tau$ in \eqref{eq:taudefin}
satisfies the bound $\E[\tau]\geq n-M$. \label{thm:tau}
\end{theorem}

\begin{figure}[t!]
    \caption{\textbf{Ratio sample paths}}\label{fig:intro}

	\begin{center}
        \includegraphics[width=0.48\textwidth]{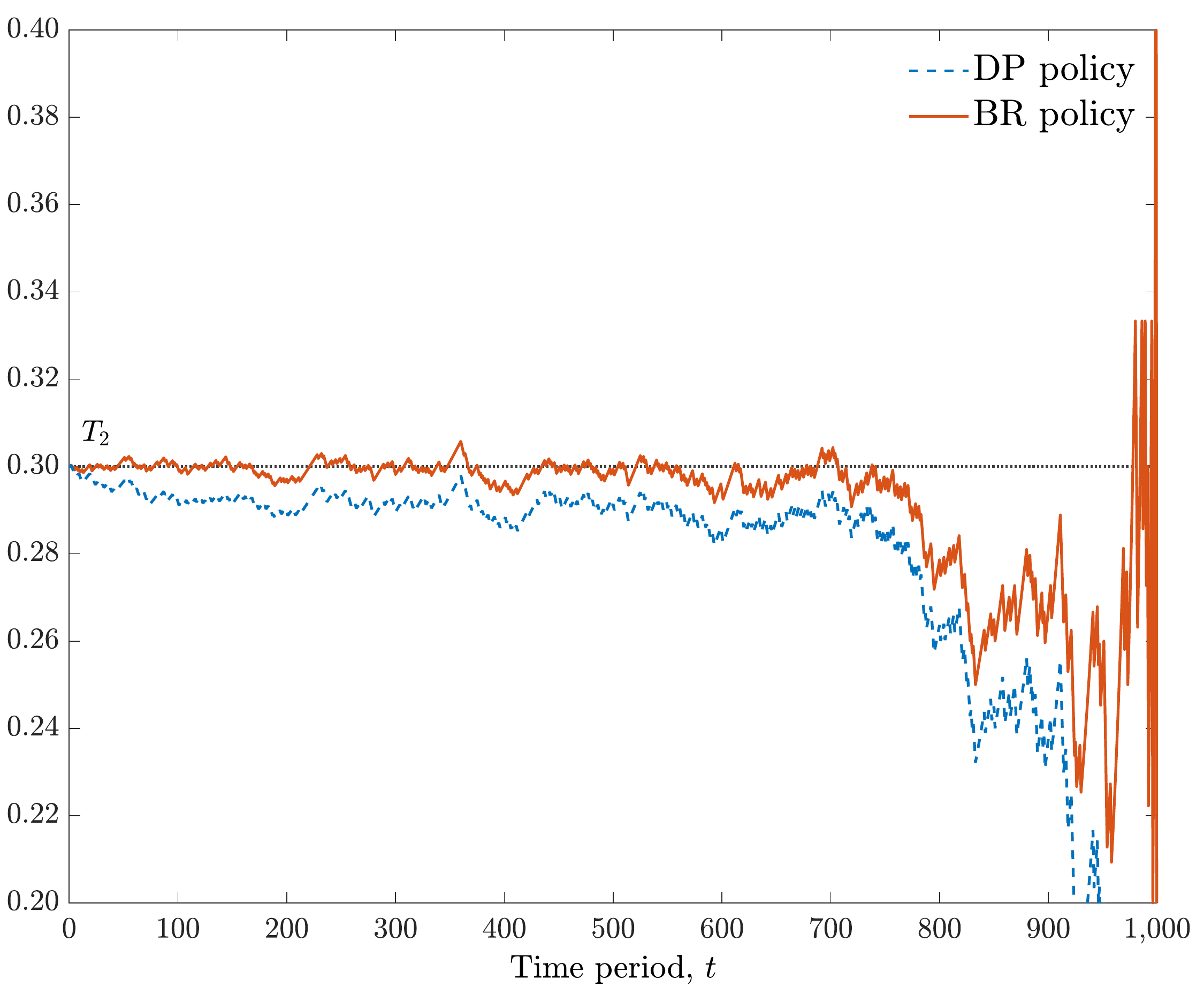}	
        \hfill
        \includegraphics[width=0.48\textwidth]{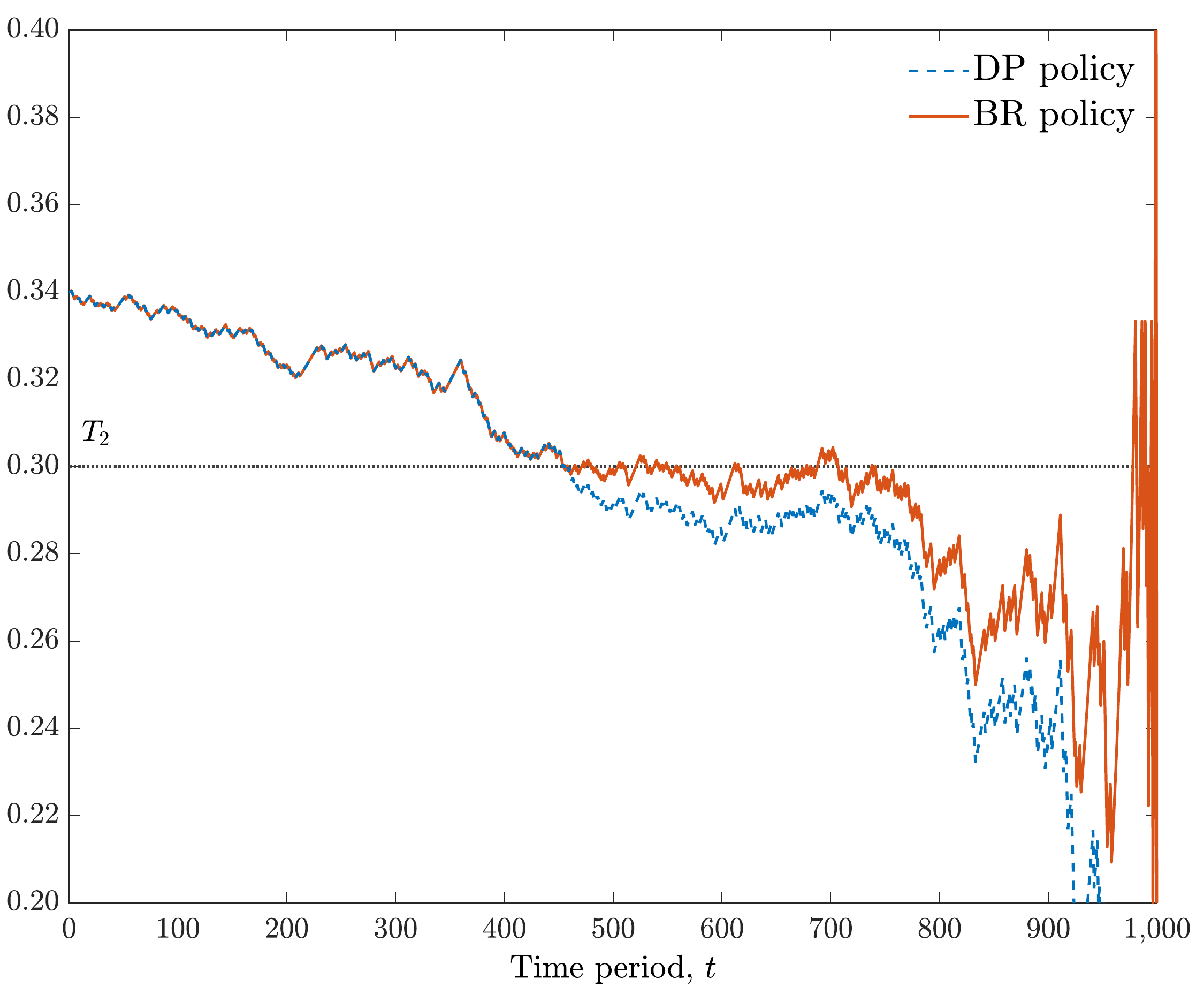}
	\end{center}

    \footnotesize{\emph{Notes.}
        Sample paths of the ratio $R^\pi_t = \tfrac{1}{n-t} \left( k - \sum_{s \in [t]} \sigma^\pi_s \right)$
        when $\pi$ is the optimal dynamic-programming (DP) policy or the Budget-Ratio (BR) policy.
        The ability distribution is uniform on $\A = \{0.20, 0.65, 1.10, 1.55, 2.00\}$.
        The two plots have common random numbers but different initial conditions.
        On the left we take $k=300$ and on the right we have $k=340$.
        In both instances, the initial ratio $k/n$ is within the orbit of $T_2 = 0.30$.
        When $k/n =0.30 = T_2$ (left) the budget ratio $R^{\rm br}_t$ stays within the orbit of $T_2$
        for most time periods. When $k/n = 0.34 > T_2$ (right) the budget ratio is first pulled towards
        the threshold $T_2$ and then it stays within its orbit until its escape time towards the end
        of the horizon. Despite just being on a sample path, the ratios $R^*_t$ and $R^{\rm br}_t$
        are remarkably close two each other.}
\end{figure}

The proof of Theorem \ref{thm:tau} (at the end of this section) is based on the mean-reversal property established above and a Lyapunov function argument. One must be careful in this analysis: when approaching the horizon's end, a small change in $K_t$ can lead to a large change in $K_t/(n-t)$. Viewed in terms of the deviation process $Y_u$, the challenge is that the target $\delta(n-\tau_0-u)$ is moving and easier to exceed when $u$ is large. Figure \ref{fig:intro} is a simulation that illustrates how this ``attraction'' to a threshold is manifested at the sample path level.

Theorem \ref{thm:tau} shows that the BR policy satisfies the sufficient condition (iv) in Proposition \ref{prop:sufficient-condition-for-policy}. Corollary \ref{cor:regret} then proceeds to show that the remaining requirements (i)--(iii) in Proposition \ref{prop:sufficient-condition-for-policy} are also satisfied.

\begin{corollary}[Uniformly bounded regret]\label{cor:regret}
Let $\epsilon = \tfrac{1}{2}\min\{ f_m, f_{m-1}, \ldots, f_1 \}$.
Then the BR policy and the stopping time $\tau$ in \eqref{eq:taudefin} satisfy the properties in Proposition \ref{prop:sufficient-condition-for-policy}.
In particular, there is a constant $M \equiv M(\epsilon)$ such that
 	$$
    \Voff^*(n,k)-\Von^*(n,k)\leq 2 a_1 M \quad \quad \mbox{ for all }(n,k)\in\mathcal{T}.
    $$	
\end{corollary}

\begin{figure}[t!]
    \caption{\textbf{Dynamic-programming and Budget-Ratio regrets}}  \label{fig:numerical1}
    \begin{center}	
        \includegraphics[width=0.48\textwidth]{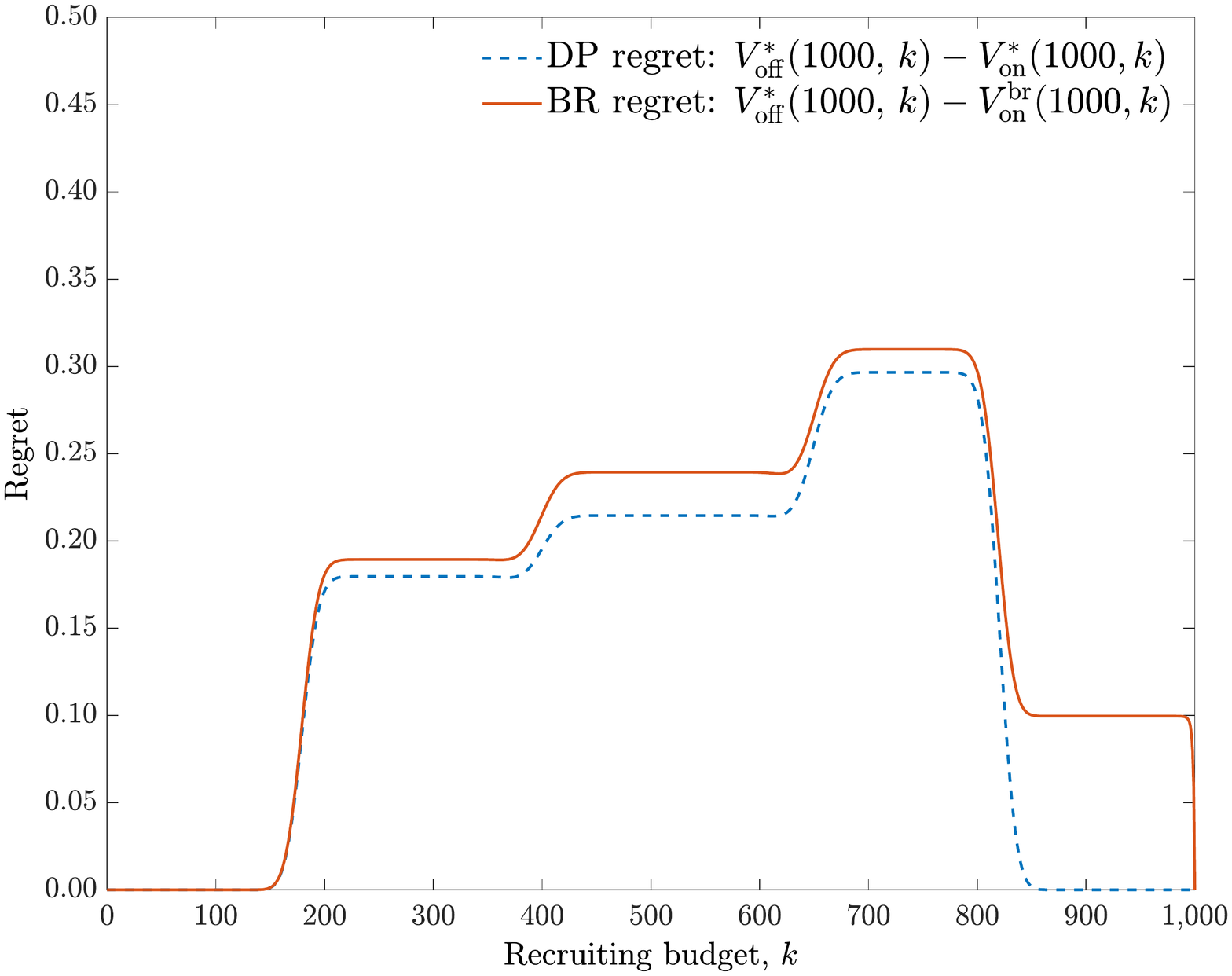}
        \hfill	
        \includegraphics[width=0.48\textwidth]{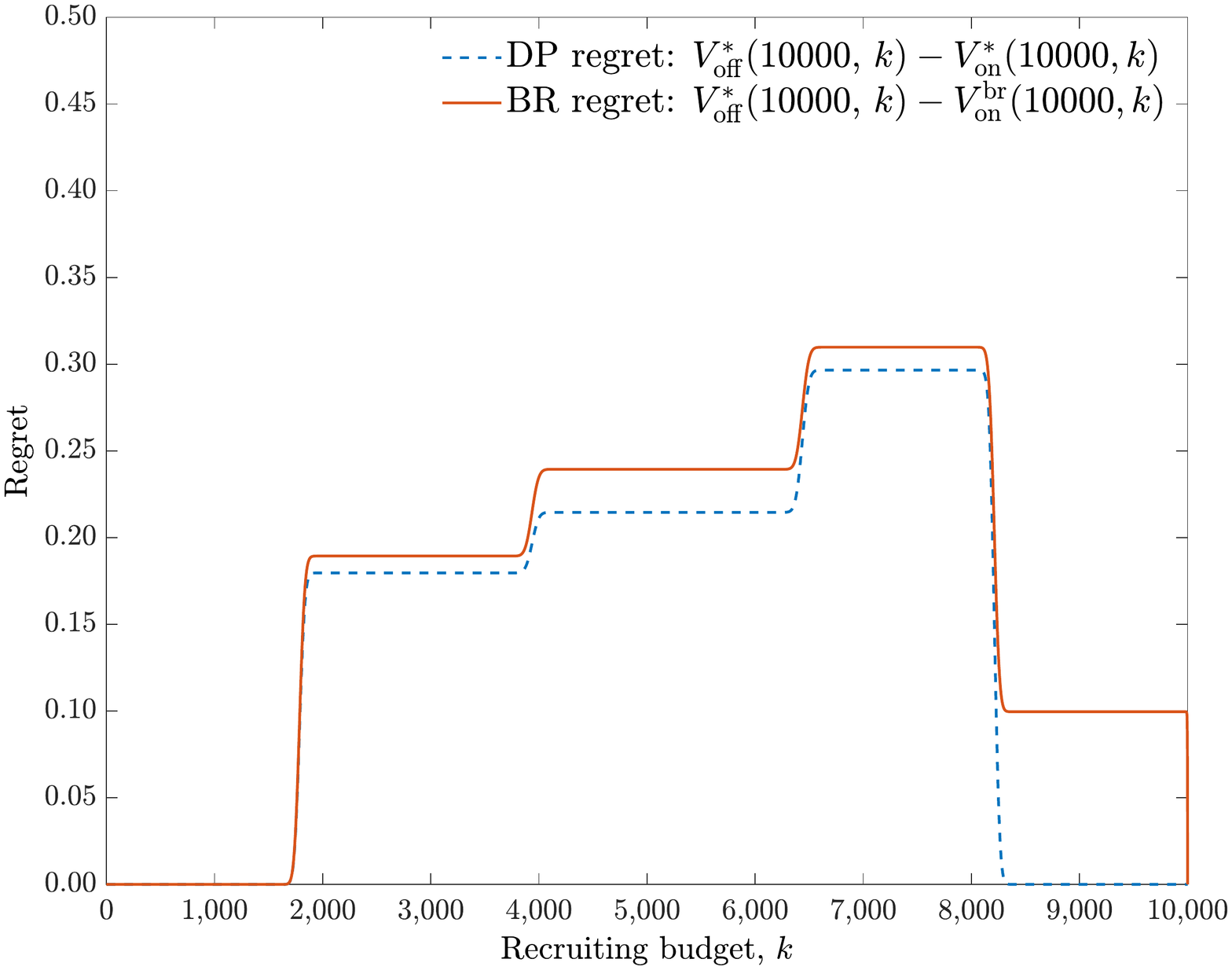}
    \end{center}

    \footnotesize{\emph{Notes.}
        The figure displays the dynamic-programming (DP) and budget-ratio (BR) regrets
        for the 5-point distribution on $\A = \{0.2, 0.5, 0.7, 0.8, 1.0\}$
        with probability mass function $(f_5, f_4, f_3, f_2, f_1) = (\tfrac{5}{28}, \tfrac{5}{28}, \tfrac{7}{28}, \tfrac{6}{28}, \tfrac{5}{28})$.
        We fix $n$ to either $n=1,000$ (left) or $n=10,000$ (right) and vary the initial budgets $k$.
        The regret is small and piecewise constant.
        The ``jumps'' in the regret function match the jump points of the distribution.}

\end{figure}

\begin{figure}[ht!]
    \caption{\textbf{Dynamic-programming and Budget-Ratio regrets for uniform distributions}}  \label{fig:numerical2}

    \begin{center}
    	\includegraphics[width=0.48\textwidth]{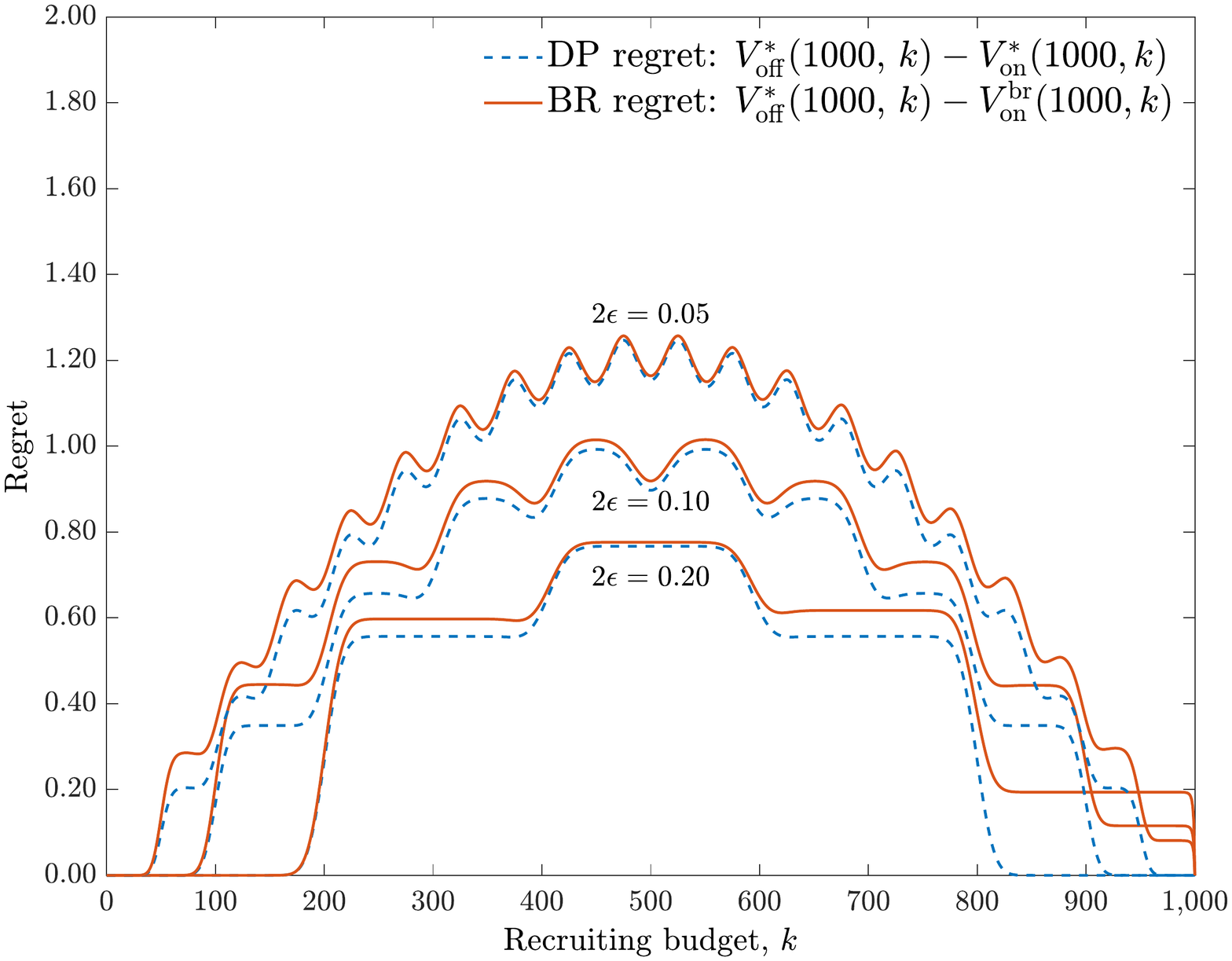}	
        \hfill
        \includegraphics[width=0.48\textwidth]{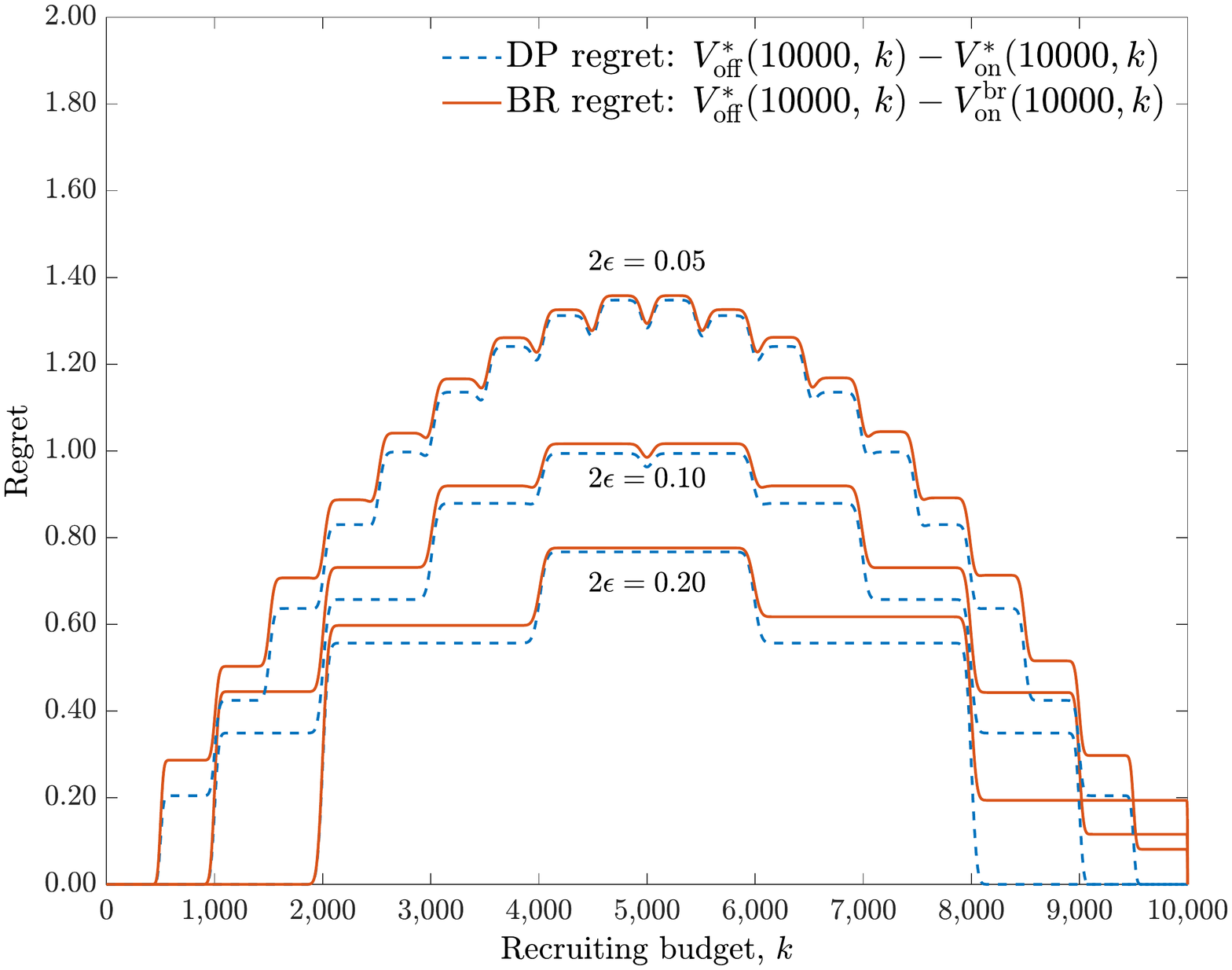}	
    \end{center}

    \footnotesize{\emph{Notes.}
    The figure displays the dynamic-programming (DP) and budget-ratio (BR) regrets
    for uniform distributions over $[0.2,0.2+\Delta,0.2+2\Delta,\ldots,2.0-\Delta,2.0]$
    where $\Delta = (2.0-0.2) 2 \epsilon/(1-2\epsilon)$ for $2 \epsilon \in \{0.05, 0.10, 0.20\}$.
    As $\epsilon$ shrinks the cardinality of the support grows.
    We fix $n=1,000$ (left) and $n=10,000$ (right) and vary the recruiting budget $k$.
    As $n$ grows, the regret approaches a piecewise constant function.
    The optimality gap of the BR policy is, regardless, very small.}

\end{figure}

The BR policy, while achieving bounded regret, is {\emph not} the optimal policy.
Considering the optimality equation \eqref{eq:g-recursion}  developed in Appendix \ref{sec:bellman},
one can see that with two periods to go ($\ell=2$ there) and one unit of budget ($k=1$)
the optimal action is to take any (and only) values $a_j\geq h_2(1)=g_1(1)-g_1(0)=\E[X_1]$.
In this state, the BR policy will instead take all values above the median.
Of course while the BR policy makes some mistakes when approaching the horizon's end,
Corollary \ref{cor:regret} is evidence that it mostly does the right thing. Figures \ref{fig:numerical1}-\ref{fig:numerical2} are numerical evidence for the performance of the BR policy.

\proof{Proof of Corollary \ref{cor:regret}.}
It suffices to prove that the BR policy satisfies the sufficient conditions stated in
Proposition \ref{prop:sufficient-condition-for-policy}.
By Theorem \ref{thm:tau} the BR policy has the associated stopping time $\tau$ in \eqref{eq:taudefin}
that satisfies condition (iv) in the proposition.
Conditions (i)--(iii) are verified by the following sample-path argument.

If $n < 2 \delta^{-1} + 1$, then the three conditions are satisfied immediately by choosing $M$
to be a suitable constant.
Otherwise if $n \geq 2 \delta^{-1} + 1$, we note that
$T_j = \Fbar(a_j) +\tfrac{1}{2}f_j = \Fbar(a_{j+1})- \tfrac{1}{2}f_{j}$ for all $j \in \{2, \ldots,m\}$,
so the definition \eqref{eq:j0} tells us that
if $k/n \in [T_{j},T_{j+1})$ then $j_0(n,k) = j$ for all $j \in [m]$.
Thus, if $j$ is the index such that $k/n\in [T_{j},T_{j+1})$,
then $j$ is the ``action'' index identified in Proposition \ref{prop:sufficient-condition-for-policy},
and we verify the three conditions
by distinguishing the case $\{\tau_0<n - 2\delta^{-1} - 1\}$ from $\{\tau_0 = n - 2\delta^{-1} - 1\}$.

As argued earlier,
on the event $\{\tau_0<n - 2\delta^{-1} - 1\}$ the index $j(\tau_0)\in\{j,j+1\}$.
In turn, for all $t< \tau_0$, the BR policy selects $X_{t+1}\geq a_j$ and skips all $X_{t+1} \leq a_{j+1}$.
If $j(\tau_0)=j$ then for time indices $t\in [\tau_0,\tau)$ all values $a_{j-1}$ and greater are selected
and all values $a_{j+1}$ or smaller are skipped. If $j(\tau_0)=j+1$, then on $t\in[\tau_0,\tau)$ all values $a_j$ and greater are selected
and all those smaller or equal than $a_{j+2}$ are skipped.
Thus, on the event $\{\tau_0<n - 2\delta^{-1} - 1 \}$ we have
\begin{equation}\label{eq:interim1}
\sum_{i\in [j - 1]}S_i^{\BR,\tau}=\sum_{i\in [j -1 ]}Z_i^{\tau},
\end{equation}
where, recall, $S_i^{\BR,\tau}$ is the number of $a_i$ candidates selected by time $\tau$.
Furthermore,
\begin{enumerate}[(i)]
    \item
        If $j(\tau_0)=j$, then
        \begin{equation}\label{eq:j(tau0)=j}
         S_j^{\BR, \tau} \!=\!  \min \{ Z_{j}^{\tau}, k-K_{\tau}-\sum_{i\in [j-1]} Z_i^{\tau} \}
                         \!=\!  \min \{ Z_{j}^{\tau}, (k-\sum_{i\in [j-1]} Z_i^{\tau})_+  -K_{\tau} \}
         \quad  \mbox{ and } \quad
         S_{j+1}^{\BR, \tau}\!=\!0.
        \end{equation}
        Here, the left equality holds because out of the total budget used by time $\tau$,
        which is given by $k-K_{\tau} = \sum_{t \in [\tau]} \sigmaBR_t$,
        the quantity $\sum_{i\in[j-1]} Z_i^{\tau}$ is allocated to values larger or equal to $a_{j-1}$
        and the remaining to $a_{j}$.
	\item
        If $j(\tau_0)=j+1$, then we have that
        \begin{equation}\label{eq:j(tau0)=j+1}
        S_j^{\BR, \tau}\!=\!Z_j^{\tau}
        \quad \mbox{ and } \quad
        S_{j+1}^{\BR, \tau}\!=\!\min \{Z_{j+1}^{\tau}, k - K_{\tau}-\sum_{i \in [j]} Z_i^{\tau}\}
                           \!=\!\min \{Z_{j+1}^{\tau}, (k -\sum_{i \in [j]} Z_i^{\tau})_+ - K_{\tau}\},
        \end{equation}
        where the right equality holds for reasons that are similar to those given just below \eqref{eq:j(tau0)=j}.
\end{enumerate}
By combining the observations in \eqref{eq:j(tau0)=j} and \eqref{eq:j(tau0)=j+1},
we find on the event $\{\tau_0<n - 2 \delta^{-1}-1\}$ that
\begin{equation} \label{eq:interim1a}
S_{j}^{\BR, \tau}\geq \min\{ (k-\sum_{i\in [j-1]} Z_i^{\tau})_+, Z_{j}^{\tau} \} - K_{\tau}
                   = \mathfrak{S}_{j}^\tau  -K_{\tau},
\end{equation}
and
\begin{equation} \label{eq:interim1b}
S_{j+1}^{\BR, \tau}\geq \min\{ (k-\sum_{i\in [j]} Z_i^{\tau})_+ , Z_{j+1}^{\tau} \} -K_{\tau}
                     = \mathfrak{S}_{j+1}^\tau  -K_{\tau},
\end{equation}
where we use the fact that if $x,y,z$ are non-negative numbers then $\min(x-y,z)\geq \min(x,z)-y$.

On the event $\{\tau_0=n - 2\delta^{-1} -1\}$ we have that $\tau = \tau_0$ and all values greater than or equal to $a_j$
are selected and all values lower than or equal to $a_{j+1}$ are skipped so that
\begin{equation}\label{eq:interim2a}
\sum_{i\in [j-1]}S_i^{\BR, \tau}=\sum_{i\in [j-1]}Z_i^{\tau}
\quad \quad \text{and}\quad \quad
S_j^{\BR, \tau} = \mathfrak{S}_j^\tau,
\end{equation}
and
\begin{equation} \label{eq:interim2b}
    0 = S_{j+1}^{\BR, \tau}=
  (k-\sum_{i\in [j]} Z_i^{\tau})_+ - K_{\tau} \geq
  \min\{Z_{j+1}^{\tau}, (k-\sum_{i\in [j]} Z_i^{\tau})_+\} -K_{\tau}
   = \mathfrak{S}_{j+1}^\tau -K_{\tau}.
\end{equation}

Finally, since the BR policy selects all remaining values as soon as there is a $t \in [n]$ such that $K_t \geq n-t$,
we have that $K_\tau \leq n - \tau$ and, consequently, that
\begin{equation}\label{eq:Ktau-upper-bound}
-\E[K_{\tau}]\geq \E[\tau]-n\geq -M
\end{equation}
where the last inequality follows from Theorem \ref{thm:tau}.

If we now recall the estimates \eqref{eq:interim1}, \eqref{eq:interim1a} and \eqref{eq:interim1b}
which hold on the event  $\{\tau_0 < n - 2\delta^{-1}-1\}$
and the relations \eqref{eq:interim2a} and \eqref{eq:interim2b} which are satisfied on $\{\tau_0 = n - 2\delta^{-1} - 1\}$,
take expectations and recall the bound \eqref{eq:Ktau-upper-bound},
we see that the sufficient conditions (i)--(iii) in Proposition \ref{prop:sufficient-condition-for-policy}
are all satisfied.
\halmos\endproof\vspace{0.25cm}

\subsection{An alternative to Budget-Ratio: the Adaptive-Index policy.}\label{se:adaptive-index}

\begin{figure}[!ht]
    \caption{\textbf{Dynamic-programming (DP), budget-ratio (BR), and adaptive-index (AI) regrets}}
    \label{fig:gap}

	\begin{center}
        \includegraphics[width=0.5\textwidth]{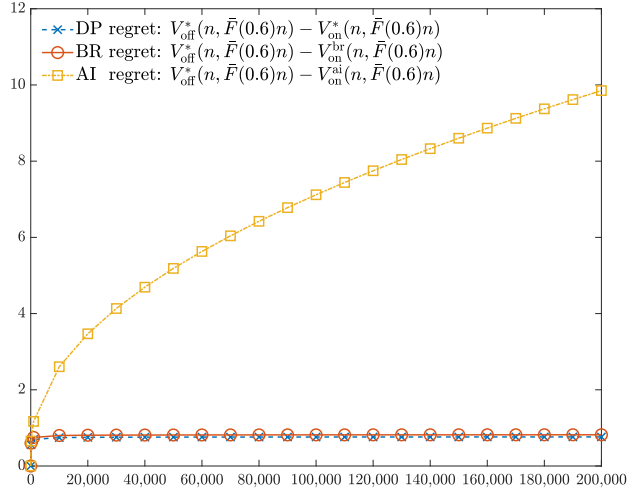}
    \end{center}	
		
    \footnotesize{\emph{Notes.}
            The graph displays the regret of different policies when the ability distribution is uniform on
            $\A = \{0.2,0.4,0.6,\ldots,1,1.2,1.4,\ldots,2.0\}$ as a function of the horizon length $n$.
            In the chart we vary $n$ from $0$ to $200,000$ and keep the initial budget at $\frac{k}{n}= \Fbar(0.6)$
            (a mass point of the distribution).
            Whereas the optimal and Budget-Ratio policies achieve bounded regret,
            the regret of the Adaptive-Index policy grows with the problem size $n$.}
\end{figure}

In closely related work, \cite{wu2015algorithms} offer an elegant adaptive-index policy that we revisit here.
Given the deterministic relaxation \eqref{eq:DR}, we re-solve in each time period the deterministic problem
$$
DR(n-t,K_t)= \varphi(\E[Z_1^{n-t}], \ldots, \E[Z_m^{n-t}], K_t),
$$
and recall from \eqref{eq:DR-optimal-solution} that the optimal solution is given by
$
s^*_{j,t}= \min \{ \E[Z_j^{n-t}], (K_t-\sum_{i\in[j-1]}\E[Z_i^{n-t}])_+ \}
     = \min \{ f_j(n-t), (K_t-\Fbar(a_j)(n-t))_+ \}
$
for all $j\in [m]$.

Then, we construct the adaptive (re-optimized) index policy by mimicking the solution of this optimization problem.
Specifically, if $s_{j,t}^*=f_j (n-t)$ and the candidate inspected at time $t+1$ has ability $a_j$
then that candidate is selected.
Otherwise, if $s_{j,t}^* = K_t-\Fbar(a_j)(n-t) > 0$ then an arriving $a_j$ candidate
is selected with probability
$$
p_{j',t+1}=\frac{K_t/(n-t)- \Fbar(a_{j'})}{f_{j'}}.
$$
Finally, if $s_{j,t}^* = 0$ then an arriving $a_j$ candidate is rejected.

This policy induces a nice martingale structure.
Using the notation introduced in Section \ref{sec:model},
we let $\sigma^{\rm ai}_1, \sigma^{\rm ai}_2, \ldots, \sigma^{\rm ai}_n$ be the sequence
of decisions of the Adaptive-Index policy and let
$K_t^{\rm ai} = k-\sum_{s\in[t]} \sigma^{\rm ai}_s$ be the associated remaining-budget process
for all $t \in [n]$ and with $K_0^{\rm ai}=k$.
In the statement and in the proof of the next proposition,
we use the standard notation $a \wedge b = \min\{a,b\}$.

\begin{proposition} \label{prop:resolve} Let $\tau_1 = \inf\{t \in [n]: K_t^{\rm ai}/(n-t)>1\}$,
then the stopped Adaptive-Index ratio process
$
\{ R^{\rm ai}_{t \wedge \tau_1} = \frac{K_{t\wedge \tau_1}^{\rm ai}}{n-(t \wedge \tau_1)}:~ t \in [n] \}
$
is a martingale.
\end{proposition}

\proof{Proof.}
Since $0\leq R^{\rm ai}_{t\wedge \tau_1}\leq 1$ it is trivially true that $\E[|R^{\rm ai}_{t\wedge \tau_1}|]<\infty$ for all $t\in [n]$.
Next, if $\F_t$ is the $\sigma$-field generated by the random variables $\sigma_1^{\rm ai},\ldots,\sigma_t^{\rm ai}$,
then we have that
$$
\E[R^{\rm ai}_{(t+1)\wedge \tau_1}-R^{\rm ai}_{t\wedge \tau_1}|\F_{t}]
 = \E[(R^{\rm ai}_{(t+1)\wedge \tau_1}-R^{\rm ai}_{t\wedge \tau_1})\1(\tau_1>t)|\F_{t}]
 = \1(\tau_1>t)\E\left[\frac{K_t^{\rm ai}-\sigma^{\rm ai}_{t+1}}{n-(t+1)}-\frac{K_t^{\rm ai}}{n-t}\Big|\F_{t}\right],
$$
where $\sigma^{\rm ai}_{t+1}=1$ with probability $\frac{K_t^{\rm ai}}{n-t}\wedge 1$ and it is zero otherwise.
Hence,
\begin{align*}
\E[R^{\rm ai}_{(t+1)\wedge \tau_1}-R^{\rm ai}_{t\wedge \tau_1}|\F_{t}]
=\left(\frac{K_t^{\rm ai}-\frac{K_t^{\rm ai}}{n-t}\wedge 1}{n-(t+1)}-\frac{K_t^{\rm ai}}{n-t}\right)\1(\tau_1>t)= 0,
\end{align*}
where we use the fact that $\1(\tau_1>t)$ implies $K_t^{\rm ai}/(n-t)\leq 1$.
\halmos\endproof\vspace{0.25cm}

Because of this martingale structure, the Adaptive-Index ratio $K_t^{\rm ai}/(n-t)$ remains ``close'' to $k/n$ so that this policy,
like the BR policy, is careful in utilizing its budget and does not run out of it until (almost) the horizon's end.
Furthermore, this martingale property guarantees bounded regret when the initial ratio $k/n$
is safely far from the masses of the discrete distribution \citep[see also][and Section \ref{sec:det_relax} in this paper]{wu2015algorithms}.

In general, however, spending the budget at the right ``rate'' is not sufficient, however.
It is also important that the budget is spent on the right candidates
and the symmetric martingale structure is too weak for that purpose: when initialized, for example, at $k/n=\Fbar(a_j)$ for some $j$, the martingale spends an equal amount of time below and above $\Fbar(a_j)$,
and the re-optimized index policy takes the values  $a_j$ and $a_{j+1}$ in equal proportions.
A good policy should start selecting $a_{j+1}$ values only after it has selected all (or most) $a_j$ values first.
Figure \ref{fig:gap} gives an illustration of this phenomenon.

\subsection{Proof of Theorem \ref{thm:tau}}

Given the deviation process
$$
Y_u = K_{\tau_0 + u} - T_{j(\tau_0)} ( n - \tau_0 - u), \quad \quad  u \in \{0,1, \ldots n-\tau_0\},
$$
the key step in the proof of Theorem \ref{thm:tau} is the derivation of an exponential tail bound for the random variable  $\abs{Y_u}$
for each $u \in \{0,1,\ldots, n-\tau_0\}$.
As before, we take $\epsilon = \tfrac{1}{2}\min\{ f_m, f_{m-1},\ldots,f_1\}$.

\begin{proposition}[Exponential tail bound]\label{pr:exponential-tail-bound-Yu}
Fix $0 < \delta < \epsilon$, $c = e^2 -3$ and $0 < \eta < (\epsilon - \delta) / c$.
Then there is a constant $M \equiv M(\epsilon)$
such that, for all $0 \leq u \leq n-\tau_0$, we have the exponential tail bound
$$
\P (\abs{ Y_u } \geq \delta  ( n - \tau_0 -u) \,|\, \widehat\F_0 )
\leq  2 \exp\{ - \tfrac{\eta \delta}{2} (n-\tau_0) \} + M \exp\{ - \eta \delta (n-\tau_0-u) \}.
$$
\end{proposition}

\proof{Proof.}
If $\tau_0 = n - 2 \delta^{-1} - 1$, then the statement is trivial. Otherwise, for $\tau_0 < n - 2 \delta^{-1} - 1$
the proof is an application---using the mean-reversal property of the BR policy---of the tail bound of \citet{hajek1982}
to the two  processes $\{Y_u: 0 \leq u \leq n-\tau_0\}$ and $\{-Y_u: 0 \leq u \leq n-\tau_0\}$.

We begin with the former. For any $a\geq 0$ we have that
$\1(Y_u > a ) = \1(Y_u > a , K_{\tau_0 + u } > 0 )$
for all $0 \leq u \leq n-\tau_0$,
so that, by replacing $f_{j(\tau_0)}$ with $\epsilon = \tfrac{1}{2} \min\{ f_m, f_{m-1}, \ldots, f_1\}$ in \eqref{eq:drift-Yupositive},
we obtain the condition
\begin{equation}\label{eq:C1}\tag{C1}
\E[ Y_{u+1} - Y_u | \widehat \F_u] \1(Y_u > a ) \leq - \epsilon \1(Y_u > a)
\quad \quad \text{for all $a \geq 0$ and $0 \leq u \leq n-\tau_0$.}
\end{equation}
Since, by definition,
$$
\abs{ Y_{u+1} - Y_u } = \abs{ -\sigmaBR_{\tau_0 + u + 1} - T_{j(\tau_0)} } \leq 2
\quad \quad \text{for all  $0 \leq u < n-\tau_0$,}
$$
we also have the condition
\begin{equation}\label{eq:C2}\tag{C2}
\E[ \exp\{  \lambda \abs{ Y_{u+1} - Y_u } \} | \widehat\F_u]  \leq e^{2 \lambda }
\quad \quad \text{for all } \lambda > 0.
\end{equation}
Conditions \eqref{eq:C1} and \eqref{eq:C2} (with $\lambda=1)$ imply, by \citet[Lemma 2.1]{hajek1982}
that for $\eta$, $\rho$ satisfying
\begin{equation} \label{eq:setting-eta-rho}
0 <  \eta < (\epsilon - \delta)/c
\quad \quad \text{and} \quad \quad
\rho = 1 - \eta ( \epsilon - \eta c),
\end{equation}
we have
$$
\E[ \exp\{ \eta ( Y_{u+1} - Y_u ) \} | \widehat\F_u] \1(Y_u > a)   \leq  \rho
\quad\quad \text{and} \quad\quad
  \E[ \exp\{ \eta ( Y_{u+1} - a ) \} | \widehat\F_u] \1(Y_u \leq a)  \leq  e^2,
$$
for all $a \geq 0$ and all $0 \leq u < n-\tau_0$. \citet[Theorem 2.3, Equation 2.8]{hajek1982} then gives
$$
\P( Y_u \geq \delta (n - \tau_0 - u) | \widehat \F_0)
\leq \rho^u \exp\{ \eta Y_0 - \eta \delta(n - \tau_0 - u) \}
     + \tfrac{e^2}{1-\rho} \exp\{ \eta a - \eta \delta( n - \tau_0 - u) \},
$$
for all $a \geq 0$ and $0 \leq u \leq n-\tau_0$. Taking $a = 0$ and using the fact that
and $ Y_0 \leq \abs{ Y_0 } = \abs{ K_{\tau_0} - T_{j(\tau_0)} (n - \tau_0) } \leq \tfrac{\delta}{2} (n - \tau_0)$ as well as the standard inequality
$\rho = 1 - \eta( \epsilon - \eta c) \leq \exp\{- \eta( \epsilon - \eta c)  \}$,
we finally have the upper bound
$$
\P( Y_u \geq \delta (n - \tau_0 - u) | \widehat \F_0)
\leq \exp\{ - \tfrac{\eta \delta }{2} (n - \tau_0) - \eta u (\epsilon - \eta c - \delta) \}
     + \tfrac{e^2}{1-\rho} \exp\{ - \eta \delta( n - \tau_0 - u) \}.
$$

The choice of $\eta$ in \eqref{eq:setting-eta-rho} tells us that $\epsilon-\eta c-\delta > 0$,
so we can drop the second term in the first exponent on the right-hand side.
By setting $M \equiv M(\epsilon)=\frac{e^2}{1-\rho}$, we then have
\begin{equation}\label{eq:Yu-tailbound}
\P( Y_u \geq \delta (n - \tau_0 - u) | \widehat \F_0)
\leq \exp\{ - \tfrac{\eta \delta }{2} (n - \tau_0)\} + M \exp\{ - \eta \delta( n - \tau_0 - u) \}.
\end{equation}

The analysis of the sequence $\{-Y_u : 0 \leq u \leq n-\tau_0\}$ follows a similar logic. For any $a \geq 0$
\begin{align*}
\E[ - Y_{u+1} + Y_u | \widehat \F_u] \1( - Y_u > a )
= & - \E[ Y_{u+1} - Y_u | \widehat \F_u] \1( Y_u < - a, K_{\tau_0 + u} > 0 )  \\
  & - \E[ Y_{u+1} - Y_u | \widehat \F_u] \1( Y_u < - a, K_{\tau_0 + u} = 0 ),
\end{align*} so that by \eqref{eq:drift-Yunegative} and \eqref{eq:drift-K=0}
we have
$$
\E[ - Y_{u+1} + Y_u | \widehat \F_u] \1( - Y_u > a )
 \leq - \tfrac{1}{2}f_{j(\tau_0)} \1( Y_u < - a, K_{\tau_0 + u} > 0 )
      - T_{j(\tau_0)} \1( Y_u < - a, K_{\tau_0 + u} = 0 ),
$$
for all $a\geq 0$. On the event $\{ Y_u < -a, K_{\tau_0 + u} = 0 \}$,
we must have that $j(\tau_0)>1$. Otherwise, if $j(\tau_0)=1$, $T_{j(\tau_0)} = 0$ and $Y_u=K_{\tau+u}=0\geq -a$ by definition.
In particular, $T_{j(\tau_0)} \geq \epsilon$ on this event, and it follows that
\begin{equation}\label{eq:C1tilde}\tag{$\widetilde{\rm C1}$}
\E[ - Y_{u+1} + Y_u | \widehat \F_u] \1( - Y_u > a )
 \leq - \epsilon \1( - Y_u >  a ).
\end{equation}
for all $a \geq 0$. It is then easily verified that
\begin{equation}\label{eq:C2tilde}\tag{$\widetilde{\rm C2}$}
\E[ \exp\{  \lambda \abs{ Y_{u+1} - Y_u } \} | \widehat\F_u]  \leq e^{2 \lambda }
\quad \quad \text{for all } \lambda > 0.
\end{equation}
As before,
\citet[Lemma 2.1 and Theorem 2.3]{hajek1982} gives---with
$c= e^2 -3$, $0 < \eta \leq (\epsilon - \delta)/c$, $\rho = 1 - \eta(\epsilon - \eta c)$,
and $M \equiv M(\epsilon) = \frac{e^2}{1-\rho}$---that
\begin{equation}\label{eq:Yu-negative-tailbound}
\P( - Y_u \geq \delta (n - \tau_0 - u) | \widehat \F_0)
\leq \exp\{ - \tfrac{\eta \delta }{2} (n - \tau_0)\} + M \exp\{ - \eta \delta( n - \tau_0 - u) \}.
\end{equation}
The statement of the lemma is now the combination of \eqref{eq:Yu-tailbound} and \eqref{eq:Yu-negative-tailbound}
with $M \equiv M(\epsilon) = \frac{2 e^2}{1-\rho}$.
\halmos\endproof\vspace{0.25cm}

The exponential tail bound in Proposition \ref{pr:exponential-tail-bound-Yu}
goes a long way for the proof of Theorem \ref{thm:tau} which follows next.

\proof{Proof of Theorem \ref{thm:tau}.}
By definition, for $\tau_0 < t \leq n - 2 \delta^{-1} - 1$, we have that
$$
\tau\leq t \mbox{ if and only if } \left| \frac{K_{\tau_0+u}}{n-\tau_0-u} - T_{j(\tau_0)} \right| > \delta,
\mbox{ for some } 0\leq u\leq t-\tau_0,
$$
or, represented in terms of $Y_u$,
$$
\tau \leq t
\quad \quad \text{if and only if} \quad \quad
\sum_{u\in[t-\tau_0]} \1( \abs{Y_u} > \delta (n - \tau_0 - u) ) \geq 1.
$$
We will represent $\E[\tau]$ as a sum of the tail probabilities
$$
\P( \tau \leq t | \widehat \F_0 )
= \P \big( \sum_{u\in[t-\tau_0]} \1( \abs{Y_u} > \delta (n - \tau_0 - u) ) \geq 1 | \widehat \F_0 \big),
\quad \quad \text{for } \tau_0 < t \leq n - 2 \delta^{-1} -1,
$$
which, by Markov's inequality, satisfy the bounds
$$
\P( \tau \leq t | \widehat \F_0 )
\leq \sum_{u \in [t-\tau_0]} \P( \abs{Y_u} > \delta (n - \tau_0 - u) | \widehat \F_0 ),
\quad \quad \text{for } \tau_0 < t \leq n - 2 \delta^{-1} -1.
$$

By integrating the exponential tail bound in Proposition \ref{pr:exponential-tail-bound-Yu}
for $0 \leq u \leq t-\tau_0$ (recall that $\tau_0 < t \leq n - 2 \delta^{-1} - 1 < n$),
we obtain
\begin{equation}\label{eq:P(tau-leq-t)}
  \P( \tau \leq t | \widehat \F_0 )
  \leq 2 (n - \tau_0) \exp\{ - \tfrac{\eta \delta}{2} (n - \tau_0) \} + M \exp\{ - \eta \delta (n - t)\},
\end{equation}
for some constant $M\equiv M(\epsilon)$.
Since $\tau>\tau_0$ by definition, then  $\P(\tau>t|\widehat \F_0)=1$
for all $t\leq \tau_0$ and we also have that
$$
\E[\tau | \widehat \F_0]
 \geq \sum_{t=0}^{\tau_0} \P( \tau > t | \widehat \F_0 )
   + \sum_{t=\tau_0+1}^{n - 2 /\delta - 1} \P( \tau > t | \widehat \F_0 )
   \geq n -  2 \delta^{-1} - 1  - \sum_{t=\tau_0+1}^{n - 2/ \delta - 1 } \P( \tau \leq t | \widehat \F_0 ).
$$
Integrating \eqref{eq:P(tau-leq-t)} then gives us another constant $M \equiv M(\epsilon)$
such that
$$
\E[\tau | \widehat \F_0]
 \geq n -  2 (n - \tau_0)^2 \exp\{ - \tfrac{\eta \delta}{2} (n - \tau_0) \} - M.
$$
The middle summand on the right-hand side is uniformly bounded, so in summary
we have that
$$
\E[\tau | \widehat \F_0] \geq n - M,
$$
for some constant $M \equiv  M(\epsilon)<\infty$,
and the proof of the theorem follows after one takes total expectations.
\halmos\endproof\vspace{0.25cm}

\section{The square-root regret of non-adaptive policies\label{sec:nonadaptive}}
		
A non-adaptive policy $\pi$ is an online feasible policy that is characterized
by a probability matrix $\{ p_{j,t} :  j \in [m] \text{ and } t \in [n]\}$.
The entry $p_{j,t}$ represents the probability of selecting the candidate
inspected at time $t$ given that the candidate's ability is $X_t = a_j$
and that the recruiting budget remaining at time $t$ is non zero.
Formally, given $\pi = \{ p_{j,t} :  j \in [m] \text{ and } t \in [n]\}$,
let
$$
B_t = \sum_{j \in [m]} \1( U_t \leq p_{j,t}, X_t = a_j)
\quad \quad \text{ for } t \in [n],
$$
and note that the sequence $B_1, B_2, \ldots, B_n$ is a sequence of $\F_t$-measurable, independent
Bernoulli random variables with success probabilities $q_1, q_2, \ldots, q_n$ such that
\begin{equation*}
    \E[B_t | X_t = a_j] = p_{j,t}
	\quad \quad \text{ and } \quad \quad
	q_t = \E[B_t] = \sum_{j\in [m]} p_{j,t} f_j.
\end{equation*}
The non-adaptive policy $\pi$ selects candidates until it runs out of budget
or reaches the end of the horizon, i.e., up to the stopping time $\nu \equiv \nu(\pi)$
given by
\begin{equation*}
	\nu = \min\big\{ r \geq 1: \sum_{t\in [r]} B_t \geq k \text{ or } r \geq n\big\}.
\end{equation*}
The expected total ability accrued by the non-adaptive policy $\pi$
is then given by
\begin{equation*}
	\Von^\pi (n,k) = \sum_{j\in [m]} a_j \E \left[\sum_{t\in [\nu]} B_t \1(X_t = a_j) \right],
\end{equation*}
and
$
\Vna^{*}(n,k)=\sup_{\pi \in \Pi_{\na}} \Von^{\pi}(n,k)
$
is the performance of the best non-adaptive policy.
		
An intuitive non-adaptive policy is the
\emph{index} policy
$\id  = \{ \p_{j,t} :  j \in [m] \text{ and } t \in [n]\}$,
that takes its probabilities from the solution \eqref{eq:DR-optimal-solution}
to the deterministic relaxation \eqref{eq:DR}.
If the residual budget at time $t$ is non-zero
and if and $j_{\id} \in [m]$ is the index such that $\Fbar(a_{j_{\id}}) \leq k/n < \Fbar(a_{j_{\id}+1})$,
then the index policy selects an arriving $a_j$-candidate with probability
\begin{align}
	\p_{j,t} = \frac{s^*_j}{\E[Z^n_j]}=\frac{s^*_j}{nf_j}= \begin{cases}
		1 & \mbox{if } j \leq j_{\id} -1, \\
		\frac{k/n - \Fbar(a_{j_{\id}})}{f_{j_{\id}}} & \mbox{if } j = j_{\id}, \\
		0 & \mbox{if } j \geq j_{\id}+1.
    \end{cases} \label{eq:index}
\end{align}
		
In the main result of this section we prove that, for a large range of $(n,k)$ pairs,
the regret of non-adaptive policies is generally of the order of $\sqrt{n}$.
Viewed in the context of Proposition \ref{prop:sufficient-condition-for-policy},
non-adaptive policies violate property (iv):
they run out of budget too early.
The proof relies on this ``greediness.'' 		
		
To build intuition, consider the problem instance
with ability distribution that is uniform on $\A = \{a_3, a_2, a_1\}$
and with $n = 2 k$. For each $t \in [n]$,
the index policy \emph{would} select $a_1$-candidates with probability $\p_{1,t} = 1$,
$a_2$-candidates with probability $\p_{2,t} = 1/2$, and no $a_3$ candidates.
The number of selections that the index policy would make by time $t$ is
then a Binomial random variables with $t$ trials and success probability equal to $k/n = 1/2$.
For $t = \lfloor n - \sqrt{n} \rfloor$ this random variable
has mean approximately equal to $\frac{1}{2}( n-\sqrt{n} ) = (1-\frac{1}{\sqrt{n}})\frac{n}{2}$
and standard deviation approximately equal to $\frac{1}{4}\sqrt{n-\sqrt{n}}$.
In this case, the recruiting budget $k = n/2$ is, in the limit, 2 standard deviations from the mean so that,
with non-negligible probability, the policy would have exhausted all its budget by time $n-\sqrt{n}$
and missed $\Theta(\sqrt{n})$\footnote{We say that a function $\phi(n) = \Theta(\psi(n))$ if
there are constants $c_1$, $c_2$, and $n_0$ such that $c_1 \leq \frac{\phi(n)}{\psi(n)} \leq c_2$ for all $n \geq n_0$.}
opportunities to select $a_1$ values.
The index policy, however, did select many $a_2$ values up to this time $t$. The Budget-Ratio policy---and the dynamic programming policy---would exchange those
$a_2$ with $a_1$ and attain a regret that does not grow with $n$.

\begin{theorem}[The regret of non-adaptive policies] \label{thm:nonadaptive}
	For $\epsilon = \tfrac{1}{2}\min\{f_m, f_{m-1}, \ldots, f_1\}$ suppose that
	$(f_1 + \epsilon) n  \leq k \leq (1-f_m - \epsilon) n$.
	Then there is a constant $M \equiv M(\epsilon, m, a_1, \ldots, a_m)$
	such that
	$$
	M \sqrt{n }\leq \Voff^{*}(n,k)-\Vna^{*}(n,k).
	$$
\end{theorem}
		
As one might expect, the non-adaptive (time-homogeneous) index policy in \eqref{eq:index}
already achieves this order of magnitude and, in this sense,
is representative of the performance of general non-adaptive policies.
		
\begin{lemma}[The regret of the non-adaptive index policy] \label{lem:index}
	The non-adaptive index policy $\id  = \{ \p_{j,t} :  j \in [m] \text{ and } t \in [n]\}$
	has a $\sqrt{n}$ regret.
	That is, for any $\varepsilon\in (0,1)$ and all pairs $(n,k)$ such that $\varepsilon \leq k/n$ we have
	$$
	\Voff^{*}(n,k)-\Vna^{\rm id}(n,k) \leq DR(n,k)- \Vna^{\id}(n,k) \leq \varepsilon^{-1} a_1 \sqrt{n}.
	$$
\end{lemma}
		
The conditions of Theorem \ref{thm:nonadaptive} are not necessary but certain budget ranges do have to be excluded.
In the small-budget range with $k \leq (f_1 - \epsilon)n$, for example, the regret is in fact constant. The offline solution
mostly takes $a_1$ values. The non-adaptive policy $\widehat \pi$ that has $p_{1,t}\equiv 1$ for all $t \in [n]$
and $p_{j,t}\equiv 0$ for all $j\neq 1$ and all $t \in [n]$ will achieve a constant regret.
Interestingly, in this same range, the index policy may not be as aggressive.
For instance, if $k = (f_1 - \epsilon)n$ we see from \eqref{eq:index} that
the index policy sets $\p_{1,t} = 1 - \epsilon/f_1$ spreading out the selection
of $a_1$ values throughout the time horizon and having a regret of order $\sqrt{n}$.

\begin{lemma}[Non-adaptive regret: small budget]\label{lm:smallbudget}
	For $ \epsilon = \tfrac{1}{2}\min\{f_m, f_{m-1}, \ldots, f_1\}$  suppose that
	$ k \leq n (f_1 - \epsilon)$. Then one has that
    \begin{equation*}
	\Voff^{*}(n,k)-\Vna^{*}(n,k) \leq \frac{a_2}{4 \epsilon}.
	\end{equation*}
\end{lemma}

Theorem \ref{thm:nonadaptive} and Lemma \ref{lm:smallbudget} are non-exhaustive.
There are ranges of $k/n$ that are covered by neither.
The purpose of these results is to show that while non-adaptivity is not universally bad (see Lemma \ref{lm:smallbudget}),
there is a non-trivial range of $(k,n)$ pairs for which the regret grows like $\sqrt{n}$
while the Budget-Ratio and the dynamic programming policy achieve $O(1)$ regret. 		
		
In the proof of Theorem \ref{thm:nonadaptive} we use three auxiliary lemmas,
the first of which provides a lower bound for the overshoot of a centered Bernoulli random walk.
		
\begin{lemma}\label{lm:martingale} 	
	Let $B_1, B_2, \ldots B_n$ be independent Bernoulli random variables
	with success probabilities $q_1, q_2, \ldots, q_n$, respectively.
    Let $N_n=\sum_{t\in[n]} \{B_t-q_t\}$ and $ \varsigma_n^2 = \Var[N_n] = \sum_{t \in [n]} q_t (1-q_t)$.
    Then, for any $\Upsilon > 0$ there exist constants
	$\beta_1 \equiv \beta_1( \Upsilon )$ and $\beta_2\equiv \beta_2(\Upsilon)$ such that
	\begin{equation}\label{eq:lower-bounds-lm-martingale}
		\E[(N_n-\Upsilon \varsigma_n )_+]
			\geq \beta_{1}\varsigma_n - (2 + 3 \sqrt{2}),
			\quad \quad
			\E[(-N_n-\Upsilon \varsigma_n )_+]\geq \beta_{1}\varsigma_n - (2 + 3 \sqrt{2}),
	\end{equation}
    and
    \begin{equation}\label{eq:upper-bound-lm-margingale}
			\E[(N_n+\Upsilon \varsigma_n )_+^2]\leq \beta_{2}  \varsigma_n^2.
	\end{equation}
\end{lemma}
		
The second auxiliary lemma shows that any good non-adaptive policy must be a perturbation of the index policy.
Specifically, if $s_j(\pi)=\sum_{t\in [n]}p_{j,t}f_j$ is the expected number of $a_j$ candidates
that policy $\pi$ selects under infinite budget,
then $s_j(\pi)$ is just a perturbation of $s_j^*$,
the solution of the deterministic relaxation \eqref{eq:DR}.
		
\begin{lemma}\label{lm:onlyindex}
	Fix $M < \infty$, $\epsilon = \tfrac{1}{2}\min\{f_m, \ldots, f_1\}$, and take any
	$(2 \epsilon \min_{j \in [m-1]}\abs{ a_j - a_{j+1}})^{-2} M^2 \leq n$.
	If $\pi = \{ p_{j,t} :  j \in [m] \text{ and } t \in [n]\}$ is a non-adaptive policy such that
	\begin{equation*} 
		DR(n,k)-\Von^{\pi}(n,k) \leq M \sqrt{n},
	\end{equation*}
	then
	$$
		s_j(\pi) =s^*_j \pm  \{\min_{j \in [m-1]}\abs{ a_j - a_{j+1}}\}^{-1} M \sqrt{n}
    	\quad \quad \text{for all } j\in [m].
	$$
\end{lemma}
		
For $ \epsilon = \tfrac{1}{2}\min\{f_m, f_{m-1}, \ldots, f_1\}$ and $(f_1 + \epsilon) n \leq k \leq  (1- f_m - \epsilon) n$, Lemma \ref{lm:multinomial-lemma} tells us that
$$
n f_1 	- \frac{1}{4 \epsilon }  \leq \E[Z_1^n] - \E[(Z_1^n - k)_+] = \E[\mathfrak{S}_1^n]
\quad \quad \text{and} \quad \quad
\E[\mathfrak{S}^n_m] \leq \frac{1}{4 \epsilon }.
$$
In words, when $k$ is bounded away from both 0 and $n$,
the offline algorithm selects---in expectation---all but a constant number of the highest values
and at most a constant number of the lowest values.
The optimal non-adaptive policy must do so as well.
In particular, it should have in ``most'' time periods $p_{1,t} = 1$ and $p_{m,t}=0$
so that the marginal probability of selection, $q_t$, is safely bounded away from zero and from one.
The last auxiliary lemma makes this intuitive idea formal.
		
\begin{lemma} \label{lm:varsigmalower}
	Let $ \epsilon = \tfrac{1}{2}\min\{f_m, f_{m-1}, \ldots, f_1\}$
	and suppose that  $(f_1 + \epsilon) n \leq k \leq  (1- f_m - \epsilon) n$.			
    Then, there is a constant $M \equiv M(\epsilon, a_1, a_2, a_m) < \infty$,
	such that an optimal non-adaptive policy must satisfy
	\begin{equation} \label{eq:qt-count-bounds}
			\sum_{t \in [n]}\1\left\{q_t\geq \frac{f_1}{2}\right\}\geq n-M\sqrt{n}
			\quad \mbox{ and } \quad
			\sum_{t \in [n]}\1\left\{1-q_t\geq \frac{f_m}{2}\right\}\geq n-M\sqrt{n}.
	\end{equation}
	Consequently,
	$$
	\sum_{t \in [n]}\1\left\{q_t\geq \frac{f_1}{2},1-q_t\geq \frac{f_m}{2}\right\}\geq n-2M\sqrt{n},
	$$
	and one has the lower bound
	$$
	\varsigma^2(\pi) = \sum_{t \in [n]}q_t(1-q_t)\geq \frac{f_1f_m}{4}(n-2M\sqrt{n}).
	$$
\end{lemma}
		
\subsection{Proof of Theorem \ref{thm:nonadaptive}}
		
For any non-adaptive policy $\pi = \{ p_{j,t} :  j \in [m] \text{ and } t \in [n]\}$
with associated stopping time $\nu = \min \{ r \geq 1 : \sum_{t \in [r]} B_t \text{ or } r \geq n\}$,
we have that policy $\pi$ does not make any selection after time $\nu$,
and we also have that $Z_j^{\nu} \leq Z_j^n$ for all $j \in [m]$.
Thus, if we recall the linear program \eqref{eq:LP}, use the
monotonicity of $\varphi(z_1, \ldots, z_m, \cdot)$ in $(z_1, \ldots, z_m)$, and recall the equivalence \eqref{eq:Voff-star},
we obtain that
$$
\Von^{\pi}(n,k) \leq \E[\varphi(Z_1^{\nu}, \ldots, Z_m^{\nu}, k)]
\leq \E[\varphi(Z_1^{n}, \ldots, Z_m^{n}, k)]=\Voff^{*}(n,k).
$$
Furthermore, since $\nu$ is the time at which $\pi$ exhausts the budget,
we also have that
\begin{align*}
	\Voff^{*}(n,k)-\Von^{\pi}(n,k)
	& \geq \E[\varphi(Z_1^{n}, \ldots, Z_m^{n}, k)]-\E[\varphi(Z_1^{\nu}, \ldots, Z_m^{\nu}, k)]\\
	& \geq  a_1 \left(\E[ \min\{ k, Z_1^n\}]-\E[\min\{ k, Z_1^{\nu}\} ]\right)\\
	& \geq  a_1(\E[(Z_1^n-Z_1^{\nu})]-n\1(Z_1^n> k)]).
\end{align*}
For $\epsilon = \tfrac{1}{2}\min\{f_m, f_{m-1}, \ldots, f_1\}$ and
$(f_1 + \epsilon) n  \leq k \leq (1-f_m - \epsilon) n$,
Hoeffding's inequality \citep[see, e.g.][Theorem 2.8]{BoucheronLugosiMassart:OXFORD2013} immediately tells us that
there is a constant $M \equiv M(\epsilon) < \infty$ such that  $n \P( Z_1^n>k ) \leq M$,
and Wald's lemma gives us that
\begin{equation}\label{eq:inter1}
	\E[Z_1^n-Z_1^{\nu}]= f_1\E[n-\nu],
\end{equation}
so the proof of the theorem is complete if we can prove an appropriate lower bound for $\E[n-\nu]$.
Lemma \ref{lem:index} tells us that for $f_1+ \epsilon \leq k/n$
the index policy has regret that is bounded above $(f_1 + \epsilon)^{-1} a_1 \sqrt{n}$, so that
it suffices to consider non-adaptive policies $\pi=\{p_{j,t},j\in[m] \mbox{ and } t\in [n]\}$ for which
$$
DR(n,k) - \Von^{\pi}(n, k) \leq (f_1 + \epsilon)^{-1} a_1 \sqrt{n}.
$$
If $(2 \epsilon)^{-2} \Mbar^2 \leq n$  for $\Mbar = \{ \min_{j \in [m-1]} \abs{a_j -a_{j+1}} (f_1 + \epsilon) \} ^{-1} a_1$,
Lemma \ref{lm:onlyindex} implies that
$$
	\sum_{t \in [n]} p_{j,t}f_j = s^*_j \pm \Mbar \sqrt{n}
    \quad \quad \text{for all } j \in [m].
$$
Since $\sum_{j \in [m]}s_j^*=k$, we also have that
\begin{equation}\label{eq:sum-qt}
	\sum_{t\in [n]}q_t
	=\sum_{t \in [n]} \sum_{j \in [m]} p_{j,t}f_j
	=\sum_{j \in [m]} s^*_j \pm m \Mbar \sqrt{n}
	=k \pm m \Mbar \sqrt{n}.
\end{equation}
Furthermore, since $0 \leq B_t \leq 1$ for all $t \in[n]$ we know that
\begin{equation}\label{eq:(sumBt-k)+equality}
	\Big(\sum_{t\in [n]} B_t - k\Big)_+
	=  \sum_{t=\nu+1}^{n} B_t  \leq n - \nu,
	\quad \quad \text{ so } \quad \quad
	\E \Big[ \Big(\sum_{t\in[n]} B_t - k\Big)_+ \Big] \leq \E[ n - \nu ].
\end{equation}
With $\varsigma^2(\pi) = \sum_{t \in [n]} q_t (1-q_t)$,
the estimate \eqref{eq:sum-qt} implies that $ k \leq \sum_{t \in [n]} q_t +  [m \Mbar \varsigma(\pi)^{-1}\sqrt{n} ] \varsigma(\pi)$,
and Lemma \ref{lm:varsigmalower} tells us that there is a constant $M \equiv M(\epsilon, m, a_1, \ldots, a_m)$
such that $m \Mbar \varsigma(\pi)^{-1}\sqrt{n} \leq M$.
In turn, we also have that  $k \leq \sum_{t \in [n]} q_t +  M \varsigma(\pi)$,
so after we subtract $\sum_{t\in[n]}B_t$ on both-sides, change sign, take the positive part,
and recall \eqref{eq:(sumBt-k)+equality}, we obtain the lower bound
$$
	\big( \sum_{t\in [n]}(B_t-q_t)- M \varsigma(\pi) \big)_+ \leq  \big( \sum_{t\in [n]}B_t-k \big)_+ \leq n - \nu.
$$
The quantity on the left-hand side is a sum of centered independent Bernoulli random variables so, when we take expectations,
Lemma \ref{lm:martingale} tells us that there is a constant $\beta_1 \equiv \beta_1 (\epsilon, m, a_1, \ldots, a_m)$ such that
$$
	\beta_1\varsigma(\pi) - (2 + 3 \sqrt{2})
	\leq \E \bigg[ \big( \sum_{t\in[n]} (B_t - q_t) - M\varsigma(\pi) \big)_+ \bigg]
	\leq \E[ n - \nu ].
$$
Plugging this last estimate back into \eqref{eq:inter1}
gives us that
$
f_1(\beta_{1}\varsigma(\pi)- 2 - 3\sqrt{2}) \leq \E[Z_1^n-Z_1^{\nu}],
$
and the theorem then follows after one uses one more time
the lower bound for $\varsigma(\pi)$ given in Lemma \ref{lm:varsigmalower}
and chooses $M \equiv M(\epsilon, m, a_1, \ldots, a_m)$ accordingly.

\section{Concluding remarks}

We have proved that in the multi-secretary problem with independent candidate abilities
drawn from a common finite-support distribution, the regret is constant and achievable by a multi-threshold policy.
In our model, the decision maker knows and makes crucial use of the distribution of candidate abilities.
Two obvious extension to consider are the problem instances in which the ability distribution is continuous
and/or unknown to the decision maker.

While one would like to think of the continuous distribution as a ``limit'' of discrete ones,
our analysis does build to a great extent on this discreteness, and our bounds depend on the cardinality of the support.
At this point, it is not clear if bounded regret is achievable also with continuous distributions.

For the case of unknown distribution, we conjecture that, with a finite support, the regret should be logarithmic in $n$.
Indeed, consider a ``stupid'' algorithm that uses the first $O(\log(n))$ steps to learn about the distribution (without concern for the objective) and,
at the end of the learning period, computes the threshold and runs with the Budget-Ratio policy thereafter.
Simple Chernoff bounds suggests that the likelihood of mis-estimation should be exponentially small.
Coupling this with the fact that the performance of the BR policy (specifically the fact that $\E[\tau]\geq n-M$)
is insensitive to small perturbations to the thresholds leads to our conjecture.

\ACKNOWLEDGMENT{This material is based upon work supported by the National Science Foundation under CAREER Award No. 1553274.}

\bibliographystyle{abbrvnat}



\newpage

\begin{APPENDICES}

\section{Regret lower bound in \citeauthor{kleinberg2017personalCommunications}'s example}\label{app:KleinbergExample}

Below we prove the regret lower bound in Lemma \ref{lem:kleinberg} in which the ability distribution
has support $\A = (a_3, a_2, a_1)$ and probability mass function $(\frac{1}{2} +  2 \epsilon, 2 \epsilon, \frac{1}{2}-4\epsilon)$.
The argument shows that, because of the special structure of this example,
there are two events (of non-negligible probability)---one being a perturbation of the other---where
the offline sort makes very different choices (taking all or none of the $a_2$ values)
but the optimal dynamic programming policy can do well only in one or the other but not in both.

\proof{Proof of Lemma \ref{lem:kleinberg}.}
In what follows we will assume that $\epsilon$ is such that $\frac{1}{\epsilon^2}$ and $\frac{1}{2\epsilon^2}$
are integers so that we can set $n_{\epsilon}=\frac{1}{\epsilon^2}$
and $k_{\epsilon}=\frac{n_\epsilon}{2}=\frac{1}{2\epsilon^2}$.
This makes the notation cleaner while not changing any of the arguments.
Next, we introduce the following events
\begin{eqnarray*}
A_\epsilon
 &=& \left\{ \frac{1}{\epsilon} \leq Z^{n_\epsilon}_2 \leq \min\big\{2 Z^{n_\epsilon/2}_2, \frac{2}{\epsilon}  \big\} \right\}
 = \left\{ \frac{1}{2}\E[Z^{n_\epsilon}_2] \leq Z^{n_\epsilon}_2 \leq \min\{2 Z^{n_\epsilon/2}_2, \E[Z^{n_\epsilon}_2]  \} \right\},
 \\
H_\epsilon
 &=& \left\{ Z^{n_\epsilon}_1 \geq k_\epsilon + \frac{2}{\epsilon} \right\}
 = \left\{ Z^{n_\epsilon}_1 \geq \E[Z^{n_\epsilon}_1] + \frac{6}{\epsilon} \right\},  \qquad \text{and}
 \\
L_\epsilon
 &=& \left\{ Z^{n_\epsilon}_1 \leq k_\epsilon - \frac{4}{\epsilon} \right\}
 = \bigg\{ Z^{n_\epsilon}_1 \leq \E[Z^{n_\epsilon}_1] \bigg\}.
\end{eqnarray*}
On the one hand, on the event $H_\epsilon$---in particular on $H_\epsilon \cap A_\epsilon$---the
offline policy takes only the highest values.
That is it only takes $a_1$ values and exhausts the budget $k_\epsilon$ while doing so.
On the other hand, on the event $L_{\epsilon}$---in particular on $L_\epsilon \cap A_\epsilon$---the
offline policy also takes some lower values. In particular, it takes
all of the $a_1$ values, all of the  $a_2$ values, and some $a_3$ values to exhaust the budget.
The events $H_\epsilon \cap A_\epsilon$ and $L_\epsilon \cap A_\epsilon$ are non-trivial as they happen with probability
that is bounded away from zero: there are constants $\epsilon_0>0$ and $\alpha>0$
such that $\P(H_\epsilon \cap A_\epsilon) \geq \alpha$ and $\P(L_\epsilon \cap A_\epsilon) \geq \alpha$
for all $\epsilon \leq \epsilon_0$.
To see this, note that the standard deviation of $Z^{n_\epsilon}_1$ is on the order of $1/\epsilon$,
so the central limit theorem implies that $\P(H_\epsilon), \P(L_\epsilon) \geq \alpha$ for all $\epsilon \leq \epsilon_0$.
Another scaling argument shows that $\P(A_\epsilon)\geq \alpha >0$ for all $\epsilon \leq \epsilon_0$
and for some re-defined $\alpha >0$.
Since $Z^{n_\epsilon}_1$  and $Z^{n_\epsilon}_2$ are dependent,
these observations and a standard condition argument that we omit
imply that there is another $\alpha >0$ such that
$\P(H_\epsilon \cap A_\epsilon) \geq \alpha$ for all $\epsilon \leq \epsilon_0$.
A similar analysis also gives that $\P(L_\epsilon \cap A_\epsilon)$
has the same property.

Next, we consider the random variable $\Lambda(n_\epsilon, k_\epsilon)$
that tracks the difference between the performance of the offline policy
and that of the dynamic programming policy:
$$
\Lambda(n_\epsilon, k_\epsilon)
    = a_1 (\frakS_1^{n_\epsilon} - S_1^{*,n_\epsilon})
    + a_2 (\frakS_2^{n_\epsilon} - S_2^{*,n_\epsilon})
    + a_3 (\frakS_3^{n_\epsilon} - S_3^{*,n_\epsilon}),
$$
and we note that the expected value
$\E[\Lambda(n_\epsilon, k_\epsilon)] = \Voff^*(n_\epsilon, k_\epsilon) - \Von^*(n_\epsilon, k_\epsilon)$
is the regret of the dynamic programming policy.
We now recall that $S_2^{*, n_\epsilon / 2}$ is the number of $a_2$-selections by the
optimal online policy up to time $n_\epsilon/2$, and we consider the event
$$
C_\epsilon = \{S_2^{*, n_\epsilon / 2} \geq  Z_2^{n_\epsilon}/4\}.
$$
On $C_{\epsilon}^c \cap L_\epsilon \cap A_\epsilon$, we then have that
\begin{equation*}
\E[\Lambda(n_\epsilon, k_\epsilon) \1(C_{\epsilon}^c \cap L_\epsilon \cap A_\epsilon)]
\geq  \frac{a_2 - a_3}{4} \E[ Z_2^{n_\epsilon} \1(C_{\epsilon}^c \cap L_\epsilon \cap A_\epsilon) ]
\geq  \frac{a_2 - a_3}{4 \epsilon} \P( C_{\epsilon}^c \cap L_\epsilon \cap A_\epsilon ).
\end{equation*}

For each sequence $\{ X_1,\ldots,X_{n_\epsilon}\} \in C_\epsilon \cap L_\epsilon \cap A_\epsilon$
we can now build a sequence in $C_\epsilon \cap H_\epsilon \cap A_\epsilon$
by keeping all the first $n_\epsilon/2$ values the same and replacing at most
$6/\epsilon$ $a_3$-values with $a_1$-values
while keeping the number of $a_2$-values constant,
and we call this resulting set $D_\epsilon$.
It is easy to see that since $f_3$ and $f_1$ are bounded away from zero,
then  $\P( D_\epsilon)\geq \widehat{\alpha} \P(C_{\epsilon}\cap L_{\epsilon}\cap A_{\epsilon})) $
for some $\widehat{\alpha}>0$ that does not depend on $\epsilon$.
Because the two sequences share the same first $n_\epsilon/2$ elements
we have that $S_2^{*,n_\epsilon/2} \geq Z_2^{n_\epsilon}/4$ also on $D_\epsilon$.
Then
\begin{align*}
    \E[\Lambda(n_\epsilon, k_\epsilon) \1(D_\epsilon )]&
    \geq (a_1 - a_2) \E[ S_2^{*,n_\epsilon/2} \1(D_\epsilon ) ]
    \geq \frac{ a_1 - a_2 }{4}\E[ Z_2^{n_\epsilon} \1( D_\epsilon ) ]
    \geq  \frac{a_1 - a_2}{4 \epsilon}\widehat{\alpha} \P(C_{\epsilon}\cap L_{\epsilon}\cap A_{\epsilon}).
\end{align*}
Since $\P(C_{\epsilon}^c\cap L_{\epsilon}\cap A_{\epsilon})=\P( L_{\epsilon}\cap A_{\epsilon})-\P(C_{\epsilon}\cap L_{\epsilon}\cap A_{\epsilon})$, we have that
$ \max\{\P(C_{\epsilon}^c\cap L_{\epsilon}\cap A_{\epsilon}),
\widehat{\alpha} \P(C_{\epsilon}\cap L_{\epsilon}\cap A_{\epsilon})\}
\geq  \P( L_{\epsilon}\cap A_{\epsilon})\min\{\widehat{\alpha},1\}$ so that we finally obtain
\begin{align*}
\Voff^*(n_\epsilon, k_\epsilon) - \Von^*(n_\epsilon, k_\epsilon)
     &   =  \E[\Lambda(n_\epsilon, k_\epsilon)]\\
     & \geq \max\{ \E[\Lambda(n_\epsilon, k_\epsilon) \1( C^c \cap L_\epsilon \cap A_\epsilon )],
                  \E[\Lambda(n_\epsilon, k_\epsilon) \1( D_\epsilon ) ] \}\\
     & \geq \min\left\{ \frac{a_2-a_3}{4\epsilon}, \frac{a_1-a_2}{4\epsilon} \right\}
            \max\{\P(C_{\epsilon}^c\cap L_{\epsilon}\cap A_{\epsilon}), \widehat{\alpha} \P(C_{\epsilon}\cap L_{\epsilon}\cap A_{\epsilon})\} \\
     & \geq \min\left\{ \frac{a_2-a_3}{4\epsilon}, \frac{a_1-a_2}{4\epsilon} \right\}
            \min\{\widehat{\alpha},1\}\P(L_{\epsilon}\cap A_{\epsilon}),
\end{align*}
so the result follows after one recalls that $\P(L_{\epsilon},A_{\epsilon})\geq \alpha$
and takes $\Gamma = \frac{1}{4} \min\{a_2 - a_3, a_1 - a_2\} \min \{\widehat{\alpha},1\} \alpha $.
\halmos\endproof\vspace{0.25cm}

\section{Proofs of auxiliary lemmas}\label{sec:proof-of-lemmas}

\proof{Proof of Lemma \ref{lm:multinomial-lemma}.}
Since $\E[B] = p n$ and $(p+\varepsilon)n \leq k$, for any $u > 0$ we have the tail bound
$$
\P ( (B-k)_+\geq u)
=\P(B - p n \geq k- p n + u)
\leq \P(B - p n \geq \varepsilon n + u) ,
$$
and Hoeffding's inequality \citep[see, e.g.][Theorem 2.8]{BoucheronLugosiMassart:OXFORD2013}
implies that
\begin{equation*}
\P ( (B-k)_+ \geq u) \leq \exp\{ -2 \varepsilon^2 n  - 4 \varepsilon u - 2 u^2 / n\},
\quad \quad \text{for all } (p+\varepsilon)n \leq k.
\end{equation*}
By integrating both sides for $u \in [0, \infty)$ we then obtain
for all $(p+\varepsilon)n \leq k$ that
\begin{align*}
\E[ (B-k)_+]
\leq \int_0^{\infty}\P ( (B-k)_+ \geq u ) \,du
\leq \exp\{ -2 \varepsilon^2 n  \} \int_0^{\infty} \exp\{ - 4 \varepsilon u  \} \, du
\leq \frac{1}{4 \varepsilon}.
\end{align*}

To prove the second bound in  \eqref{eq:binomial-lemma-E(X-k)bound},
one applies the first bound to the binomial random variable $B' = n- B$
and the budget $k' = n - k$.
\halmos\endproof\vspace{0.25cm}

\proof{Proof of Lemma \ref{lem:index}.}
The index policy $\id  = \{ \p_{j,t} :  j \in [m] \text{ and } t \in [n]\}$
defined by \eqref{eq:index} is such that all values greater than or equal to $a_{j_{\id}-1}$ are selected
together with a fraction of the $a_{j_{\id}}$ values.
Formally, by Wald's lemma, we have that
$$
\E[S_j^{\id, \nu}] = \begin{cases}
                        f_j\E[\nu] & \mbox{if } j \leq j_{\id} - 1 \\
                        (k/n - \Fbar(a_{j_{\id}})) \E[\nu] & \mbox{if } j = j_{\id} \\
                        0 & \mbox{if } j \geq j_{\id} + 1,
                      \end{cases}
$$
where $S_j^{\id,\nu}=\sum_{t=1}^{\nu}   \sigma_t^{\id} \1(X_t=a_j)$ is the number of $a_j$-candidates
selected by the index policy by time $\nu$. Recall that now that for any policy $\pi$ we have the inequalities
$$
\Von^{\pi}(n,k) \leq \Voff^{*}(n,k) \leq  DR(n,k)=\sum_{j=1}^{j_{\id} - 1}a_j f_j n +a_{j_{\id}}(k-n\Fbar(a_{j_{\id}})),
$$
so that
$$
\Voff^*(n,k) - \Von^{\id}(n,k) \leq DR(n,k)- \Von^{\id}(n,k)\leq a_1\E[n-\nu],
$$
and to bound the regret as desired it suffices to obtain an upper bound for $\E[n-\nu]$.

Let $\{B_t: t\in[n]\}$ be i.i.d Bernoulli random variables with success probability
$q = \sum_{i=1}^{j_{\id}-1}f_j+ k/n-\Fbar(a_{j_{\id}})= k/n$.
If
$$
N_r = \sum_{t\in[r]}B_t - q r
$$
is the centered number of candidates the index policy selects by time $r$,
we then have for $t\geq 0$ and $q = k/n$ that,
$$
n-\nu\geq t
\quad \text{ if and only if } \quad
N_{n-\lceil t\rceil} \geq k - q(n - \lceil t\rceil) = q \lceil t\rceil.
$$
In turn, Kolmogorov's maximal inequality \citep[See, e.g.][Theorem 22.4]{Bil:WILEY1995}
tells us that for any $t > 0$
$$
\P\left\{n-\nu\geq t\sqrt{n}\right\}
      =  \P\left\{N_{n-\lceil t\sqrt{n}\rceil } \geq q\lceil t\sqrt{n}\rceil \right\}
    \leq \P( \sup_{t\in [n]} \abs{ N_t } \geq q t\sqrt{n} )
    \leq \frac{\E[N_n^2]}{q^2  t^2 n} = \frac{1-q}{q t^2}.
$$
It then follows that
$$
\E\left[\frac{n-\nu}{\sqrt{n}}\right]\leq 1+\int_1^{\infty} \P( n-\nu\geq t\sqrt{n} ) \, dy \leq \frac{1}{q},
$$
so that $\E[n-\nu] \leq q^{-1}\sqrt{n} \leq \varepsilon^{-1}\sqrt{n}$
for any $\varepsilon \in (0,1)$ and all pairs $(n,k)$ such that $\varepsilon \leq k/n$.
\halmos\endproof\vspace{0.25cm}

\proof{Proof of Lemma \ref{lm:smallbudget}.}
The non-adaptive policy $\widehat \pi$ that takes all values $a_1$ and rejects all others achieves bounded regret.
To see this, notice that
$$
\sum_{j=2}^{m}\frakS_j^n
    = \sum_{j=2}^{m} \min\big\{ Z_j^n , ( k-\sum_{i\in[j-1]} Z_i^n )_+ \big\}
    = ( k - Z_1^n )_+.
$$
Since the random variable $Z_1^n$ is Binomial with parameters $n, f_1$,  and $0 \leq k \leq n(f_1 - \epsilon)$,
then Lemma \ref{lm:multinomial-lemma} implies that
$$
\E\left[ \sum_{j=2}^{m} \frakS_j^n \right]  = \E[( k - Z_1^n )_+] \leq \frac{1}{4 \epsilon}.
$$
In turn, the value of the offline solution for $ k \leq n(f_1 - \epsilon)$ satisfies the bound
$$
\Voff^{*}(n,k)
    \leq a_1\E[\frakS_1^n] + \frac{a_2}{4 \epsilon}
    = a_1\E[\min\{Z_1^n, k\}] + \frac{a_2}{4 \epsilon}.
$$
The non-adaptive policy $\widehat \pi$.
takes all $a_1$ values and none of the others, until it runs out of budget at time
$\nu = \min\{r\geq 1: \sum_{t\in [r]}  B_t \geq k \text{ or } r \geq n\}$,
so that $\E[S_j^{\widehat \pi,n}]=0$ for all $j\geq 2$ and
$$
\E[S_1^{\widehat \pi,n}] = \E[\sum_{t\in [\nu]} B_t \1(X_t = a_1)].
$$
Furthermore, we have that $S_1^{\widehat \pi,n}=Z_1^n$ if $Z_1^n<k$, and it equals $k$ otherwise.
Thus,
$$
\E[S_1^{\widehat \pi,n}]=\E[Z_1^{\nu}]=\E[\min\{Z_1^n, k\}] = \E[\frakS_1^n],
$$
so we finally have the bound
$$
\Voff^{*}(n,k) - \frac{a_2}{4 \epsilon} \leq a_1 \E[\min\{Z_1^n, k\}] = \Von^{\pi}(n,k),
$$
just as needed.
\halmos\endproof\vspace{0.25cm}

\proof{Proof of Lemma \ref{lm:martingale}.}
Let $\mathcal{Z}$ denote a  normal random variable with mean zero and variance 1.
Then,
\begin{equation}\label{eq:absolute-difference-bound}
\abs{ \E[ (N_n-\Upsilon \varsigma_n )_+ ]
	- \varsigma_n \E[ ( \mathcal{Z} -\Upsilon )_+ ] }
\leq d_W(N_n,  \varsigma_n \mathcal{Z} ) = \sup_{h \in \rm{Lip_1}} \abs{ \E[ h(N_n)] - \E[h( \varsigma_n \mathcal{Z} )]},
\end{equation}
where the $\rm{Lip_1}$ is the set of Lipschitz-1 functions on $\R$. The random variables $B_1, B_2, \ldots, B_n$ are
independent Bernoulli with success probabilities $q_1, q_2, \ldots, q_n$
so that
$$
\E[ \abs{ B_t - q_t }^3] \leq 2q_t(1-q_t)
\quad \quad \text{and} \quad \quad
\E[ ( B_t - q_t )^4] \leq 2q_t(1-q_t)
$$
and
$$
\sum_{t\in [n]}\E[ \abs{ B_t - q_t }^3] \leq 2  \varsigma_n^2
\quad \quad \text{and} \quad \quad
\sum_{t\in[n]}\E[ ( B_t - q_t )^4] \leq  2  \varsigma_n^2.
$$
A version of Stein's lemma \citep[see, e.g.,][Theorem 3.6]{ross2011fundamentals} for the sum of independent (not necessarily identically distributed)
random variables tells us that
$$
d_W(N_n/ \varsigma_n, \mathcal{Z})
\leq \frac{2}{ \varsigma_n} + \frac{3 \sqrt{2}}{ \varsigma_n},
$$
so by the homogeneity of the distance function $d_W$, we also have that
$$
d_W(N_n,  \varsigma_n \mathcal{Z} ) =  \varsigma_n d_W(N_n/ \varsigma_n, \mathcal{Z}) \leq 2 + 3 \sqrt{2}.
$$
The inequality \eqref{eq:absolute-difference-bound} then implies that
$$
\varsigma_n \E[ ( \mathcal{Z} -\Upsilon )_+ ] - (2+3\sqrt{2})\leq
\E[ (N_n-\Upsilon \varsigma_n )_+ ],
$$
and one can obtain an immediate lower bound for the left-hand side is by
$$
\varsigma_n \Upsilon \P( \mathcal{Z} \geq 2 \Upsilon )\leq \varsigma_n \, \E[ ( \mathcal{Z}- \Upsilon) \1(  \mathcal{Z}\geq 2 \Upsilon) ]\leq
\varsigma_n \E[ ( \mathcal{Z} -\Upsilon )_+ ].
$$
The left inequality in  \eqref{eq:lower-bounds-lm-martingale} then follows setting $\beta_1 = \Upsilon \P( \mathcal{Z} \geq  2 \Upsilon )$.
For the right inequality we combine the earlier argument with the symmetry of the normal distribution and the Lipschitz-1 continuity of the map
$x \mapsto (- x - \Upsilon \varsigma_n)_+$.

The argument for the second inequality in \eqref{eq:upper-bound-lm-margingale} is standard.
Since $\E[N_n] = 0$ and $\E[N_n^2] =  \varsigma_n^2$, we have that
$$
\E[( N_n+\Upsilon  \varsigma_n )_+^2]
\leq \E[(N_n + \Upsilon  \varsigma_n)^2]
=  \varsigma_n^2 +  \Upsilon^2  \varsigma_n^2,
$$
so taking $\beta_2 = 1+\Upsilon^2$ concludes the proof.
\halmos\endproof\vspace{0.25cm}

\proof{Proof of Lemma \ref{lm:onlyindex}.}
Given a non-adaptive policy $\pi = \{ p_{j,t}: j \in [m] \text{ and } t \in [n] \}$
and a constant $M < \infty$ we let
\begin{equation}\label{eq:Mbar}
\Mbar = \frac{M}{\min_{j \in [m-1]} \abs{ a_j - a_{j+1}}}
\quad \quad \text{and} \quad \quad
\iota = \argmax_{j \in [m]}\abs{s_j(\pi)-s^*_j },
\end{equation}
and we show that if
\begin{equation}\label{eq:max-assumption}
\abs{ s_{\iota}(\pi)-s^*_{\iota} }  = \max_{j \in [m]}\abs{ s_j(\pi)-s^*_j } > \Mbar \sqrt{n},
\end{equation}
then
\begin{equation*}
DR(n,k)-\Von^{\pi}(n,k)> M \sqrt{n}.
\end{equation*}
We let $j_{\id}\in[m]$ be the index such that $\Fbar(a_{j_{\id}}) \leq k/n < \Fbar(a_{j_{\id}+1})$.
There are two cases to consider: (i)
\eqref{eq:max-assumption} is attained at $\iota \geq j_{\id}+1$ ($k/n<\Fbar(a_{j_{id}})<\Fbar(a_{\iota})$) in which case $0=s_{\iota}^*<s_{\iota}(\pi)$, and (ii) $\iota \leq j_{\id}$ ($\Fbar(a_{\iota}\leq \Fbar(a_{j_{id}})\leq k/n$) in which case $0<s_{\iota}^*$ and $s_{\iota}(\pi)<s_{\iota}^*$.

We begin with the first case.
That is, we assume that \eqref{eq:max-assumption} is attained for some
$\iota \geq j_{\id}+1$ when $s_{\iota}^* = 0$,
and we obtain that
$s_{\iota}(\pi) \geq s^*_{\iota} +\Mbar \sqrt{n}= \Mbar \sqrt{n}$.
To estimate the gap between $DR(n,k)$ and $\Von^\pi(n,k)$,
consider a version of the deterministic relaxation \eqref{eq:DR}
that requires the selection of at least $\Mbar \sqrt{n}$ candidates with ability $a_{\iota}$
out of the $f_\iota n$ available.
Since $\Mbar \leq 2 \epsilon \sqrt{n}$
we have that $ \Mbar \sqrt{n} \leq f_j n $ for all $j \in [m]$,
so we write the optimization problem as
\begin{eqnarray*}
	DRC(n,k, \iota)  = & \underset{s_1, \ldots, s_m}{\max} &  \sum_{j\in[m]} a_j s_j\\
	& \text{s.t.} &  0 \leq s_j\leq \E[Z_j^n] \mbox{ for all } j\in[m] \\
    &  & \Mbar \sqrt{n} \leq s_{\iota}   \notag\\
	&  & \sum_{j\in [m]}s_j \leq k, \notag
\end{eqnarray*}
and we obtain that
\begin{equation*}
\Von^{\pi}(n,k) \leq DRC(n,k, \iota) \leq DR(n,k).
\end{equation*}
The unique maximizer $(\check{s}_1,\ldots,\check{s}_{m})$ of $DRC(n,k, \iota)$
is given by
$$
\check{s}_j = \begin{cases}
s^*_j         & \mbox{if } j \leq j_{\id} - 2 \\
s^*_{j_{\id}-1} - (\Mbar \sqrt{n} - s^*_{j_{\id}})_+     & \mbox{if } j = j_{\id}-1 \\
(s^*_{j_{\id}} - \Mbar \sqrt{n})_+, & \mbox{if } j = j_{\id} \\
\Mbar \sqrt{n} & \mbox{if } j = \iota\\
0 & \mbox{otherwise,}
\end{cases}
$$
so that the difference between the value of the deterministic relaxation and the value of its constrained version
is given by
$$
DR(n,k)-DRC(n,k,\iota)
= a_{j_{\id}-1} (\Mbar \sqrt{n} - s^*_{j_{\id}})_+
+ a_{j_{\id}} [s^*_{j_{\id}} - (s^*_{j_{\id}} - \Mbar \sqrt{n})_+]
- a_{\iota} \Mbar \sqrt{n}.
$$
Because $j_{\id} - 1 < j_{\id} \leq \iota - 1 < \iota$, the monotonicity $a_{\iota}< a_{\iota-1} \leq a_{j_{\id}} <  a_{j_{\id}-1}$ gives us the lower bound
$$
DR(n,k)-DRC(n,k,\iota)
\geq a_{j_{\id}-1} \Mbar \sqrt{n} - a_{\iota} \Mbar \sqrt{n}
> [ a_{\iota -1} -a_{\iota} ] \Mbar \sqrt{n},
$$
and since $V_{on}^{\pi}(n,k)\leq DRC(n,k,\iota)$ we have
\begin{equation}\label{eq:DR-Von-iota1}
DR(n,k)-\Von^{\pi}(n,k)  > [ a_{\iota -1} -a_{\iota} ] \Mbar \sqrt{n}
\quad \quad \text{for all } \iota \geq j_{\id}+1.
\end{equation}

A similar inequality can be obtained for second case in which when \eqref{eq:max-assumption} is attained at some index $\iota \leq j_{\id}$.
In this case we would consider a version of the deterministic relaxation \eqref{eq:DR} with the additional constraint
$s_{\iota} \leq n f_{\iota} - \Mbar \sqrt{n}$ or $s_{j_{\id}} \leq k - n \Fbar(a_{j_{\id}}) - \Mbar \sqrt{n}$.
This analysis then implies the bound
$$
DR(n,k)-\Von^{\pi}(n,k)  > [ a_{\iota} -a_{\iota + 1} ] \Mbar \sqrt{n}
\quad \quad \text{for all } \iota \leq j_{\id},
$$
so if we recall \eqref{eq:DR-Von-iota1} and use the definition of $\Mbar$ in \eqref{eq:Mbar}
we finally obtain that
$$
DR(n,k)-\Von^{\pi}(n,k) > M \sqrt{n},
$$
concluding the proof of the lemma.
\halmos\endproof\vspace{0.25cm}

\proof{Proof of Lemma \ref{lm:varsigmalower}.}
Fix any non-adaptive policy $\pi \in \Pi_{\na}$ and recall that $S_1^{\pi,\nu}$ is the number
of $a_1$-candidates that the policy selects.
If
\begin{eqnarray*}
	\varphi_{-1}(z_2, \ldots, z_m, k)  = & \underset{s_2, \ldots, s_m}{\max} &  \sum_{j=2}^m a_j s_j\\
	& \text{s.t.} &  0 \leq s_j\leq z_j \mbox{ for all } j=2,\ldots,m \notag\\
	&  & \sum_{j=2}^m s_j \leq k,  \notag
\end{eqnarray*}
then
\begin{equation}\label{eq:Von-varphi-bound}
\Von^{\pi}(n,k)\leq a_1\E[S_1^{\pi,\nu}]+\E[\varphi_{-1}(Z_2^n,\ldots,Z_m^n,k-S_1^{\pi,\nu})].
\end{equation}
Since $S_1^{\pi,\nu}\leq \mathfrak{S}_1^n=\min\{Z_1^n,k\}$, the monotonicity in $k$ of $\varphi_{-1}(\cdot, k)$
gives us the upper bound
$$
\varphi_{-1}(Z_2^n,\ldots,Z_m^n,k-S_1^{\pi,\nu})
    \leq \varphi_{-1}(Z_2^n,\ldots,Z_m^n,k-\mathfrak{S}_1^{n})+a_2 (\mathfrak{S}_1^n-S_1^{\pi,\nu}),
$$
so when we take expectations and recall \eqref{eq:Von-varphi-bound}, we obtain that
$$
\Von^{\pi}(n,k)\leq a_1\E[S_1^{\pi,\nu}]+\E[\varphi_{-1}(Z_2^n,\ldots,Z_m^n,k-\mathfrak{S}_1^{n})]
        + a_2 (\E[\mathfrak{S}_1^n] - \E[S_1^{\pi,\nu}]).
$$
We now note that  $\Voff^{*}(n,k)=a_1\E[\mathfrak{S}_1^n] + \E[\varphi_{-1}(Z_2^n,\ldots,Z_m^n,k-\mathfrak{S}_1^n)]$,
so when we use this last decomposition in the displayed equation above we conclude that 	
$$
\Von^{\pi}(n,k)\leq \Voff^{*}(n,k)-(a_1-a_2)\left(\E[\mathfrak{S}_1^{n}]-\E[S_1^{\pi,\nu}]\right).
$$

Let $M \equiv M(\epsilon, a_1, a_2)$ be such that the index policy achieves an $(a_1-a_2)M\sqrt{n}$
regret (see Lemma \ref{lem:index} with $\varepsilon = f_1 + \epsilon$).
Then, if $\pi$ is such that $\E[S_1^{\pi,\nu}] < \E[\mathfrak{S}_1^n]-M\sqrt{n}$, then
$$
\Von^{\pi}(n,k)\leq \Voff^{*}(n,k)-(a_1-a_2)M\sqrt{n},
$$
so that $\pi$ cannot be optimal.
In other words, the optimal policy must satisfy
\begin{equation} \label{eq:optimal-napolicy-S1-constraint}
\E[S_1^{\pi,\nu}]\geq \E[\mathfrak{S}_1^n]-M\sqrt{n}.
\end{equation}
Since $(f_1 + \epsilon) n \leq k$, we have from Lemma \ref{lm:multinomial-lemma} that
$$
\E[\mathfrak{S}_1^n] = \E[Z_1^n] - \E[(Z_1^n - k )_+] \geq  nf_1 - \frac{1}{4 \epsilon},
$$
which, together with \eqref{eq:optimal-napolicy-S1-constraint}, implies that the optimal non-adaptive policy
must have
$$
\sum_{t\in[n]}p_{1,t}f_1 \geq \E[S_1^{\pi,\nu}]\geq nf_1 - \frac{1}{ 4 \epsilon }-M\sqrt{n}.
$$
Since $ p_{1,t}f_1 \leq q_t$ for all $t \in [n]$, it follows that
there is another constant $M \equiv M(\epsilon, a_1, a_2)$ such that
$$
\sum_{t \in [n]}\1\left\{q_t\leq \frac{f_1}{2}\right\} \leq \sum_{t \in [n]}\1\left\{p_{1,t}\leq \frac{1}{2}\right\}\leq M\sqrt{n},
$$
concluding the proof of the left inequality of \eqref{eq:qt-count-bounds}.

The argument for the right inequality of \eqref{eq:qt-count-bounds} is similar.
We let
\begin{eqnarray*}
	\varphi_{-m}(z_1, \ldots, z_{m-1}, k)  = & \underset{s_1, \ldots, s_{m-1}}{\max} &  \sum_{j=1}^{m-1} a_j s_j\\
	& \text{s.t.} &  0 \leq s_j\leq z_j \mbox{ for all } j=1,\ldots,m-1 \notag\\
	&  & \sum_{j=1}^{m-1} s_j \leq k,  \notag
\end{eqnarray*}
and note that
$$
\Von^{\pi}(n,k)
    \leq a_m\E[S_m^{\pi,\nu}]+\E[\varphi_{-1}(Z_1^n,\ldots,Z_m^n,k-S_m^{\pi,\nu})]
    \leq a_m\E[S_m^{\pi,\nu}]+\E[\varphi_{-1}(Z_1^n,\ldots,Z_m^n,k)].
$$
The decomposition $\Voff^{*}(n,k)= \E[\varphi_{-m}(Z_1^n,\ldots,Z_{m-1}^n,k)] + a_m \E[\mathfrak{S}_m^n]$,
then implies that
$$
\Von^{\pi}(n,k) \leq \Voff^{*}(n,k) -  a_m (\E[\mathfrak{S}_m^n] - \E[S_m^{\pi,\nu}]).
$$
Here we have that $k \leq (1 - f_m - \epsilon)n$ so Lemma \ref{lm:multinomial-lemma} tells us that
$ \E[\mathfrak{S}_m^n] \leq (4 \epsilon)^{-1}$ and if $M \equiv M(\epsilon, a_m)$ is the
constant such that the index policy achieves $a_m M \sqrt{n}$ regret,
then we see that policy $\pi$ cannot be optimal if
$ \E[S_m^{\pi,\nu}] > M \sqrt{n}$.
Since $(1- p_{m,t})f_m \leq 1 - q_t$, this observation then gives us
another constant $M \equiv M(\epsilon, a_m)$ such that
$$
\sum_{t \in [n]}\1\left\{1 - q_t \leq \frac{f_m}{2}\right\}
\leq \sum_{t \in [n]}\1\left\{1 - p_{m,t} \leq \frac{1}{2}\right\}
  =  \sum_{t \in [n]}\1\left\{ p_{m,t} \geq \frac{1}{2}\right\}
\leq M\sqrt{n},
$$
which completes the proof.
\halmos\endproof\vspace{0.25cm}

\section{The dynamic programming formulation}\label{sec:bellman}

Let $v_\ell(w,\kappa)$ denote the expected total ability when $\ell$ candidates are yet to be inspected
(the are $\ell$ ``periods'' to-go), the current level of accrued ability is $w$,
and at most $\kappa$ can be selected ($\kappa$ is the residual recruiting budget).
The sequence of functions $\{ v_\ell : \R_+ \times \Z_+ \rightarrow \R_+ \text{ for } 1 \leq \ell < \infty\}$
then satisfy the Bellman recursion
\begin{equation}\label{eq:v-recursion}
v_\ell(w,\kappa) = \sum_{j\in [m]} \max\{v_{\ell-1}(w+a_j,\kappa-1), v_{\ell-1}(w,\kappa)\} f_j,
\end{equation}
with the boundary conditions
$$
v_0(w,\kappa) = w \text{ for all $\kappa \in\bbZ_+$, and $w\geq 0$}
\quad \quad \text{and} \quad \quad
v_\ell(w,0) = w   \text{ for all $\ell \in [n]$ and $w\geq 0$}.
$$
If $n$ is the number of available candidates and $k$ is the initial recruiting budget,
the optimality principle of dynamic programming then implies that
$$
\Von^*(n,k) = v_n(0,k) \quad \quad \text{for all } n,k \in\bbZ_+.
$$
The Bellman equation \eqref{eq:v-recursion} implies that
the optimal online policy takes the form of a threshold policy,
and the next proposition proves that the optimal thresholds depend only on the residual budget, $\kappa$,
and not on the accrued ability $w$.

\begin{proposition}[State-space Reduction]\label{prop:affine-decomposition}
Let  $\{g_\ell: \Z_+ \rightarrow \R_+ \text{ for } 1 \leq \ell < \infty\}$
satisfy the recursion
\begin{equation}\label{eq:g-recursion}
    g_\ell(\kappa) = \sum_{j\in [m]} \max\{ a_j + g_{\ell-1}(\kappa-1), g_{\ell-1}(\kappa) \} f_j,
\end{equation}
with the boundary conditions
\begin{equation}\label{eq:g-boundary-condition}
	g_0(\kappa)=0 \mbox{ for al } \kappa\in\bbZ_+ \quad \text{and} \quad
	g_\ell(0) = 0 \text{ for all } 1 \leq \ell < \infty.
\end{equation}
Then,
\begin{equation}\label{eq:affine-decomposition}
	v_\ell(w,\kappa) = w + g_\ell(\kappa) \quad \quad \text{for all $w\in \R_+$ and $\kappa\in\bbZ_+$}.
\end{equation}
Consequently, with $\ell$ periods to-go and a residual budget of $\kappa$,
it is optimal to select a candidate with ability $a_j$ if and only if
$$
    a_j \geq h_\ell(\kappa) \equiv g_{\ell-1}(\kappa) - g_{\ell-1}(\kappa-1).
$$
\end{proposition}

Under the optimal online policy, ability levels strictly smaller than $h_\ell(\kappa)$ are skipped
and the rest are selected.
Since the Markov decision problem (MDP) associated with the Bellman equation \eqref{eq:g-recursion}
is a finite-horizon problem with finite state space and a finite number of actions,
this deterministic policy is also optimal within the larger class of non-anticipating policies
\citep[See, e.g.][Corollary 8.5.1.]{BeS:78}.

\proof{Proof of Proposition \ref{prop:affine-decomposition}.}
This is an induction proof.
For $\ell=1$, we have by \eqref{eq:v-recursion} that
$$
v_1(w,0) = w   \quad \text{and} \quad
v_1(w,\kappa) = w + \sum_{j\in[m]} a_j f_j=w+\E[X_1] \text{ for all } \kappa \geq 1.
$$
Taking
$$
g_1(\kappa) = \begin{cases}
0                 & \mbox{if } \kappa=0 \\
\E[X_1]           & \mbox{if } \kappa \geq 1,
\end{cases}
$$
one has that
$$
v_1(w,\kappa) = w + g_1(\kappa) \quad \quad \text{for all $w\geq 0$ and all $\kappa\in\bbZ_+$}.
$$
Next, as induction hypothesis suppose that we have the decomposition
$$
v_{\ell-1}(w,\kappa) = w + g_{\ell-1}(\kappa) \quad \quad \text{for all $w \geq 0$ and all $\kappa \in\bbZ_+$}.
$$
If $\kappa=0$ the decomposition is satisfied with $g_\ell(0)=0$ and
\begin{equation}\label{eq:gk-zero}
v_\ell(w, 0) = w+g_\ell(0)=w.
\end{equation}
For $\kappa \geq 1$, the Bellman equation \eqref{eq:v-recursion} and the
induction assumption imply
\begin{align*}
v_\ell(w,\kappa)  & = w + \sum_{j\in [m]} \max\{a_j + v_{\ell-1}(w+a_j,\kappa-1) - w - a_j, v_{\ell-1}(w,\kappa)-w\} f_j\\
& = w + \sum_{j\in[m]} \max\{a_j + g_{\ell-1}(\kappa-1) , g_{\ell-1}(\kappa)\} f_j.
\end{align*}
Defining recursively
$$
g_\ell(\kappa) = \sum_{j\in [m]} \max\{a_j + g_{\ell-1}(\kappa-1) , g_{\ell-1}(\kappa)\} f_j \quad \quad \text{for all } \kappa \geq 1,
$$
we then have together with \eqref{eq:gk-zero} that \eqref{eq:affine-decomposition} holds for all $w \geq 0$   and all $\kappa \in\bbZ_+$. From this argument it also follows that $g_\ell(\kappa)$ satisfy the recursion \eqref{eq:g-recursion} with the boundary condition \eqref{eq:g-boundary-condition}.
\halmos\endproof\vspace{0.25cm}

\section{Deterministic relaxation revisited\label{sec:det_relax}}

The subset $\T'$ is identified by taking
$0 < \epsilon' < \epsilon=\frac{1}{2} \min\{f_m, f_{m-1}, \ldots, f_1\}$, and letting
$$
\T' = \{(n,k) \in \T: \Fbar(a_{j}) + \epsilon' \leq k / n \leq \Fbar(a_{j+1}) - \epsilon' \mbox{ for some } j\in[m]\}.
$$

\begin{proposition}[Deterministic Relaxation Gap] \label{prop:DR-bounds}
There exists a constant $M \equiv M( \epsilon, m, a_m, a_{m-1}, \ldots, a_1)$
such that
$$
	0 \leq DR(n,k) - \Voff^{*}(n,k)\leq M\sqrt{n}
	\quad \quad \text{for $(n, k)\in \T$.}
$$	
Furthermore,
$$
	0 \leq DR(n,k) - \Voff^{*}(n,k)\leq \frac{a_1 m}{ 4 \epsilon'}
	\quad \quad \text{for $(n, k)\in \T'$.}
$$	
\end{proposition}

The proof of Proposition \ref{prop:DR-bounds} requires the following auxiliary lemma.

\begin{lemma}\label{lm:minmax-deviations}
There exists a constant $M \equiv M( \epsilon, m, a_m, a_{m-1}, \ldots, a_1)$ such that, for all $j\in [m]$,
$$
    \E[\frakS_j^n]
        = \E\left[\min \Big\{ Z_j^n, \Big(k - \sum_{i \in [j-1]} Z_i^n\Big)_+\Big\}\right]
	    = \min \left\{ \E[Z_j^n], \Big(k - \sum_{i \in [j-1]} \E[Z_i^n]\Big)_+\right\}  \pm M \sqrt{n}.
$$
In turn,
$$
    \E[\frakS_j^n]=s_j^{*}\pm M\sqrt{n}, \mbox{ for all } j\in[m].
$$
\end{lemma}

\proof{Proof.}
Let $\xi^n_j = Z_j^n - \E[Z_j^n]$ for all $j \in [m]$
and obtain the representation
$$
\min \Big\{ Z_j^n, \Big(k - \sum_{i \in [j-1]} Z_i^n\Big)_+\Big\}
= \min \Big\{ \E[Z_j^n] + \xi^n_j, \Big(k - \sum_{i \in [j-1]} \E[Z_i^n] - \sum_{i \in [j-1]} \xi^n_i\Big)_+\Big\}.
$$
For any real numbers $a, b, c, d$, we have that
$$
\min\{ a, \max\{c, 0\}\} - \abs{b} - \abs{d}
\leq \min\{ a +b , \max\{c + d, 0\}\}
\leq \min\{ a, \max\{c, 0\}\} + \abs{b} + \abs{d},
$$
so that by setting,
$a = \E[Z_j^n]$, $b = \xi^n_j$, $c = k - \sum_{i \in [j-1]} \E[Z_i^n]$,
$d = - \sum_{i \in [j-1]} \xi^n_i$, and taking expectations,
we have
$$
\E\left[\min \Big\{ Z_j^n, \Big(k - \sum_{i \in [j-1]} Z_i^n\Big)_+\Big\}\right]
= \E\left[\min \Big\{ \E[Z_j^n], \Big(k - \sum_{i \in [j-1]} \E[Z_i^n]\Big)_+\Big\}\right]
\pm \sum_{i \in [j]} \E[ \abs{\xi^n_i}].
$$
For the random variables $Z_1^n,\ldots,Z_m^n$ we have
$\Var[Z_j^n] = \E[ (\xi^n_j)^2] = n f_j (1-f_j)$ for all $j \in [m]$,
and by the Cauchy-Schwartz inequality
\begin{align*}
\E\left[\min \Big\{ Z_j^n, \Big(k - \sum_{i \in [j-1]} Z_i^n\Big)_+\Big\}\right]&
= \E\left[\min \Big\{ \E[Z_j^n], \Big(k - \sum_{i \in [j-1]} \E[Z_i^n] \Big)_+\Big\} \right]\\& \quad \pm \sum_{i \in [j]} \sqrt{n f_i (1-f_i)},
\end{align*}
and the result follows by setting $M = \sum_{j \in [m]}\sqrt{f_j(1-f_j)}$.
\halmos\endproof\vspace{0.25cm}

\proof{Proof of Proposition \ref{prop:DR-bounds}.}
The first part follows immediately from Lemma \ref{lm:minmax-deviations}.
Next, for $(n,k) \in \T'$ and $j_{\id}$ be the index such that
$\Fbar(a_{j_{\id}}) \leq k/n < \Fbar(a_{j_{\id}+1})$,
the optimal solution \eqref{eq:DR-optimal-solution}
of the deterministic relaxation $DR(n,k)$ is given by
\begin{equation}\label{eq:DR-optimal-solution-j0}
s^*_{j} =
\begin{cases}
\E[Z_j^n]                           & \mbox{if } j \leq j_{\id}-1 \\
k - \sum_{i \in [j-1]}\E[Z_i^n]   & \mbox{if } j = j_{\id} \\
0                                   & \mbox{otherwise}.
\end{cases}
\end{equation}
To estimate the deterministic-relaxation gap when $(n,k) \in \T'$
it suffices to study the differences $s_j^*-\E[\frakS_j^n]=\E[Z_j^n-\frakS^n_j]$ for $j \leq j_{\id} - 1$
and $ s_{j_{\id}}^{*}-\E[\frakS_{j_{\id}}^n]=\E[k - \sum_{i \in [j_{\id}-1]}Z_i^n  - \frakS^n_{j_{\id}}]$.
Recalling the definition of $\frakS^n_j$ in \eqref{eq:Sj-definition} we have
$$
0 \leq Z_j^n - \frakS_j^n
\leq \begin{cases}
    0                           & \mbox{if } k-\sum_{i \in [j-1]}Z_i^n\geq Z_j^n \\
    \sum_{i \in [j]}Z_i^n-k     & \mbox{if }  k-\sum_{i \in [j-1]}Z_i^n<Z_j^n.
  \end{cases}
$$
In particular,
$$
0 \leq Z_j^n-\min\big\{ Z_j^n , ( k-\sum_{i\in [j-1]}Z_i^n )_+ \big\}
  \leq \Big( \sum_{i \in [j]} Z_i^n - k \Big)_+.
$$
Taking expectations on both sides, and recalling \eqref{eq:DR-optimal-solution-j0} for $j \leq j_{\id} - 1$,
we have
$$
0 \leq s^*_j - \E[\frakS^n_j] \leq  \E\Big[ \Big(  \sum_{i \in [j]} Z_i^n - k \Big)_+ \Big].
$$
For such $j \leq j_{\id} -1 $ the sum $\sum_{i \in [j]} Z_i^n$ is a Binomial
random variable with $n$ trials and success probability $\Fbar(a_{j+1})\leq \Fbar(a_{j_{\id}})$,
so because $(n, k) \in \T'$ and $ (\Fbar(a_{j_{\id}}) + \epsilon') n \leq k$,
we obtain from \eqref{eq:binomial-lemma-E(X-k)bound} that
\begin{equation}\label{eq:sj-estimate}
s^*_j = \E[\frakS^n_j] \pm \frac{1}{4 \epsilon'}, \quad \quad \text{for all } j \leq j_{\id} - 1.
\end{equation}

Similarly,
$$
\{ k - \sum_{i \in [j_{\id}-1]}Z_i^n\} - \frakS^n_{j_{\id}}
= \max \{ k - \sum_{i \in [j_{\id}]}Z_i^n,  \min\{ 0, k - \sum_{i \in [j_{\id}-1]}Z_i^n\}\},
$$
so we have the two inequalities
$$
\min\left\{ 0, k - \sum_{i \in [j_{\id}-1]}Z_i^n\right\}
\leq \left\{ k - \sum_{i \in [j_{\id}-1]}Z_i^n \right\} - \frakS^n_{j_{\id}}
\leq \max \left\{ k - \sum_{i \in [j_{\id}]}Z_i^n, 0 \right\}.
$$
Here, the two sums $\sum_{i \in [j_{\id}-1]}Z_i^n$ and $\sum_{i \in [j_{\id}]}Z_i^n$
are again Binomial random variables with $n$ trials and
success probabilities given, respectively, by $\Fbar(a_{j_{\id}})$ and $\Fbar( a_{j_{\id}+1})$.
Taking expectations and using $\Fbar( a_{j_{\id}}) + \epsilon' \leq k / n \leq \Fbar( a_{j_{\id}+1}) - \epsilon'$,
Lemma \ref{lm:multinomial-lemma} guarantees that
$$
- \frac{1}{4 \epsilon'} \leq k - \sum_{i \in [j_{\id}-1]}\E[Z_i^n]  - \E[\frakS^n_{j_{\id}}] \leq \frac{1}{4 \epsilon'}.
$$
In turn, the representation \eqref{eq:DR-optimal-solution-j0} for $s^*_{j_{\id}}$
implies that
\begin{equation}\label{eq:sj0-estimate}
s^*_{j_{\id}} = \E[\frakS^n_{j_{\id}}]  \pm \frac{1}{4 \epsilon'}.
\end{equation}
Combining the estimates \eqref{eq:sj-estimate} and \eqref{eq:sj0-estimate}
then give us that
$$
DR(n,k) - \Voff^*(n,k) \leq \frac{a_1 m}{4 \epsilon'}
\quad \quad \text{for all } (n,k) \in \T',
$$
just as needed to complete the proof of the proposition.
\halmos\endproof\vspace{0.25cm}

\end{APPENDICES}

\end{document}